\documentclass[12pt]{article}
\usepackage[margin=1in]{geometry}
\usepackage{amsmath,amssymb,amsthm,mathtools,mathrsfs}
\usepackage{times}
\usepackage{xcolor}
\usepackage{booktabs,tabularx}
\usepackage{needspace}
\usepackage[colorlinks,pdfstartview=FitV,linkcolor=blue,anchorcolor=blue,citecolor=blue,urlcolor=blue]{hyperref}
\hypersetup{pdftitle={The inviscid limit near a self-similar three-vortex collapse},pdfauthor={Trinh T. Nguyen}}
\numberwithin{equation}{section}
\allowdisplaybreaks
\newtheorem{theorem}{Theorem}[section]
\newtheorem{proposition}[theorem]{Proposition}
\newtheorem{lemma}[theorem]{Lemma}
\newtheorem{corollary}[theorem]{Corollary}

\newtheorem{definition}[theorem]{Definition}
\theoremstyle{definition}

\newcommand{\R}{\mathbb R}
\newcommand{\Lop}{\mathcal L}
\newcommand{\app}{\mathrm{app}}
\newcommand{\dd}{\,\mathrm d}
\newcommand{\norm}[1]{\left\|#1\right\|}
\newcommand{\av}[1]{\langle #1\rangle}
\newcommand{\lsim}[1]{\mathrel{\underset{#1}{\lesssim}}}
\title{The inviscid limit of 2D Navier-Stokes near a self-similar vortex collapse}
\author{Trinh T. Nguyen\thanks{School of Mathematical Sciences,
Shanghai Jiao Tong University, Shanghai, China.\newline
Email: \texttt{q0004349@sjtu.edu.cn}.}}
\date{}
\begin{document}
\maketitle
\begin{abstract}
We study the inviscid limit of the two-dimensional Navier--Stokes equations with initial vorticity given by three point vortices whose inviscid trajectories collapse self-similarly at $t=1$. For every fixed $0<\gamma<1$, we prove that the solution is approximated by a superposition of Lamb--Oseen vortices centered at the point-vortex positions, with an $L^1$ error of order $O(\nu^{1-\gamma})$, uniformly for $0<t\le T_\nu=1-\nu^\gamma$. We use coordinates adapted to the self-similar collapse and study a long-time inviscid limit problem on the interval $0<s\le\gamma\log(1/\nu)$, with a time-dependent rescaled diffusion coefficient $\eta(s)$ that tends to zero as $s\to\infty$.
\end{abstract}
\tableofcontents
\section{Introduction}\label{sec:introduction}
The Helmholtz--Kirchhoff system \cite{He58,Ki76} describes the motion
of $N$ point vortices in a two-dimensional incompressible ideal fluid.
It is given by
\begin{equation}\label{eq:point-vortex-system}
 \dot z_k(t)=\frac1{2\pi}\sum_{\substack{j=1\\j\ne k}}^N\Gamma_j
       \frac{(z_k(t)-z_j(t))^\perp}{|z_k(t)-z_j(t)|^2},
 \qquad z_k(0)=z_k^0,\qquad k=1,\ldots,N,
\end{equation}
where $z_k(t)\in\R^2$ are the vortex positions, $\Gamma_k$ are their
circulations, and $x^\perp=(-x_2,x_1)$. Vortex dynamics has attracted great interest among
mathematicians and physicists
(see \cite{Li41,MP83,MZM88,MP93,MP94,Ne01,LV02,GW02,MLL05,GW05,Ga11,GM11,Ga18,NN19,LWZ20,KN25,DoGa26,DG26,Nguyen26,DoGh26}).
Formally, one may regard
$\omega_{\mathrm{pv}}(t)=\sum_{k=1}^N\Gamma_k\delta_{z_k(t)}$
as a solution of the two-dimensional Euler equation, disregarding the
self-interaction of each vortex. Such initial data are
not covered by the usual existence theory for weak Euler solutions:
the velocity is singular at the vortex positions, and the
self-interaction terms in the vorticity equation are undefined
\cite{Ga11}. For fixed viscosity $\nu>0$, the Navier--Stokes
Cauchy problem with initial vorticity $\omega_{\mathrm{pv}}(0)$ has
a unique global solution $\omega^\nu$, smooth for positive time
\cite{GG05}.

Gallay \cite{Ga11} proves that this viscous solution follows the
point-vortex motion on every fixed interval $[0,T]$ on which the
vortices remain separated. More precisely, if
$d_T:=\min_{0\le t\le T}\min_{j\ne k}|z_j(t)-z_k(t)|>0$,
then, for sufficiently small $\nu$, we have
\begin{equation}\label{eq:Gallay-fixed-interval}
 \left\|\omega^\nu(t)-\sum_{k=1}^N\frac{\Gamma_k}{4\pi\nu t}
       e^{-\frac{|\cdot-z_k(t)|^2}{4\nu t}}\right\|_{L^1(\R^2)}
 \le C\frac{\nu t}{d_T^2},\qquad 0<t\le T.
\end{equation}
A natural question is whether this approximation remains valid as the
vortices approach collision. At collision, the Helmholtz--Kirchhoff dynamics
become singular because the Biot--Savart kernel is singular at zero. The dependence
on the minimum separation in \cite{Ga11} prevents direct application
of the fixed-interval result in this regime. The inviscid limit
through a point-vortex collision remains a delicate open problem.
We refer to \cite{GS26} for a discussion of the
difficulty of describing the limit after collision.

In this paper, we study the inviscid limit of the Navier--Stokes
equations up to times arbitrarily close to, but before, a self-similar
collapse of the Helmholtz--Kirchhoff system. For three vortices, such a
collapse contracts and rotates
the vortex triangle while preserving its shape \cite{Ar10,Ki90,LKZ00,TT88}.
For simplicity, we restrict the analysis to $N=3$. The same method
also applies to self-similar collapses of finitely many vortices under
suitable spectral assumptions on the linearized center dynamics in
collapse coordinates.

\subsection{Setting and main results}
We consider the vorticity equation
\begin{equation}\label{eq:NS}
 \partial_t\omega^\nu+u^\nu\cdot\nabla\omega^\nu=\nu\Delta\omega^\nu,
 \qquad u^\nu=v^{\omega^\nu},\qquad
 \omega^\nu(0)=\sum_{k=1}^3\Gamma_k\delta_{z_k^0}.
\end{equation}
The Biot--Savart velocity and kernel are
\begin{equation}\label{eq:kernel}
 v^f(x)=\int_{\R^2}K(x-y)f(y)\dd y,\qquad
 K(x)=\frac{x^\perp}{2\pi|x|^2}.
\end{equation}
Throughout the paper, $\Gamma_k$ and $z_k^0$ are fixed.
The normalized Oseen profile and its velocity are
\begin{equation}\label{eq:Gaussian}
 \begin{gathered}
 G(\xi)=\frac1{4\pi}e^{-|\xi|^2/4},\qquad
 v^G(\xi)=\phi(|\xi|)\xi^\perp,\\
 \phi(r)=\frac{1-e^{-r^2/4}}{2\pi r^2},\qquad \phi(0)=\frac1{8\pi}.
 \end{gathered}
\end{equation}
Without loss of generality, we assume that the vortices collapse
at time $t=1$.
We denote counterclockwise rotation through $\theta$ by $\mathfrak R_\theta$.
\begin{definition}\label{def:collapse}
A normalized collapsing configuration consists of the circulations
$\Gamma_k$, distinct positions $z_k^0\in\R^2$, and $\beta\in\R$ such that
\begin{equation}\label{eq:configuration}
 \frac12z_k^0-\beta(z_k^0)^\perp+
 \sum_{j\ne k}\Gamma_jK(z_k^0-z_j^0)=0,\qquad k=1,2,3.
\end{equation}
\end{definition}
\begin{proposition}\label{prop:reference-solution}
For such a configuration, the trajectories
\begin{equation}\label{eq:reference-centers}
 z_k(t)=\sqrt{1-t}\,\mathfrak R_{\beta s(t)}z_k^0,
 \qquad s(t)=-\log(1-t),
\end{equation}
form the unique solution of \eqref{eq:point-vortex-system} on $[0,1)$.
In particular,
\begin{equation}\label{eq:reference-distance}
 |z_j(t)-z_k(t)|=\sqrt{1-t}|z_j^0-z_k^0|.
\end{equation}
\end{proposition}
We decompose $\omega^\nu=\sum_k\omega_k^\nu$, where
\begin{equation}\label{eq:components}
\partial_t\omega_k^\nu+u^\nu\cdot\nabla\omega_k^\nu
 =\nu\Delta\omega_k^\nu,\qquad
 \omega_k^\nu(0)=\Gamma_k\delta_{z_k^0}.
\end{equation}
Each component has the sign of $\Gamma_k$ and mass $\Gamma_k$.

Our first result concerns the vorticity approximation.

\begin{theorem}\label{thm:main}
Fix a normalized collapsing configuration. For every $0<\gamma<1$,
there exists $\nu_\gamma>0$ such that, for $0<\nu\le\nu_\gamma$
and $T_\nu=1-\nu^\gamma$,
\begin{equation}\label{eq:main-L1}
 \sup_{0<t\le T_\nu}
 \left\|\omega^\nu(t)-\sum_{k=1}^3\frac{\Gamma_k}{\nu t}
 G\!\left(\frac{\cdot-z_k(t)}{\sqrt{\nu t}}\right)\right\|_{L^1(\R^2)}
 \lsim{\gamma} \nu^{1-\gamma}.
\end{equation}
\end{theorem}

We also estimate the displacement of the vortex centers.

\begin{theorem}\label{thm:centers}
Under the assumptions of Theorem~\ref{thm:main}, after decreasing
$\nu_\gamma$ if necessary, for $0<\nu\le\nu_\gamma$ and $0<t\le T_\nu$,
\begin{equation}\label{eq:main-centers}
 \max_{1\le k\le3}
 \left|\frac1{\Gamma_k}\int_{\R^2}x\omega_k^\nu(t,x)\dd x-z_k(t)\right|
 \lsim{\gamma}\frac{\nu^2t^2}{(1-t)^{\frac{3}{2}}}.
\end{equation}
In particular,
\begin{equation}\label{eq:main-centers-uniform}
 \sup_{0<t\le T_\nu}\max_{1\le k\le3}
 \left|\frac1{\Gamma_k}\int_{\R^2}x\omega_k^\nu(t,x)\dd x-z_k(t)\right|
 \lsim{\gamma}\nu^{2-\frac{3}{2}\gamma}.
\end{equation}
\end{theorem}

\subsection{Idea of the proof}\label{sec:proof-roadmap}

\paragraph{Collapse coordinates.}
We use $s=-\log(1-t)$, $X=(1-t)^{-\frac{1}{2}}\mathfrak R_{-\beta s}x$,
and $\omega^\nu=(1-t)^{-1}W$. The reference configuration is stationary
at $z^0$, and $T_\nu$ becomes $S_\nu=\gamma\log(1/\nu)$.
We define $Z_k=\frac1{\Gamma_k}\int_{\R^2}XW_k\dd X$,
$\delta=\nu(e^s-1)$, and $\eta=\frac\nu\delta$, and write
\[
 W_k(s,X)=\frac{\Gamma_k}{\delta}
 w_k\!\left(s,\frac{X-Z_k}{\sqrt\delta}\right).
\]
The profiles $w_k$ have unit mass and zero first moments. Both $\delta(S_\nu)$ and
$\eta(S_\nu)$ tend to zero as the viscosity vanishes. In collapse coordinates,
the vortices remain separated while the rescaled diffusion weakens.
By Proposition~\ref{lem:reduction}, we obtain the vorticity approximation from
the estimates $\sum_k\|w_k-G\|_1\lesssim_\gamma\delta$ and
$|Z-z^0|\lesssim_\gamma\delta^2$ on $[\log2,S_\nu]$.

In Section~\ref{prof:section}, we construct the approximate solution
\[
 w_{\app}^k=G+\delta F_k^\nu+\delta^{\frac{3}{2}}H_k^\nu
                   +\delta^2P_k^\nu+\delta^{\frac{5}{2}}J_k^\nu.
\]
For the prescribed center velocity $\mathcal F(Z)+\delta^2V_2^\nu$,
the residual $\mathcal R_k$ satisfies
$\|\sqrt{p_k}\mathcal R_k\|_2\lesssim\delta(\delta+\nu)$.

\paragraph{Adapted energy.}
To adapt the energy estimate to the deformed vortex profiles, we use
the Arnold-energy method developed in 
\cite{GS24,DoGa26}. We briefly explain
the underlying cancellation in Section~\ref{wgt:Arnold-structure}.
We define the inner weight directly from
the approximate Euler profile, avoiding the separate nonsymmetric
correction to the weight used in \cite{Ga11,Nguyen26}. The global weight
$p_k=e^{Q_k}$ has three regimes: it is adapted to the deformed profile
near the vortex center, constant on a plateau containing the other
vortex centers, and exponentially increasing beyond the plateau.
The plateau eliminates transport contributions involving
$\nabla p_k$ near the other centers, while the exponential growth
controls the moments of the remainders.

\paragraph{Nonlinear transport.}
In collapse coordinates, the diffusion coefficient in the remainder
equation is $\eta=\nu/\delta=(1-t)/t$, which tends to zero at
$t=T_\nu=1-\nu^\gamma$ as $\nu\to0$. To control nonlinear transport
with this weak dissipation, we use weighted $L^q$ estimates for the
remainder velocities, with $q$ finite and sufficiently large depending
on $\gamma$. These estimates use only a small fraction of the available
dissipation and allow us to close the energy estimate for every fixed
$0<\gamma<1$; see Section~\ref{sec:nonlinear}.

\paragraph{Center estimates.}
In Section~\ref{ctr:section}, we combine Taylor expansion of the pair interactions
with moment cancellations and weighted tail estimates to obtain
$|Z'-\mathcal F(Z)-\delta^2V_2^\nu|
\lesssim\delta^3+\delta^2\sqrt{\mathcal E}+\delta^4\mathcal E$,
where $V_2^\nu$ is the leading correction coefficient defined in
\eqref{prof:center-coefficient}. 
We use a center expansion related to those in
\cite{NN19,Ga11,DoGa26}.

\paragraph{Notation.}
We write $\av x=(1+|x|^2)^{\frac{1}{2}}$ for a scalar or vector $x$.
In particular, we have $\av r=\av\xi$ when $r=|\xi|$.
Primes denote derivatives with respect to the function's argument.
We denote angular averaging by $\Pi_0$ and set $\Pi_{\ne0}=I-\Pi_0$.
In the rest of the paper, we fix $0<\gamma<1$ and
$0<\kappa<\frac{1}{6}$.

\section{Point-vortex collapse and center dynamics}
\label{sec:geometry}
The linear stability of self-similar three-vortex motions was studied
by Aref \cite{Ar10}. We use the following spectral and semigroup estimates.

\begin{proposition}\label{geo:spectrum}
Fix a normalized collapsing configuration as in
Definition~\ref{def:collapse}. For configurations with distinct centers,
define
\begin{equation}\label{eq:configuration-field}
 \mathcal F_k(Z)=\frac12Z_k-\beta Z_k^\perp+
              \sum_{j\ne k}\Gamma_jK(Z_k-Z_j),
 \qquad k=1,2,3,
\end{equation}
and its linearization at the equilibrium $z^0$ by
\begin{equation}\label{eq:center-linearization}
 \mathscr A=D\mathcal F(z^0):(\R^2)^3\longrightarrow(\R^2)^3.
\end{equation}
Then
\begin{equation}\label{geo:semigroup}
 \|e^{s\mathscr A}\|\lesssim(1+s)e^s,\qquad s\ge0,
\end{equation}
where the operator norm is induced by the Euclidean norm and the implicit
constant depends only on the fixed configuration. More precisely,
\begin{equation}\label{geo:characteristic-polynomial}
 \det(\lambda I_6-\mathscr A)
 =\bigl((\lambda-\tfrac12)^2+\beta^2\bigr)
                         \lambda^2(\lambda-1)^2.
\end{equation}
\end{proposition}

\subsection{Reference motion and linearization}
\begin{proof}[Proof of Proposition~\ref{prop:reference-solution}]
We consider a solution $y_k(t)$ of the point-vortex system
on a collision-free interval with $t<1$ and define
$Y_k(s)=e^{s/2}\mathfrak R_{-\beta s}y_k(1-e^{-s})$.
Differentiating and using $t'(s)=e^{-s}$, we obtain
\begin{equation}\label{geo:autonomous-center-system}
 Y_k'=\tfrac12Y_k-\beta Y_k^\perp
       +\sum_{j\ne k}\Gamma_jK(Y_k-Y_j)=\mathcal F_k(Y),
\end{equation}
where we use $K(r\mathfrak R_\theta x)=r^{-1}\mathfrak R_\theta K(x)$.
By \eqref{eq:configuration}, we have $\mathcal F(z^0)=0$, so $Y(s)=z^0$ is stationary.
Transforming back, we obtain \eqref{eq:reference-centers}.
Since rotations preserve distances, we also obtain \eqref{eq:reference-distance}.
The centers remain distinct for $t<1$ and tend to zero as $t\to1$.
Since the point-vortex field is smooth away from collisions, we obtain
uniqueness on each $[0,T]$, $T<1$.
\end{proof}

For sufficiently small $h$, we define the nonlinear remainder by
\begin{equation}\label{geo:center-nonlinearity-definition}
 \mathcal N(h):=\mathcal F(z^0+h)-\mathscr Ah.
\end{equation}

\begin{lemma}
\label{geo:linearization-lemma}
The linearization \eqref{eq:center-linearization} acts on
$h=(h_1,h_2,h_3)\in(\R^2)^3$ by
\begin{equation}\label{geo:A-action}
 (\mathscr Ah)_k
 =\frac12h_k-\beta h_k^\perp+
       \sum_{j\ne k}\Gamma_jDK(z_k^0-z_j^0)(h_k-h_j),
 \qquad k=1,2,3.
\end{equation}
For $x\ne0$ and $h\in\R^2$, the derivative of the Biot--Savart
kernel is
\begin{equation}\label{geo:DK-action}
 DK(x)h
 =\frac{h^\perp}{2\pi|x|^2}
  -\frac{(x\cdot h)x^\perp}{\pi|x|^4}
 =\frac{|x|^2h^\perp-2(x\cdot h)x^\perp}{2\pi|x|^4}.
\end{equation}
Equivalently,
\begin{equation}\label{geo:kernel-jacobian}
 DK(x)=\frac1{2\pi|x|^4}
 \begin{pmatrix}
  2x_1x_2&x_2^2-x_1^2\\
  x_2^2-x_1^2&-2x_1x_2
 \end{pmatrix}.
\end{equation}
Moreover, the remainder \eqref{geo:center-nonlinearity-definition}
satisfies
\begin{equation}\label{geo:nonlinear-linearization}
 \mathcal F(z^0+h)=\mathscr Ah+\mathcal N(h),
 \qquad |\mathcal N(h)|\lesssim |h|^2
\end{equation}
for $|h|$ sufficiently small.
\end{lemma}

\begin{proof}
Differentiating \eqref{eq:configuration-field}, we obtain
\eqref{geo:A-action}. By differentiating $K$, we also obtain
\eqref{geo:DK-action}--\eqref{geo:kernel-jacobian}.
For $|h|\le\frac14\min_{i\ne j}|z_i^0-z_j^0|$, all segments
$z_k^0-z_j^0+u(h_k-h_j)$, $0\le u\le1$, stay a fixed positive
distance from zero. Since the frame terms are linear, we apply
Taylor's formula to obtain
\begin{equation}\label{geo:explicit-center-remainder}
 \mathcal N_k(h)=\sum_{j\ne k}\Gamma_j
 \int_0^1(1-u)D^2K(z_k^0-z_j^0+u(h_k-h_j))
                 [h_k-h_j,h_k-h_j]\dd u.
\end{equation}
Since $D^2K$ is bounded along these segments, we obtain
$|\mathcal N(h)|\lesssim|h|^2$, and hence
\eqref{geo:nonlinear-linearization}.
\end{proof}

\subsection{Collapse invariants and symmetry directions}
\begin{lemma}\label{geo:invariants}
A collapsing configuration as above satisfies
\begin{equation}\label{geo:identities}
 \sum_{i<j}\Gamma_i\Gamma_j=0,
 \qquad \sum_i\Gamma_i z_i^0=0,
 \qquad \sum_i\Gamma_i|z_i^0|^2=0.
\end{equation}
In particular, $\sum_i\Gamma_i\ne0$, and both signs of circulation
occur.
\end{lemma}
\begin{proof}
Pairing the interaction terms and using $K(-x)=-K(x)$ and
$x\cdot K(x)=0$, we obtain
$\frac{\dd}{\dd t}\sum_i\Gamma_i z_i=0$ and
$\frac{\dd}{\dd t}\sum_i\Gamma_i|z_i|^2=0$.
Both conserved quantities tend to zero as $t\to1$, by
\eqref{eq:reference-centers}. They therefore vanish at $t=0$,
and we obtain the last two identities in \eqref{geo:identities}.

For the first identity, we consider the point-vortex Hamiltonian
\[
 H_{\rm pv}(Z)=-\frac1{2\pi}
 \sum_{i<j}\Gamma_i\Gamma_j\log|Z_i-Z_j|
\]
and use $\dot z_i=-\Gamma_i^{-1}(\nabla_{Z_i}H_{\rm pv}(z))^\perp$ to obtain
$\frac{\dd}{\dd t}H_{\rm pv}(z(t))=0$.
Substituting \eqref{eq:reference-distance}, we obtain
\[
 H_{\rm pv}(z(t))=H_{\rm pv}(z^0)
 -\frac{\log(1-t)}{4\pi}\sum_{i<j}\Gamma_i\Gamma_j.
\]
Since the Hamiltonian is constant, we conclude that $\sum_{i<j}\Gamma_i\Gamma_j=0$. Hence we have
$(\sum_i\Gamma_i)^2=\sum_i\Gamma_i^2>0$.
If all circulations had the same sign, the first identity would be
impossible.
\end{proof}

We set $\Gamma_{\rm tot}=\sum_{k=1}^3\Gamma_k$ and define the
translation and centered subspaces by
\begin{equation}\label{geo:translation-centered-spaces}
 \mathcal T=\{(u,u,u):u\in\R^2\},\qquad
 \mathcal V=\left\{h\in(\R^2)^3:
                         \sum_{k=1}^3\Gamma_kh_k=0\right\}.
\end{equation}

\begin{lemma}\label{geo:translation-splitting}
The decomposition $(\R^2)^3=\mathcal T\oplus\mathcal V$ is
invariant under $\mathscr A$; the two subspaces have dimensions two and four,
respectively.
On the translation subspace,
\begin{equation}\label{geo:translation-matrix}
 \mathscr A|_{\mathcal T}\simeq
 \begin{pmatrix}1/2&\beta\\-\beta&1/2\end{pmatrix},
 \qquad
 \|e^{s\mathscr A}|_{\mathcal T}\|=e^{s/2}.
\end{equation}
\end{lemma}
\begin{proof}
Since $\Gamma_{\rm tot}\ne0$, the projection onto $\mathcal T$ along
$\mathcal V$ is $h\mapsto(\bar h,\bar h,\bar h)$, where
$\bar h=\Gamma_{\rm tot}^{-1}\sum_k\Gamma_kh_k$.
Using $DK(-x)=DK(x)$ in \eqref{geo:A-action}, we obtain
\begin{align*}
 \mathscr A(u,u,u)&=(u/2-\beta u^\perp)_{k=1}^3,\\
 \sum_k\Gamma_k(\mathscr Ah)_k
 &=\frac12\sum_k\Gamma_kh_k
       -\beta\left(\sum_k\Gamma_kh_k\right)^\perp.
\end{align*}
Thus both subspaces are invariant, and the translation flow is
$u\mapsto e^{s/2}\mathfrak R_{-\beta s}u$, from which we obtain the stated norm.
\end{proof}

\begin{lemma}
\label{geo:symmetry-directions}
The vectors $r_0=(z^0)^\perp$ and $r_1=z^0-2\beta(z^0)^\perp$
are linearly independent. They satisfy $r_0,r_1\in\mathcal V$ and
\begin{equation}\label{geo:symmetry-eigenvectors}
 \mathscr Ar_0=0,\qquad \mathscr Ar_1=r_1.
\end{equation}
\end{lemma}
\begin{proof}
Using rotation equivariance and homogeneity, we obtain
\[
 \mathcal F(\mathfrak R_\theta z^0)=0,\qquad
 \mathcal F(rz^0)=\left(r-r^{-1}\right)
       \left(\frac12z^0-\beta(z^0)^\perp\right).
\]
Differentiating at $\theta=0$ and $r=1$, we obtain
$\mathscr Ar_0=0$ and $\mathscr Az^0=r_1$, and hence \eqref{geo:symmetry-eigenvectors}.
Using the identity $\sum_k\Gamma_kz_k^0=0$, we obtain
$r_0,r_1\in\mathcal V$. Since $z^0\ne0$ and
$z^0\cdot(z^0)^\perp=0$, we also see that they are independent.
\end{proof}

\subsection{The linearized spectrum}
\begin{proof}[Proof of Proposition~\ref{geo:spectrum}]
On the centered subspace, we use the circulation-weighted pairing
$\Omega_\Gamma(u,v)=\sum_{k=1}^3\Gamma_k u_k^\perp\cdot v_k$.
This form is skew-symmetric and nondegenerate on $(\R^2)^3$.
Since $\mathcal T$ and $\mathcal V$ are orthogonal for this pairing,
its restriction to $\mathcal V$ is also nondegenerate.
Using the second-moment identity in \eqref{geo:identities}, we obtain
\begin{equation}\label{geo:symmetry-pairings}
 \Omega_\Gamma(r_0,r_1)
       =-\sum_k\Gamma_k|z_k^0|^2=0.
\end{equation}

Using \eqref{geo:A-action}, $DK(-x)=DK(x)$, and
$\operatorname{tr}DK(x)=0$, we obtain by a direct calculation
\begin{equation}\label{geo:bilinear-identity}
 \Omega_\Gamma(\mathscr Au,v)+\Omega_\Gamma(u,\mathscr Av)
       =\Omega_\Gamma(u,v).
\end{equation}
By the nondegeneracy of the pairing on $\mathcal V$ and the independence of $r_0,r_1$, we see that
$\Omega_\Gamma(r_0,\cdot),\Omega_\Gamma(r_1,\cdot)$ are independent.
We choose $q_0,q_1\in\mathcal V$ such that
$\Omega_\Gamma(r_i,q_j)=\delta_{ij}$ for $i,j\in\{0,1\}$.
By \eqref{geo:symmetry-pairings}, we see that $(r_0,r_1,q_0,q_1)$ is a basis, and the coefficient of $q_i$ is
$\Omega_\Gamma(r_i,h)$ for $h\in\mathcal V$.

Using \eqref{geo:bilinear-identity} and
\eqref{geo:symmetry-eigenvectors}, we therefore obtain
\begin{equation}\label{geo:centered-block-matrix}
 \mathscr A|_{\mathcal V}\simeq
 \begin{pmatrix}
 0&0&*&*\\
 0&1&*&*\\
 0&0&1&0\\
 0&0&0&0
 \end{pmatrix}.
\end{equation}
Thus the centered characteristic polynomial is
$\lambda^2(\lambda-1)^2$. Multiplying by the translation factor
$(\lambda-1/2)^2+\beta^2$, we obtain
\eqref{geo:characteristic-polynomial}.

Both centered eigenvalues have algebraic multiplicity two,
so the Jordan blocks have size at most two and we have
$\|e^{s\mathscr A}|_{\mathcal V}\|\lesssim(1+s)e^s$.
The translation evolution has norm $e^{s/2}$ by
\eqref{geo:translation-matrix}. Applying the fixed change of basis, we obtain
\eqref{geo:semigroup}.
\end{proof}

\section{Collapse coordinates and reduction}
\label{sec:coordinates}
\subsection{Rotating collapse coordinates}
We use $s(t)=-\log(1-t)$ and $t(s)=1-e^{-s}$ for $0\le t<1$ and $s\ge0$.
For fixed $\nu>0$, we set
\begin{equation}\label{eq:delta-definition}
 \delta(s):=\nu(e^s-1)
       =\frac{\nu t(s)}{1-t(s)},\qquad s\ge0,
\end{equation}
and
\begin{equation}\label{eq:eta-definition}
 \eta(s):=\frac{\nu}{\delta(s)}
       =\frac1{e^s-1}
       =\frac{1-t(s)}{t(s)},\qquad s>0.
\end{equation}
The parameter $\delta$ is the squared ratio of the diffusion length
$\sqrt{\nu t}$ to the collapse length $\sqrt{1-t}$, and $\eta$ is
the diffusion coefficient in \eqref{eq:inner-equation}. They satisfy
\begin{equation}\label{eq:scale-identities}
 \delta'=\delta+\nu=(1+\eta)\delta,\qquad
 \eta'=-\eta(1+\eta),\qquad \eta\delta=\nu
 \quad(s>0).
\end{equation}
The long-time argument starts at
\begin{equation}\label{eq:scale-startup-values}
 t_0=\frac12,\qquad s_0=\log2,\qquad
 \delta(s_0)=\nu,\qquad \eta(s_0)=1.
\end{equation}
Thus we have $\nu\le\delta$ and $0<\eta\le1$ for $s\ge s_0$.
The endpoint $T_\nu=1-\nu^\gamma$
corresponds to $S_\nu=\gamma\log(1/\nu)$, with
$\delta(S_\nu)=\nu^{1-\gamma}-\nu$ and
$\eta(S_\nu)=\frac{\nu^\gamma}{1-\nu^\gamma}$.
We conclude that $\delta(S_\nu)\to0$ as $\nu\to0$.
For estimates at fixed time, we may treat $\delta,\nu$ independently
before setting $\delta=\delta(s)$; the equations hold for $s>0$.

\begin{lemma}\label{lem:rotating-equations}
For $0<t<1$, define
\begin{equation}\label{eq:global-change}
s=-\log(1-t),\qquad X=\frac1{\sqrt{1-t}}\mathfrak R_{-\beta s}x,\qquad
    \omega^\nu(t,x)=\frac1{1-t}W(s,X).
\end{equation}
Let $W_k$ denote the transform of $\omega_k^\nu$ under the same change of variables,
so that $W=\sum_kW_k$. Then \eqref{eq:NS} becomes
\begin{equation}\label{eq:global-equation}
 \partial_sW+(V+X/2-\beta X^\perp)\cdot\nabla_XW+W
 =\nu\Delta_XW,\qquad V=v^W.
\end{equation}
Equivalently,
\begin{equation}\label{eq:global-conservative}
 \partial_sW+\nabla_X\cdot[(V+X/2-\beta X^\perp)W]
 =\nu\Delta_XW.
\end{equation}
Each $W_k$ satisfies the same linear equation with the total velocity
$V$. The physical velocity and the $L^1$ norm transform as
\begin{align}
 u^\nu(t,x)&=\frac1{\sqrt{1-t}}\mathfrak R_{\beta s}V(s,X),
 \label{eq:velocity-change}\\
 \norm{\omega^\nu(t)}_{L^1_x}&=\norm{W(s)}_{L^1_X}.
 \label{eq:L1-isometry}
\end{align}
The norm identity also holds for differences of two fields transformed
by \eqref{eq:global-change}.

The reference trajectories $z_k(t)$ become the fixed centers $z_k^0$.
\end{lemma}
\begin{proof}
By a direct calculation at fixed $x$, we obtain
\begin{equation}\label{eq:global-chain-rule}
 \partial_t\omega^\nu(t,x)
  =\frac1{(1-t)^2}
   \left\{\partial_sW+
       \left(\frac X2-\beta X^\perp\right)\cdot\nabla_XW+W\right\}.
\end{equation}
Using homogeneity and rotation covariance of $K$, we obtain
\eqref{eq:velocity-change}. By orthogonality of the rotation, we also have
$\Delta_x\omega^\nu=\frac1{(1-t)^2}\Delta_XW$ and
$u^\nu\cdot\nabla_x\omega^\nu
=\frac1{(1-t)^2}V\cdot\nabla_XW$.
Substituting into \eqref{eq:NS}, we obtain \eqref{eq:global-equation}.
Using the identities $\nabla_X\cdot V=\nabla_X\cdot X^\perp=0$ and
$\nabla_X\cdot(X/2)=1$, we obtain \eqref{eq:global-conservative}.
We use the same computation for \eqref{eq:components}.
Using the Jacobian $\dd x=(1-t)\dd X$, we obtain \eqref{eq:L1-isometry},
including for differences. Substituting \eqref{eq:reference-centers} into
\eqref{eq:global-change}, we see that $z_k(t)$ is mapped to $z_k^0$.
\end{proof}

\subsection{Centered profiles and their equations}
For $s>0$, the exact center in collapse coordinates is
$Z_k(s)=\Gamma_k^{-1}\int_{\R^2}XW_k(s,X)\dd X$.
We put $a_{kj}(s)=Z_k(s)-Z_j(s)$, so $a_{kk}=0$.
We fix the separation constant
\begin{equation}\label{eq:fixed-gap}
 d_*=\frac14\min_{j\ne k}|z_j^0-z_k^0|>0.
\end{equation}
The Oseen operators used below are
\begin{equation}\label{eq:linearized-operators}
 \Lop f=\Delta_\xi f+\tfrac12\xi\cdot\nabla_\xi f+f,
 \qquad
 \Lambda f=v^G\cdot\nabla_\xi f+v^f\cdot\nabla_\xi G.
\end{equation}
They satisfy $\Lop G=0$ and $v^G\cdot\nabla G=0$.
The Gaussian moments are $\int_{\R^2}G\dd\xi=1$,
$\int_{\R^2}\xi G\dd\xi=0$, and $\int_{\R^2}|\xi|^2G\dd\xi=4$.

Following \cite{CS24,DoGa26}, we define the centered physical profile by
\begin{equation}\label{eq:physical-profile}
 \widehat w_k^\nu(t,\xi)=\frac{\nu t}{\Gamma_k}
 \omega_k^\nu\!\left(t,\sqrt{1-t}\,\mathfrak R_{\beta s}Z_k(s)
                         +\sqrt{\nu t}\,\xi\right).
\end{equation}

\begin{lemma}\label{lem:vortex-equations}
Define the profiles by
\begin{equation}\label{eq:inner-change}
 \xi=\frac{X-Z_k(s)}{\delta^{\frac12}},\qquad
 W_k(s,X)=\frac{\Gamma_k}{\delta}w_k(s,\xi).
\end{equation}
Then
\begin{align}\label{eq:inner-equation}
 \partial_sw_k&+\frac{\Gamma_k}{\delta}v^{w_k}\cdot\nabla_\xi w_k
 +\sum_{j\ne k}\frac{\Gamma_j}{\delta}
 v^{w_j}(\xi+\frac{a_{kj}}{\delta^{\frac12}})\cdot\nabla_\xi w_k\notag\\
 &+\frac{Z_k/2-\beta Z_k^\perp-Z_k'}{\delta^{\frac12}}\cdot\nabla_\xi w_k
 -\beta \xi^\perp\cdot\nabla_\xi w_k=\eta\Lop w_k.
\end{align}
Moreover, $w_k\ge0$, $\int_{\R^2}w_k\dd\xi=1$, and
$\int_{\R^2}\xi w_k\dd\xi=0$. The physical center satisfies
\begin{equation}\label{eq:physical-centroids}
 \frac1{\Gamma_k}\int_{\R^2}x\omega_k^\nu(t,x)\dd x
 =\sqrt{1-t}\,\mathfrak R_{\beta s}Z_k(s).
\end{equation}
The physical and centered coordinates are related by
\begin{equation}\label{eq:physical-inner-relations}
x=\sqrt{1-t}\,\mathfrak R_{\beta s} \bigl(Z_k(s)+\delta^{\frac12}\,\xi\bigr),\qquad w_k(s,\xi)=\widehat
    w_k^\nu(t,\mathfrak R_{\beta s}\xi).
\end{equation}
\end{lemma}
\begin{proof}
By homogeneity of $K$, we have
\[
 V(s,Z_k+\delta^{\frac12}\,\xi)
   =\sum_j\frac{\Gamma_j}{\delta^{\frac12}}v^{w_j}(\xi+\frac{a_{kj}}{\delta^{\frac12}}).
\]
Substituting \eqref{eq:inner-change} into
\eqref{eq:global-equation} and using $\delta'/\delta=1+\eta$,
we obtain \eqref{eq:inner-equation} by a direct calculation.

Using the Jacobians $\dd x=(1-t)\dd X$ and $\dd X=\delta\dd\xi$,
we obtain \eqref{eq:physical-centroids} and
\[
 \begin{aligned}
 \int_{\R^2}w_k\dd\xi
 &=\frac1{\Gamma_k}\int_{\R^2}W_k\dd X=1,\\
 \int_{\R^2}\xi w_k\dd\xi
 &=\frac1{\Gamma_k\delta^{\frac12}}
       \int_{\R^2}(X-Z_k)W_k\dd X=0.
 \end{aligned}
\]
We obtain nonnegativity from $\omega_k^\nu/\Gamma_k\ge0$.
Finally, we derive \eqref{eq:physical-inner-relations} using
$\delta=\nu t/(1-t)$ and the inverse change of variables.
\end{proof}

The other centers have coordinates $-a_{kj}/\sqrt\delta$. If
$\max_k|Z_k(s)|\le C$ and $\min_{j\ne k}|Z_k(s)-Z_j(s)|\ge c>0$,
with $C,c$ independent of $s,\nu$, their distances
from the origin are of order $\delta^{-\frac{1}{2}}$, tending to
infinity uniformly on $s\le S_\nu$ as $\nu\to0$.

\subsection{Reduction to profile and center estimates}

\begin{proposition}\label{lem:reduction}
Take $\nu$ sufficiently small that $S_\nu>s_0$.
If the estimates
\begin{equation}\label{eq:rescaled-targets}
 \sum_{k=1}^3\|w_k(s)-G\|_{L^1(\R^2)}\lsim{\gamma}\delta(s),\qquad
 \max_k|Z_k(s)-z_k^0|\lsim{\gamma}\delta(s)^2,
\end{equation}
hold for $s_0\le s\le S_\nu$, then Theorem~\ref{thm:main} follows.
\end{proposition}
\begin{proof}
The intervals $[1/2,1-\nu^\gamma]$ and $[s_0,S_\nu]$
correspond under $t=1-e^{-s}$, and we have
$\delta(s)\le\nu^{1-\gamma}\le1$ for small $\nu$.
For any $h\in\R^2$, we obtain by a direct calculation
\begin{equation}\label{eq:Gaussian-translation-calculation}
 \|G(\cdot-h)-G\|_{L^1}
 \le |h|\|\nabla G\|_{L^1}\lesssim |h|.
\end{equation}

Using the radial symmetry of $G$ and the
$L^1$ isometry \eqref{eq:L1-isometry}, we obtain
\[
 \left\|\omega^\nu(t)-\sum_{k=1}^3\frac{\Gamma_k}{\nu t} G\!\left(\frac{\cdot-z_k(t)}{\sqrt{\nu t}}\right)\right\|_{L^1(\R^2)} =\left\|W(s)-\sum_{k=1}^3\frac{\Gamma_k}{\delta(s)} G\!\left(\frac{\cdot-z_k^0}{\sqrt{\delta(s)}}\right)\right\|_{L^1(\R^2)}.
\]
Combining \eqref{eq:inner-change} with
\eqref{eq:Gaussian-translation-calculation}, we obtain
\begin{align*}
 &\left\|W(s)-\sum_{k=1}^3\frac{\Gamma_k}{\delta}
       G\!\left(\frac{\cdot-z_k^0}{\delta^{\frac12}}\right)
   \right\|_{L^1(\R^2)}\\
 &\quad\le\sum_{k=1}^3|\Gamma_k|
   \left(\|w_k(s)-G\|_{L^1(\R^2)}
    +\|\nabla G\|_{L^1(\R^2)}
        \frac{|Z_k(s)-z_k^0|}{\delta^{\frac12}}\right)\\
 &\quad\lsim{\gamma} (\delta+\delta^{\frac{3}{2}})
 \lsim{\gamma} \nu^{1-\gamma},
\end{align*}
for $s_0\le s\le \gamma\log(\frac1\nu)$. Here we used \eqref{eq:rescaled-targets} and $\delta\le1$.

On the fixed interval $0<t\le1/2$, using Gallay's approximation and
regularized-center estimates \cite{Ga11}, we obtain
\[
 \left\|\omega^\nu(t)-\sum_{k=1}^3\frac{\Gamma_k}{\nu t}
       G\!\left(\frac{\cdot-z_k(t)}{\sqrt{\nu t}}\right)
   \right\|_{L^1(\R^2)}
 \lesssim \left(\nu+\frac{e^{-c/(\nu t)}}{\sqrt{\nu t}}\right)
 \lesssim \nu.
\]
The exponential term accounts for replacing the regularized centers
by $z_k(t)$ using the same Gaussian translation estimate. Since
$\nu\le\nu^{1-\gamma}$ for $0<\nu\le1$, we obtain
\eqref{eq:main-L1} by combining the two intervals.
\end{proof}

\section{Approximate solutions and the remainder equation}\label{prof:section}

We construct \eqref{prof:approx} by collecting powers of $\delta$ in
\eqref{eq:inner-equation}. For
$0<s\le S_\nu=\gamma\log(1/\nu)$, we have
\begin{equation}\label{prof:two-parameter-range}
 0<\delta(s)=\nu(e^s-1)\le\nu^{1-\gamma},\qquad
 \eta(s)=\nu/\delta(s).
\end{equation}
For the Taylor expansions, we assume
\begin{equation}\label{prof:shifted-velocity-region}
 0<\delta\le1,\qquad |\xi|\le\frac{d_*}{4\delta^{\frac12}}.
\end{equation}
Using \eqref{ctr:closed-force} and \eqref{eq:energy-growth},
proved below, we obtain
\begin{equation}\label{prof:configuration-neighborhood}
 |Z(s)-z^0|\le M\delta(s)^2,\qquad 0<s\le S_\nu,
\end{equation}
\begin{equation}\label{prof:center-velocity-bounds}
 |Z'(s)|+|Z'(s)-\mathcal F(Z(s))|\le C_Z\delta(s)^2,
 \qquad 0<s\le S_\nu,
\end{equation}
with $M,C_Z$ independent of $s$ and $\nu$ (see also \eqref{ctr:refined-center-equation}).
We work under these bounds, taking $\nu$ small after fixing $M,C_Z$.
With $d_*$ as in \eqref{eq:fixed-gap} and $M\nu^{2(1-\gamma)}\le d_*$, we obtain
\begin{equation}\label{prof:gap}
 \min_{j\ne k}|Z_k(s)-Z_j(s)|\ge d_*>0,\qquad 0<s\le S_\nu.
\end{equation}
We set
\begin{equation}\label{prof:approx}
 w_{\app}^k(\xi,s)=G+\delta F_k^\nu+\delta^{\frac32}H_k^\nu
                   +\delta^2P_k^\nu+\delta^{\frac52}J_k^\nu.
\end{equation}
For smooth profiles $q_j(\xi,s)$ and velocities
$v_j(\xi,s)\in\R^2$, we define
\begin{equation}\label{prof:construction-residual}
 \begin{split}
 \mathcal R_k[q;v](\xi,s):={}&(\partial_s-\eta\Lop)q_k
 +\sum_j\frac{\Gamma_j}{\delta}
       v_j\!\left(\xi+\frac{a_{kj}}{\delta^{\frac12}}\right)
       \cdot\nabla q_k\\
 &+\left\{\frac{Z_k/2-\beta Z_k^\perp-Z_k'}{\delta^{\frac12}}
                    -\beta\xi^\perp\right\}\cdot\nabla q_k.
 \end{split}
\end{equation}
In this notation, \eqref{eq:inner-equation} reads $\mathcal R_k[w;v^w]=0$.
Substituting $w_{\app}$, we obtain the error
$\mathcal R_k[w_{\app};v^{w_{\app}}]$.

\begin{theorem}\label{prof:construction-theorem}
Assume that $Z(s)$ is a smooth center path satisfying \eqref{prof:gap}
and that $\nu\ge0$ is sufficiently small. The profile equations
\eqref{prof:profile-system} admit smooth solutions
$F_k^\nu,H_k^\nu,P_k^\nu,J_k^\nu$ satisfying the Gaussian bounds
\eqref{prof:derivativebounds}, with zero mass and first moments.
Consequently, the approximation \eqref{prof:approx} has unit mass and zero first moments.
For $\nu>0$, under \eqref{prof:two-parameter-range} and
\eqref{prof:center-velocity-bounds}, the error evaluated with the prescribed
center velocity $\mathcal F(Z)+\delta^2V_2^\nu$, as defined in
\eqref{prof:prescribed-residual}, satisfies
$\mathfrak n_{0,\vartheta}(\mathcal R_k)
\lesssim\delta(\delta+\nu)$ for $0<\vartheta<1$.
\end{theorem}

\subsection{Error expansion}
The center-velocity correction enforces first-moment compatibility
for $J_k^\nu$. For a fixed path $Z(s)$ and smooth, rapidly decaying
profiles $q=(q_j)_j$, we obtain by substitution
\begin{equation}\label{prof:gaussian-perturbation-identity}
 \begin{split}
 &\mathcal R_k[G+q;v^G+v^q]\\
 &\quad=\left(\partial_s-\eta\Lop+\frac{\Gamma_k}{\delta}\Lambda\right)q_k\\
 &\qquad+\sum_{j\ne k}\frac{\Gamma_j}{\delta}
       \left[
       v^G\!\left(\xi+\frac{a_{kj}}{\delta^{\frac12}}\right)
                         \cdot\nabla q_k
       +v^{q_j}\!\left(\xi+\frac{a_{kj}}{\delta^{\frac12}}\right)
                         \cdot\nabla G\right]\\
 &\qquad+\left\{\frac{Z_k/2-\beta Z_k^\perp-Z_k'}{\delta^{\frac12}}
                    -\beta\xi^\perp\right\}\cdot\nabla q_k\\
 &\qquad+\sum_j\frac{\Gamma_j}{\delta}
       v^{q_j}\!\left(\xi+\frac{a_{kj}}{\delta^{\frac12}}\right)
                         \cdot\nabla q_k\\
 &\qquad+\left\{
       \sum_{j\ne k}\frac{\Gamma_j}{\delta}
       v^G\!\left(\xi+\frac{a_{kj}}{\delta^{\frac12}}\right)
       +\frac{Z_k/2-\beta Z_k^\perp-Z_k'}{\delta^{\frac12}}
       \right\}\cdot\nabla G.
 \end{split}
\end{equation}

We use the following standard Biot--Savart estimates; see \cite{DG26,Nguyen26}.
Here $D^nK(x)$ denotes the $n$th multilinear derivative, with $D^0K=K$.

\begin{lemma}\label{prof:kernel-taylor}
For $n\ge1$, $x\ne0$, and $h_1,\ldots,h_n\in\R^2$,
\begin{equation}\label{prof:kernel-derivatives}
 |D^nK(x)[h_1,\ldots,h_n]|
 \lsim{n}|x|^{-n-1}\prod_{j=1}^n|h_j|.
\end{equation}
For an integer $N\ge0$, $0<\theta<1$, $x\ne0$, and $|h|\le\theta|x|$,
\begin{equation}\label{prof:kernel-expansion}
 K(x+h)=\sum_{n=0}^N\frac1{n!}D^nK(x)[h^{\otimes n}]
                         +\mathcal R_N(x,h),
\end{equation}
where the remainder satisfies
\begin{equation}\label{prof:kernel-uniform-remainder}
 |\mathcal R_N(x,h)|
 \lsim{N,\theta}\frac{|h|^{N+1}}{|x|^{N+2}}.
\end{equation}
\end{lemma}

\begin{lemma}\label{prof:shifted-velocity}
Assume that $|a|\ge d_*>0$ and $q$ is smooth with
$|q(y)|\le A e^{-\vartheta|y|^2/4}$ for $y\in\R^2$,
where $A\ge0$ and $0<\vartheta<1$.
Under \eqref{prof:shifted-velocity-region},
\begin{equation}\label{prof:shifted-velocity-expansion}
 \begin{split}
 v^q\!\left(\xi+\frac{a}{\delta^{\frac12}}\right)
 ={}&\delta^{\frac12}K(a)\int_{\R^2}q(y)\dd y
       +\delta\int_{\R^2}DK(a)[\xi-y]q(y)\dd y\\
 &+\frac{\delta^{\frac32}}{2}\int_{\R^2}
       D^2K(a)[\xi-y,\xi-y]q(y)\dd y\\
 &+\frac{\delta^2}{6}\int_{\R^2}
       D^3K(a)[\xi-y,\xi-y,\xi-y]q(y)\dd y\\
 &+\frac{\delta^{\frac52}}{24}\int_{\R^2}
       D^4K(a)[\xi-y,\xi-y,\xi-y,\xi-y]q(y)\dd y\\
 &+O_{d_*,\vartheta}\!\left(A\delta^3\av\xi^5\right).
 \end{split}
\end{equation}
\end{lemma}
\begin{proof}
By homogeneity of $K$, we have the exact identity
\begin{equation}\label{prof:shifted-velocity-integral}
 v^q\!\left(\xi+\frac{a}{\delta^{\frac12}}\right)
 =\delta^{\frac12}\int_{\R^2}
       K\bigl(a+\delta^{\frac12}(\xi-y)\bigr)q(y)\dd y.
\end{equation}
We split the integral at $|y|=d_*/(4\delta^{\frac12})$.
For $|y|\le d_*/(4\delta^{\frac12})$, using \eqref{prof:shifted-velocity-region}, we obtain
$|a+u\delta^{\frac12}(\xi-y)|
\ge |a|-\delta^{\frac12}(|\xi|+|y|)\ge d_*/2$
for $0\le u\le1$.
Using \eqref{prof:kernel-expansion}--\eqref{prof:kernel-uniform-remainder}
to expand through degree four, we obtain the displayed terms with truncated
moments and bound the remainder by
$C_{d_*}\delta^3\int_{\R^2}(|\xi|+|y|)^5|q(y)|\dd y
\le C_{d_*,\vartheta}A\delta^3\av\xi^5$.
Replacing the truncated moments by the full moments, we incur an error bounded by
$C_{d_*,\vartheta}Ae^{-c/\delta}$.
Splitting $\{y:|y|>d_*/(4\delta^{\frac12})\}$ at
$|\xi+a/\delta^{\frac12}-y|=1$ and using local integrability of $K$, we obtain
\[
 \begin{aligned}
 &\delta^{\frac12}\left|\int_{|y|>d_*/(4\delta^{\frac12})} K\bigl(a+\delta^{\frac12}(\xi-y)\bigr)q(y)\dd y\right|\\
 &\quad\le\sup_{|y|>d_*/(4\delta^{\frac12})}|q(y)| +\frac1{2\pi}\int_{|y|>d_*/(4\delta^{\frac12})}|q(y)|\dd y \le\left(1+\frac2\vartheta\right) A\exp\!\left(-\frac{\vartheta d_*^2}{64\delta}\right).
 \end{aligned}
\]
Combining these estimates, we obtain the stated remainder bound.
\end{proof}

\begin{corollary}\label{prof:shifted-Gaussian-velocity}
Assume that $|a|\ge d_*>0$.
Under \eqref{prof:shifted-velocity-region},
\begin{equation}\label{prof:shifted-Gaussian-expansion}
 \begin{split}
 v^G\!\left(\xi+\frac{a}{\delta^{\frac12}}\right) ={}&\delta^{\frac12}K(a)+\delta DK(a)[\xi] +\frac{\delta^{\frac32}}2D^2K(a)[\xi,\xi]\\
 &+\frac{\delta^2}{6}D^3K(a)[\xi,\xi,\xi] +\frac{\delta^{\frac52}}{24}D^4K(a)[\xi,\xi,\xi,\xi] +O_{d_*}(\delta^3\av\xi^5).
 \end{split}
\end{equation}
\end{corollary}
\begin{proof}
We apply Lemma~\ref{prof:shifted-velocity} with $q=G$,
$A=1/(4\pi)$, and $\vartheta=1/2$.
For $1\le n\le4$, the polynomial
$y\mapsto D^nK(a)[\xi-y,\ldots,\xi-y]$ has harmonic components,
because $K$ is harmonic away from zero.
Using the mean-value property, radial symmetry of $G$, and $\int_{\R^2}G=1$,
we therefore obtain
\[
 \int_{\R^2}D^nK(a)[\xi-y,\ldots,\xi-y]G(y)\dd y
 =D^nK(a)[\xi,\ldots,\xi].
\]
Substituting this identity into \eqref{prof:shifted-velocity-expansion}, we obtain
\eqref{prof:shifted-Gaussian-expansion}.
\end{proof}

We define the strain matrix and the Taylor fields by
\begin{align}
 \Sigma_k(s)&=\sum_{j\ne k}\Gamma_jDK(a_{kj}(s)),\label{eq:strain-matrix}\\
 T_{k,n}(\xi,s)&=\frac1{(n-1)!}\sum_{j\ne k}\Gamma_j
 D^{n-1}K(a_{kj}(s))[\xi,\ldots,\xi],\qquad n=2,3,4,5.\label{prof:taylorfields}
\end{align}
\begin{lemma}\label{prof:linear-expansion}
Assume that $Z(s)$ is a $C^1$ path satisfying \eqref{prof:gap} and that the
corrections in \eqref{prof:approx} are smooth in $(\xi,s)$.
For every $k,s$ and each $Q_k\in\{F_k^\nu,H_k^\nu,P_k^\nu,J_k^\nu\}$,
assume
\begin{equation}\label{prof:linear-moment-conditions}
 \int_{\R^2}Q_k(y,s)\dd y=0,\qquad
 \int_{\R^2}yQ_k(y,s)\dd y=0,
\end{equation}
and
\begin{equation}\label{prof:linear-profile-bound}
 |Q_k(\xi,s)|+|\nabla Q_k(\xi,s)|
 \le C_0e^{-\vartheta|\xi|^2/4},\qquad \xi\in\R^2,
\end{equation}
for fixed $C_0>0$ and $0<\vartheta<1$. Set
$q_k(\xi,s)=\delta F_k^\nu+\delta^{\frac32}H_k^\nu
+\delta^2P_k^\nu+\delta^{\frac52}J_k^\nu$.
Under \eqref{prof:shifted-velocity-region}, the linear part of
\eqref{prof:gaussian-perturbation-identity} is
\begin{equation}\label{prof:linear-collected}
 \begin{aligned}
 &[\Gamma_k\Lambda+\nu(1-\Lop)]F_k^\nu +\delta^{\frac12}[\Gamma_k\Lambda+\nu(\tfrac32-\Lop)]H_k^\nu\\
 &+\delta\bigl\{[\Gamma_k\Lambda+\nu(2-\Lop)]P_k^\nu+F_k^\nu +(\Sigma_k\xi-\beta\xi^\perp)\cdot\nabla F_k^\nu\bigr\}\\
 &+\delta^{\frac32}\bigl\{[\Gamma_k\Lambda+\nu(\tfrac52-\Lop)]J_k^\nu +\tfrac32H_k^\nu\\
 &\hspace{24mm}+(\Sigma_k\xi-\beta\xi^\perp)\cdot\nabla H_k^\nu +T_{k,3}\cdot\nabla F_k^\nu +\mathcal U_k^\nu\cdot\nabla G\bigr\}\\
 &+\delta\partial_sF_k^\nu+\delta^{\frac32}\partial_sH_k^\nu +\delta^2\partial_sP_k^\nu+\delta^{\frac52}\partial_sJ_k^\nu +\frac{\mathcal F_k(Z)-Z_k'}{\delta^{\frac12}}\cdot\nabla q_k +O\!\left(\delta^2\av\xi^6e^{-\vartheta|\xi|^2/4}\right),
 \end{aligned}
\end{equation}
where $T_{k,n}$ is defined in \eqref{prof:taylorfields} and
\begin{equation}\label{prof:remote-center-coefficient}
 \mathcal U_k^\nu(s)=\frac12\sum_{j\ne k}\Gamma_j
       \int_{\R^2}D^2K(a_{kj}(s))[y,y]F_j^\nu(y,s)\dd y.
\end{equation}
The remainder constant depends on $C_0,\vartheta,d_*$, the fixed
circulations, and $\beta$, but is independent of $Z',s,\nu,\delta$.
\end{lemma}
\begin{proof}
Using \eqref{eq:scale-identities} and $\Lop G=0$,
we differentiate \eqref{prof:approx} to obtain
\begin{equation}\label{prof:scale-expansion}
 \begin{split}
 (\partial_s-\eta\Lop)w_{\app}^k ={}&\delta F_k^\nu+\tfrac32\delta^{\frac32}H_k^\nu +2\delta^2P_k^\nu+\tfrac52\delta^{\frac52}J_k^\nu\\
 &+\nu(1-\Lop)F_k^\nu +\nu\delta^{\frac12}(\tfrac32-\Lop)H_k^\nu +\nu\delta(2-\Lop)P_k^\nu +\nu\delta^{\frac32}(\tfrac52-\Lop)J_k^\nu\\
 &+\delta\partial_sF_k^\nu+\delta^{\frac32}\partial_sH_k^\nu +\delta^2\partial_sP_k^\nu+\delta^{\frac52}\partial_sJ_k^\nu.
 \end{split}
\end{equation}

Using Corollary~\ref{prof:shifted-Gaussian-velocity} and
\eqref{eq:configuration-field}, we obtain
\[
 \begin{aligned}
 &\sum_{j\ne k}\frac{\Gamma_j}{\delta}
       v^G\!\left(\xi+\frac{a_{kj}}{\delta^{\frac12}}\right)
       +\frac{Z_k/2-\beta Z_k^\perp-Z_k'}{\delta^{\frac12}}
       -\beta\xi^\perp\\
 &\quad=\Sigma_k\xi-\beta\xi^\perp+\delta^{\frac12}T_{k,3}
       +\frac{\mathcal F_k(Z)-Z_k'}{\delta^{\frac12}}
       +O(\delta\av\xi^5).
 \end{aligned}
\]
Using the moment conditions \eqref{prof:linear-moment-conditions},
we remove the constant and linear terms in
\eqref{prof:shifted-velocity-expansion}. From the quadratic term, we obtain
\begin{equation}\label{prof:linear-remote-expansion}
 \sum_{j\ne k}\frac{\Gamma_j}{\delta}
       v^{q_j}\!\left(\xi+\frac{a_{kj}}{\delta^{\frac12}}\right)
 =\delta^{\frac32}\mathcal U_k^\nu+O(\delta^2\av\xi^5).
\end{equation}
Collecting powers of $\delta$, including $\Gamma_k\Lambda q_k/\delta$,
we obtain \eqref{prof:linear-collected}. We bound the remainder using
\eqref{prof:linear-profile-bound},
$|\nabla q_k|\lesssim\delta e^{-\vartheta|\xi|^2/4}$, and
$|\nabla G|\lesssim\av\xi e^{-|\xi|^2/4}$.
\end{proof}

\begin{lemma}\label{prof:nonlinear-expansion}
Assume that the corrections in \eqref{prof:approx} are smooth and satisfy
\eqref{prof:linear-moment-conditions} and \eqref{prof:linear-profile-bound}
for every index and time. Assume also \eqref{prof:gap} and set
$q_k(\xi,s)=\delta F_k^\nu+\delta^{\frac32}H_k^\nu
+\delta^2P_k^\nu+\delta^{\frac52}J_k^\nu$.
Under \eqref{prof:shifted-velocity-region}, the nonlinear part of
\eqref{prof:gaussian-perturbation-identity} is
\begin{equation}\label{prof:nonlinear-collected}
 \begin{split}
 &\sum_j\frac{\Gamma_j}{\delta} v^{q_j}\!\left(\xi+\frac{a_{kj}}{\delta^{\frac12}}\right) \cdot\nabla q_k =\delta\Gamma_kv^{F_k^\nu}\cdot\nabla F_k^\nu +\delta^{\frac32}\Gamma_k \bigl\{v^{F_k^\nu}\cdot\nabla H_k^\nu +v^{H_k^\nu}\cdot\nabla F_k^\nu\bigr\}\\
 &\qquad+O\!\left(\delta^2\av\xi^5e^{-\vartheta|\xi|^2/4}\right).
 \end{split}
\end{equation}
The remainder constant depends only on $C_0,\vartheta,d_*$ and the
fixed circulations; it is uniform in $s,\nu,\delta$.
\end{lemma}
\begin{proof}
Expanding the self-interaction, we obtain the two displayed terms
in \eqref{prof:nonlinear-collected}. The remainder is
$O(\delta^2e^{-\vartheta|\xi|^2/4})$, by
\eqref{prof:Gaussian-velocity-bound} and \eqref{prof:linear-profile-bound}.
Using the gap and \eqref{prof:linear-profile-bound}, we obtain
$|\mathcal U_k^\nu|\lesssim1$ and
$|\nabla q_k|\lesssim\delta e^{-\vartheta|\xi|^2/4}$.
Using \eqref{prof:linear-remote-expansion}, we therefore obtain
\[
 \left|\sum_{j\ne k}\frac{\Gamma_j}{\delta} v^{q_j}\!\left(\xi+\frac{a_{kj}}{\delta^{\frac12}}\right) \cdot\nabla q_k\right| \lesssim\bigl(\delta^{\frac52}+\delta^3\av\xi^5\bigr) e^{-\vartheta|\xi|^2/4} \lesssim\delta^2\av\xi^5e^{-\vartheta|\xi|^2/4}.
\]
\end{proof}

\begin{proposition}\label{prof:residual-expansion}
Assume that $Z(s)$ is a $C^1$ path satisfying \eqref{prof:gap} and that the
smooth profiles in \eqref{prof:approx} satisfy
\eqref{prof:linear-moment-conditions} and \eqref{prof:linear-profile-bound}.
Under \eqref{prof:shifted-velocity-region}, the error
\eqref{prof:construction-residual} of $w_{\app}$ has the expansion
\begin{equation}\label{prof:collected-residual}
 \begin{aligned}
 &\mathcal R_k[w_{\app};v^{w_{\app}}] =\frac{\mathcal F_k(Z)-Z_k'}{\delta^{\frac12}} \cdot\nabla w_{\app}^k\\
 &\qquad+\bigl\{(\Sigma_k\xi)\cdot\nabla G +[\Gamma_k\Lambda+\nu(1-\Lop)]F_k^\nu\bigr\}\\
 &\qquad+\delta^{\frac12}\bigl\{T_{k,3}\cdot\nabla G +[\Gamma_k\Lambda+\nu(\tfrac32-\Lop)]H_k^\nu\bigr\}\\
 &\qquad+\delta\bigl\{C_k^\nu+[\Gamma_k\Lambda+\nu(2-\Lop)]P_k^\nu +\partial_sF_k^\nu\bigr\}\\
 &\qquad+\delta^{\frac32}\bigl\{\widetilde S_k^\nu +[\Gamma_k\Lambda+\nu(\tfrac52-\Lop)]J_k^\nu +\partial_sH_k^\nu\bigr\}\\
 &\qquad+\delta^2\partial_sP_k^\nu+\delta^{\frac52}\partial_sJ_k^\nu +O\!\left(\delta^2\av\xi^6e^{-\vartheta|\xi|^2/4}\right),
 \end{aligned}
\end{equation}
where
\begin{equation}\label{prof:Csource}
 C_k^\nu(\xi,s)=T_{k,4}\cdot\nabla G+F_k^\nu
 +(\Sigma_k\xi-\beta\xi^\perp)\cdot\nabla F_k^\nu
 +\Gamma_kv^{F_k^\nu}\cdot\nabla F_k^\nu,
\end{equation}
\begin{equation}\label{prof:uncentered-odd-source}
 \begin{split}
 \widetilde S_k^\nu(\xi,s)={}&(T_{k,5}+\mathcal U_k^\nu)\cdot\nabla G
       +\tfrac32H_k^\nu
 +(\Sigma_k\xi-\beta\xi^\perp)\cdot\nabla H_k^\nu
       +T_{k,3}\cdot\nabla F_k^\nu\\
 &+\Gamma_k\{v^{F_k^\nu}\cdot\nabla H_k^\nu
                         +v^{H_k^\nu}\cdot\nabla F_k^\nu\}.
 \end{split}
\end{equation}
Here $\Sigma_k,T_{k,n}$ and $\mathcal U_k^\nu$ are defined in
\eqref{eq:strain-matrix}, \eqref{prof:taylorfields}, and
\eqref{prof:remote-center-coefficient}, respectively.
The remainder constant depends only on $C_0,\vartheta,d_*$, the fixed
circulations, and $\beta$; it is uniform in $Z',s,\nu,\delta$.
\end{proposition}
\begin{proof}
Using \eqref{prof:shifted-Gaussian-expansion}, we obtain the Gaussian error in
\eqref{prof:gaussian-perturbation-identity} in the form
\[
 \mathcal R_k[G;v^G] ={}\frac{\mathcal F_k(Z)-Z_k'}{\delta^{\frac12}}\cdot\nabla G +\sum_{n=2}^5\delta^{\frac{n-2}{2}}T_{k,n}\cdot\nabla G +O\!\left(\delta^2\av\xi^6e^{-\vartheta|\xi|^2/4}\right).
\]
Adding \eqref{prof:linear-collected} and \eqref{prof:nonlinear-collected},
with $q_k=w_{\app}^k-G$ and $T_{k,2}=\Sigma_k\xi$, we obtain
\eqref{prof:collected-residual}.
\end{proof}

\subsection{Profile equations and the first three profiles}
The leading correction to the center velocity is
\begin{equation}\label{prof:center-coefficient}
 V_{2,k}^\nu(s)=\frac12\sum_{j\ne k}\Gamma_j\int_{\R^2}
 D^2K(a_{kj})[y,y]
       \bigl(F_k^\nu(y,s)+F_j^\nu(y,s)\bigr)\dd y.
\end{equation}
A zero superscript denotes the Euler profiles evaluated at $Z(s)$.
We obtain this coefficient by substituting \eqref{prof:approx} into
\eqref{ctr:exact-law}; see \eqref{ctr:corrected-approximate-force}.
By Proposition~\ref{ctr:energy-force}, we have
\begin{equation}\label{prof:center-law-definition}
 Z_k'=\mathcal F_k(Z)+\delta^2V_{2,k}^\nu+e_k.
\end{equation}
Using \eqref{ctr:refined-center-equation} and \eqref{eq:energy-growth},
we obtain $e_k=O_\gamma(\delta^3)$ for $s_0\le s\le S_\nu$. Thus we have
\[
 \frac{\mathcal F_k(Z)-Z_k'}{\delta^{\frac12}}
                         \cdot\nabla w_{\app}^k
 =-\delta^{\frac32}V_{2,k}^\nu\cdot\nabla G
   -\frac{e_k}{\delta^{\frac12}}\cdot\nabla w_{\app}^k
   +O\!\left(\delta^{\frac52}e^{-\vartheta|\xi|^2/4}\right).
\]
On this interval, the $e_k$ term is
$O_\gamma(\delta^{\frac52}e^{-\vartheta|\xi|^2/4})$.

We set
\begin{equation}\label{prof:receiving-center-coefficient}
 b_k^\nu(s)=\int_{\R^2}T_{k,3}(\xi,s)F_k^\nu(\xi,s)\dd\xi
          =V_{2,k}^\nu-\mathcal U_k^\nu,
\end{equation}
and define
\begin{equation}\label{prof:odd-source}
 \begin{split}
 S_k^\nu(\xi,s)={}&\widetilde S_k^\nu-V_{2,k}^\nu\cdot\nabla G\\
 ={}&(T_{k,5}-b_k^\nu)\cdot\nabla G+\tfrac32H_k^\nu
 +(\Sigma_k\xi-\beta\xi^\perp)\cdot\nabla H_k^\nu
 +T_{k,3}\cdot\nabla F_k^\nu\\
 &+\Gamma_k\{v^{F_k^\nu}\cdot\nabla H_k^\nu
                  +v^{H_k^\nu}\cdot\nabla F_k^\nu\}.
 \end{split}
\end{equation}

We choose the profiles to satisfy, for all $\xi\in\R^2$,
\begin{subequations}\label{prof:profile-system}
\begin{align}
 [\Gamma_k\Lambda+\nu(1-\Lop)]F_k^\nu
     &=-(\Sigma_k\xi)\cdot\nabla G,\label{prof:firstthree}\\
 [\Gamma_k\Lambda+\nu(\tfrac32-\Lop)]H_k^\nu
     &=-T_{k,3}\cdot\nabla G,\label{prof:third-profile}\\
 \Gamma_k\Lambda P_k^\nu&=-\Pi_{\ne0}C_k^\nu,\label{prof:remaining-profile}\\
 \Gamma_k\Lambda J_k^\nu
     &=-S_k^\nu.\label{prof:fifth-profile}
\end{align}
\end{subequations}
Under \eqref{prof:shifted-velocity-region}, we substitute the profile equations into
\eqref{prof:collected-residual} to obtain
\begin{equation}\label{prof:cancelled-residual}
 \begin{split}
 \mathcal R_k[w_{\app};v^{w_{\app}}]
 ={}&-\frac{e_k}{\delta^{\frac12}}\cdot\nabla w_{\app}^k
       +\delta\Pi_0C_k^\nu+\nu\delta(2-\Lop)P_k^\nu\\
 &+\nu\delta^{\frac32}(\tfrac52-\Lop)J_k^\nu
       +\delta\partial_sF_k^\nu+\delta^{\frac32}\partial_sH_k^\nu\\
 &+\delta^2\partial_sP_k^\nu+\delta^{\frac52}\partial_sJ_k^\nu
       +O\!\left(\delta^2\av\xi^6e^{-\vartheta|\xi|^2/4}\right).
 \end{split}
\end{equation}

Using \eqref{prof:derivativebounds} below, we bound the four time-derivative
terms in \eqref{prof:cancelled-residual} by $C\delta|Z'|$ in
$\mathfrak n_{0,\vartheta}$. Under \eqref{prof:center-velocity-bounds}, we therefore
obtain an $O(\delta^3)$ bound for their sum.

\paragraph{Function spaces and inverse estimates.}
We define $Y=L^2(\R^2,e^{|\xi|^2/4}\dd\xi;\R)$,
with $\|f\|_Y=\|e^{|\xi|^2/8}f\|_{L^2}$, and
$Y_n=\{a(r)\cos(n\theta)+b(r)\sin(n\theta)\}\cap Y$ for $n\ge1$.
We denote by $\Pi_n$ the orthogonal projection in $Y$ onto $Y_n$.
The radial sector is denoted by $Y_0$, and we define
$Y_1^{\mathrm c}=\{f\in Y_1:\ \int_{\R^2}\xi f\dd\xi=0\}$.
We  set
\begin{equation}\label{prof:Gaussian-derivative-norm}
 \mathfrak n_{m,\vartheta}(f)
 =\sum_{|\alpha|\le m}\sup_\xi
       e^{\vartheta|\xi|^2/4}|\partial^\alpha f(\xi)|,
 \qquad m\ge0,\quad 0<\vartheta<1.
\end{equation}
For angular modes $n\ge2$, we use the inverse estimates in
\cite{Nguyen26,GW07,Ga11,GS24,DoGa26}.

\begin{lemma}\label{prof:angular-input}
Assume that $n\ge2$, $\alpha\in\{1,3/2\}$, and $g\in Y_n$ is smooth,
with all derivatives of $e^{|\xi|^2/4}g$ of at most polynomial growth.
For sufficiently small $\nu\ge0$, there is a unique smooth solution
$q^\nu\in Y_n$ to
\begin{equation}\label{prof:regularized-inverse}
 [\Gamma_k\Lambda+\nu(\alpha-\Lop)]q^\nu=g.
\end{equation}
The same weighted regularity holds for $q^\nu$.
For $m=0,1,2$ and $0<\vartheta<\vartheta_1<1$,
\begin{equation}\label{prof:unregularized-bound}
 \mathfrak n_{m,\vartheta}(q^0)
 \lsim{n,\Gamma_k,m,\vartheta,\vartheta_1} \mathfrak n_{m,\vartheta_1}(g).
\end{equation}
For each fixed source $g$, $m=0,1,2,3$, and $0<\vartheta<1$,
\begin{equation}\label{prof:inverse-pointwise}
 \begin{split}
 \mathfrak n_{m,\vartheta}(q^\nu)
       &\lsim{n,\Gamma_k,g,m,\vartheta}1,\\
 \mathfrak n_{m,\vartheta}(q^\nu-q^0)
       &\lsim{n,\Gamma_k,g,m,\vartheta}\nu.
 \end{split}
\end{equation}
For every integer $m\ge0$, $0<\vartheta<1$, and
$q\in C^m(\R^2)$ with $\mathfrak n_{m,\vartheta}(q)<\infty$,
\begin{equation}\label{prof:Gaussian-velocity-bound}
 \sum_{|\alpha|\le m}\|\partial^\alpha v^q\|_{L^\infty}
       \lsim{m,\vartheta} \mathfrak n_{m,\vartheta}(q).
\end{equation}
For $0<\vartheta<1$ and $q\in Y_n$ with
$\mathfrak n_{0,\vartheta}(q)<\infty$,
\begin{equation}\label{prof:mode-velocity-input}
 |v^q(\xi)|\lsim{n,\vartheta} \av\xi^{-n-1}
                \mathfrak n_{0,\vartheta}(q),\qquad n\ge2.
\end{equation}
For every $0<\vartheta<1$ and every $q$ satisfying
$\mathfrak n_{0,\vartheta}(q)<\infty$,
$\int_{\R^2}q\dd\xi=0$, and $\int_{\R^2}\xi q\dd\xi=0$,
\begin{equation}\label{prof:centered-Gaussian-velocity}
 |v^q(\xi)|\lsim{\vartheta} \av\xi^{-3}\mathfrak n_{0,\vartheta}(q).
\end{equation}
\end{lemma}

\paragraph{Solving for $F_k^\nu$.}
\begin{lemma}\label{prof:F-profile}
Equation \eqref{prof:firstthree} has a unique solution
$F_k^\nu\in Y_2$, with zero mass and first moments.
For $0<\vartheta<1$,
\begin{equation}\label{prof:F-bounds}
 \begin{split}
 \mathfrak n_{3,\vartheta}(F_k^\nu)&\lesssim1,\\
 \mathfrak n_{1,\vartheta}(F_k^\nu-F_k^0)&\lesssim\nu,\\
 \mathfrak n_{1,\vartheta}(\partial_sF_k^\nu)&\lesssim |Z'(s)|.
 \end{split}
\end{equation}
\end{lemma}
\begin{proof}
Since $\Sigma_k$ is symmetric and trace free, we have
\[
 -(\Sigma_k\xi)\cdot\nabla G
 =\frac{r^2G}{2}\bigl\{(\Sigma_k)_{11}\cos2\theta
                         +(\Sigma_k)_{12}\sin2\theta\bigr\}.
\]
The forcing has frequency two, and we bound its
coefficients using the separation assumption. We then use the identity
\begin{equation}\label{prof:strain-time-dependence}
 \Sigma_k'=\sum_{j\ne k}\Gamma_j
       D^2K(a_{kj})[Z_k'-Z_j',\,\cdot\,]
\end{equation}
to bound their time derivatives by $C|Z'|$.
Applying Lemma~\ref{prof:angular-input}, we obtain the unique solution
$F_k^\nu\in Y_2$, together with the estimates \eqref{prof:F-bounds}.
\end{proof}

\paragraph{Solving for $H_k^\nu$.}
\begin{lemma}\label{prof:H-profile}
Equation \eqref{prof:third-profile} has a unique solution
$H_k^\nu\in Y_3$, with zero mass and first moments.
For $0<\vartheta<1$,
\begin{equation}\label{prof:H-bounds}
 \mathfrak n_{3,\vartheta}(H_k^\nu)\lesssim1,\qquad
 \mathfrak n_{1,\vartheta}(\partial_sH_k^\nu)\lesssim |Z'(s)|.
\end{equation}
\end{lemma}
\begin{proof}
With $a_{kj}=|a_{kj}|(\cos\phi_{kj},\sin\phi_{kj})$,
we expand $K$ and use the identity $\nabla G=-\xi G/2$ to obtain
\begin{equation}\label{prof:Gaussian-expansion}
 -T_{k,n}\cdot\nabla G
 =\frac{(-1)^{n-1}r^nG}{4\pi}
       \sum_{j\ne k}\frac{\Gamma_j}{|a_{kj}|^n}
                         \sin(n(\theta-\phi_{kj})),\qquad n=3,4,5.
\end{equation}
We have
$-T_{k,3}\cdot\nabla G\in\operatorname{span}\{r^3G\cos3\theta,r^3G\sin3\theta\}$. We bound
its coefficients using the separation assumption and use the identity
\begin{equation}\label{prof:T-time-dependence}
 \partial_sT_{k,n}
 =\frac1{(n-1)!}\sum_{j\ne k}\Gamma_j
       D^nK(a_{kj})[Z_k'-Z_j',\xi,\ldots,\xi]
\end{equation}
to bound their time derivatives by $C|Z'|$.
Applying Lemma~\ref{prof:angular-input} with $n=3$ and $\alpha=3/2$,
we obtain the unique solution $H_k^\nu\in Y_3$, together with
\eqref{prof:H-bounds}.
\end{proof}

\paragraph{Solving for $P_k^\nu$.}
\begin{lemma}\label{prof:P-profile}
The forcing in \eqref{prof:remaining-profile} satisfies
$-\Pi_{\ne0}C_k^\nu\in Y_2\oplus Y_4$. The profile equation has a
unique solution $P_k^\nu\in Y_2\oplus Y_4$, with zero mass and
first moments. For $0<\vartheta<1$,
\begin{equation}\label{prof:P-bounds}
 \mathfrak n_{2,\vartheta}(P_k^\nu)\lesssim1,\qquad
 \mathfrak n_{0,\vartheta}(\partial_sP_k^\nu)\lesssim |Z'(s)|.
\end{equation}
\end{lemma}
\begin{proof}
In \eqref{prof:Csource}, the first term has frequency four by
\eqref{prof:Gaussian-expansion}, the terms
$F_k^\nu-\beta\partial_\theta F_k^\nu$ have frequency two,
and the two remaining products have frequencies zero and four.
Thus we have $\Pi_{\ne0}C_k^\nu\in Y_2\oplus Y_4$.
For $0<\vartheta_1<\vartheta_2<1$, using \eqref{prof:Gaussian-velocity-bound}, we obtain
\begin{equation}\label{prof:transport-product-bound}
 \mathfrak n_{m,\vartheta_1}(v^f\cdot\nabla g)
 \lsim{m,\vartheta_2}
       \mathfrak n_{m,\vartheta_2}(f)
       \mathfrak n_{m+1,\vartheta_2}(g).
\end{equation}
Likewise, we have
\begin{equation}\label{prof:linear-transport-bound}
 \mathfrak n_{m,\vartheta_1}
       ((\Sigma_k\xi-\beta\xi^\perp)\cdot\nabla g)
 \lsim{m,\vartheta_1,\vartheta_2}
       \mathfrak n_{m+1,\vartheta_2}(g).
\end{equation}
Choosing $\vartheta<\vartheta_1<\vartheta_2<1$ and using
\eqref{prof:F-bounds}, \eqref{prof:strain-time-dependence}, and
\eqref{prof:T-time-dependence}, we obtain
$\mathfrak n_{2,\vartheta_1}(\Pi_{\ne0}C_k^\nu)\lesssim1$ and
$\mathfrak n_{0,\vartheta_1}
(\partial_s\Pi_{\ne0}C_k^\nu)\lesssim |Z'|$.
Applying Lemma~\ref{prof:angular-input} on modes two and four, we obtain existence
and uniqueness. We then apply \eqref{prof:unregularized-bound} with
$m=2$ to the equation and with $m=0$ to its $s$ derivative to obtain
\eqref{prof:P-bounds}.
\end{proof}

\paragraph{Solving for $J_k^\nu$.}
We seek $J_k^\nu$ with only angular
modes one, three, and five, and with zero first moments.
At $s=0$, we evaluate the coefficients using
$a_{kj}(0)=z_k^0-z_j^0$.

\begin{lemma}\label{prof:odd-solvability}
The source satisfies $S_k^\nu\in Y_1^{\mathrm c}\oplus Y_3\oplus Y_5$;
in particular, $\Pi_0S_k^\nu=0$ and
$\int_{\R^2}\xi S_k^\nu\dd\xi=0$.
\end{lemma}
\begin{proof}
Using the profile estimates and \eqref{prof:Gaussian-velocity-bound}, we find
that $S_k^\nu$ has the weighted regularity stated in
Lemma~\ref{prof:angular-input}. In polar components, $\Sigma_k\xi$ and
$v^{F_k^\nu}$ have frequency two, whereas $T_{k,3}$ and
$v^{H_k^\nu}$ have frequency three. Their products with derivatives
of $H_k^\nu$ and $F_k^\nu$, respectively, therefore have only
frequencies one and five. By \eqref{prof:Gaussian-expansion},
$T_{k,5}\cdot\nabla G$ has frequency five; the remaining terms
$-b_k^\nu\cdot\nabla G$ and
$\tfrac32H_k^\nu-\beta\partial_\theta H_k^\nu$ have frequencies
one and three. Thus we have $S_k^\nu\in Y_1\oplus Y_3\oplus Y_5$.

All transport fields in \eqref{prof:odd-source} are divergence free.
Integrating by parts, we obtain
$\int_{\R^2}\xi(B\cdot\nabla q)\dd\xi
=-\int_{\R^2}Bq\dd\xi$.
We also have
$\int_{\R^2}T_{k,5}G\dd\xi=0$ and
$\int_{\R^2}\xi H_k^\nu\dd\xi=0$ by angular orthogonality.
By oddness of $K$, we have
$\int_{\R^2}
\bigl(v^{F_k^\nu}H_k^\nu+v^{H_k^\nu}F_k^\nu\bigr)\dd\xi=0$.
Using these identities in \eqref{prof:odd-source}, we obtain
$\int_{\R^2}\xi S_k^\nu\dd\xi
=b_k^\nu-\int_{\R^2}T_{k,3}F_k^\nu\dd\xi=0$
by \eqref{prof:receiving-center-coefficient}. We have proved the claim.
\end{proof}
\begin{lemma}[{\cite{GS24,DoGa26}}]\label{prof:J-profile}
Equation \eqref{prof:fifth-profile} has a unique solution
$J_k^\nu\in Y_1^{\mathrm c}\oplus Y_3\oplus Y_5$.
For every $0<\vartheta<1$,
\begin{equation}\label{prof:J-bounds}
 \mathfrak n_{2,\vartheta}(J_k^\nu)\lesssim1,\qquad
 \mathfrak n_{0,\vartheta}(\partial_sJ_k^\nu)\lesssim |Z'(s)|.
\end{equation}
\end{lemma}
\begin{proof}[Proof sketch]
Using \eqref{prof:F-bounds}, \eqref{prof:H-bounds}, and the product bounds
\eqref{prof:transport-product-bound}--\eqref{prof:linear-transport-bound}
in \eqref{prof:odd-source}, we obtain, for every $0<\vartheta_1<1$,
\begin{equation}\label{prof:J-source-bounds}
 \mathfrak n_{2,\vartheta_1}(S_k^\nu)\lesssim1,\qquad
 \mathfrak n_{0,\vartheta_1}(\partial_sS_k^\nu)\lesssim |Z'(s)|.
\end{equation}
For the second bound, we also use \eqref{prof:strain-time-dependence},
\eqref{prof:T-time-dependence}, and
$b_k^\nu=\int_{\R^2}T_{k,3}F_k^\nu\dd\xi$.
By Lemma~\ref{prof:odd-solvability}, both sources have zero first moments.
Using the inverse of $\Lambda$ on $Y_1^{\mathrm c}\oplus Y_3\oplus Y_5$,
we obtain existence and uniqueness. On modes three and five, we obtain
the required bounds from \eqref{prof:unregularized-bound}.

For the centered solution of $\Lambda q=g$ on $Y_1^{\mathrm c}$,
each radial coefficient $b$ of $g$ satisfies
$\int_0^r t^2b(t)\dd t=-\int_r^\infty t^2b(t)\dd t$.
Differentiating the radial representation of the centered inverse in
Cartesian coordinates near the origin and using this tail identity
for $r\ge1$, we obtain
\begin{equation}\label{prof:J-centered-Euler-bound}
 \mathfrak n_{m,\vartheta}(q)
 \lsim{m,\vartheta,\vartheta_1}\mathfrak n_{m,\vartheta_1}(g),
 \qquad m\ge0,\quad \max\{\vartheta,\tfrac12\}<\vartheta_1<1.
\end{equation}
The centering condition removes the translation kernel.
Since the inverse is independent of $s$ and $\nu$, we have
$\Gamma_k\Lambda\partial_sJ_k^\nu=-\partial_sS_k^\nu$.
Applying \eqref{prof:J-centered-Euler-bound} with $m=2$ and $m=0$
to $-\Pi_1S_k^\nu/\Gamma_k$ and
$-\Pi_1\partial_sS_k^\nu/\Gamma_k$, respectively, and using
\eqref{prof:J-source-bounds}, we obtain \eqref{prof:J-bounds}
uniformly in $\nu$.
\end{proof}

\subsection{Gaussian bounds and radial cancellation}
\begin{proposition}\label{prof:localization}
Under \eqref{prof:gap}, for sufficiently small $\nu\ge0$, the profiles
satisfying \eqref{prof:profile-system} in the stated spaces obey,
for $0<\vartheta<1$,
\begin{equation}\label{prof:derivativebounds}
 \begin{split}
 \sum_{Q\in\{F_k^\nu,H_k^\nu,P_k^\nu,J_k^\nu\}}
             \mathfrak n_{2,\vartheta}(Q)&\lesssim1,\\
 \sum_{Q\in\{F_k^\nu,H_k^\nu,P_k^\nu,J_k^\nu\}}
             \mathfrak n_{0,\vartheta}(\partial_sQ)&\lesssim |Z'(s)|.
 \end{split}
\end{equation}
The coefficients $b^\nu,\mathcal U^\nu,V_2^\nu$ are uniformly
bounded, and $\partial_s b^\nu=O(|Z'|)$.

The associated velocities satisfy
\begin{equation}\label{prof:profile-velocities}
 \begin{gathered}
 |v^{F_k^\nu}(\xi)|+|v^{P_k^\nu}(\xi)|+|v^{J_k^\nu}(\xi)|
       \lesssim\av\xi^{-3},\\
 |v^{H_k^\nu}(\xi)|\lesssim\av\xi^{-4}.
 \end{gathered}
\end{equation}
\end{proposition}
\begin{proof}
Using Lemmas~\ref{prof:F-profile}, \ref{prof:H-profile},
\ref{prof:P-profile}, and \ref{prof:J-profile}, we obtain
\eqref{prof:derivativebounds}.
We obtain the coefficient bounds by differentiating
\eqref{prof:receiving-center-coefficient},
\eqref{prof:remote-center-coefficient}, and
\eqref{prof:center-coefficient}, using the gap,
\eqref{prof:kernel-derivatives}, and \eqref{prof:F-bounds}.

We apply \eqref{prof:mode-velocity-input}, recalled from \cite{Nguyen26}
(see also \cite{Ga11}),
to the angular modes of $F_k^\nu,H_k^\nu,P_k^\nu$, which are
$2$, $3$, and $\{2,4\}$, respectively. For $P_k^\nu$, we retain
the slower bound $\av\xi^{-3}$.
Applying \eqref{prof:centered-Gaussian-velocity} to $J_k^\nu$,
we obtain the same decay bound from its zero mass and first moments; see also \cite{DoGa26}.
Combining these bounds, we obtain \eqref{prof:profile-velocities}.
\end{proof}

By \eqref{prof:approx}, we have
$w_{\app}^k-G-\delta F_k^0
=\delta(F_k^\nu-F_k^0)+\delta^{\frac32}H_k^\nu
+\delta^2P_k^\nu+\delta^{\frac52}J_k^\nu$.
By \eqref{prof:derivativebounds}, \eqref{prof:F-bounds},
\eqref{prof:profile-velocities}, and \eqref{prof:centered-Gaussian-velocity},
we obtain, for fixed $\vartheta\in(1/2,1)$,
\begin{equation}\label{lin:profileinputs}
 \begin{split}
 &\int_{\R^2}w_{\app}^k\dd\xi=1,\qquad
   \int_{\R^2}\xi w_{\app}^k\dd\xi=0,\\
 &\mathfrak n_{1,\vartheta}(w_{\app}^k-G)\lesssim\delta,\\
 &\mathfrak n_{1,\vartheta}(w_{\app}^k-G-\delta F_k^0)
                      \lesssim\delta(\nu+\delta^{\frac12}),\\
 &|v^{w_{\app}^k}(\xi)-v^G(\xi)|
                      \lesssim\delta\av\xi^{-1},\\
 &|v^{w_{\app}^k}(\xi)-v^G(\xi)-\delta v^{F_k^0}(\xi)|
                      \lesssim\delta(\nu+\delta^{\frac12})\av\xi^{-1}.
 \end{split}
\end{equation}

\begin{lemma}\label{prof:radialcancel}
The source \eqref{prof:Csource} satisfies $\Pi_0C_k^0=0$, and for
$0<\vartheta<1$,
\begin{equation}\label{prof:radial-defect-bound}
 \mathfrak n_{0,\vartheta}(C_k^\nu-C_k^0)
 +\mathfrak n_{0,\vartheta}(\Pi_0C_k^\nu)\lesssim \nu.
\end{equation}
\end{lemma}
\begin{proof}
Since $\Sigma_k(s)$ is symmetric and trace free, we fix $s$
and rotate coordinates so that
$\Sigma_k=\operatorname{diag}(\sigma,-\sigma)$. The operators
$\Lambda$ and angular averaging commute with this rotation.
The right-hand side of the equation for $F_k^0$ is
$(\sigma/2)r^2G\cos2\theta$. On mode two, $\Lambda$ sends
the sine sector to the cosine sector; by uniqueness of its inverse,
we therefore have $F_k^0=f(r)\sin2\theta$. The polar components
of its velocity have the form $v^{F_k^0}\cdot \mathbf{e}_r=a(r)\cos2\theta$
and $v^{F_k^0}\cdot \mathbf{e}_\theta=b(r)\sin2\theta$.
Since
$(\Sigma_k\xi)\cdot \mathbf{e}_r=\sigma r\cos2\theta$ and
$(\Sigma_k\xi)\cdot \mathbf{e}_\theta=-\sigma r\sin2\theta$, we obtain
\[
 \begin{aligned}
 (\Sigma_k\xi)\cdot\nabla F_k^0
 &=\sigma(rf'-2f)\sin2\theta\cos2\theta,\\
 v^{F_k^0}\cdot\nabla F_k^0
 &=(af'+2bf/r)\sin2\theta\cos2\theta.
 \end{aligned}
\]
Both angular averages vanish.
The other terms in $C_k^0$ are $T_{k,4}\cdot\nabla G$,
of frequency four, and $F_k^0-\beta\partial_\theta F_k^0$,
of frequency two, so $\Pi_0C_k^0=0$.
This also covers the case $\Sigma_k=0$, for which $F_k^0=0$.

For $\vartheta<\vartheta_1<1$, using bilinearity,
\eqref{prof:transport-product-bound},
\eqref{prof:linear-transport-bound}, and \eqref{prof:F-bounds}, we obtain
\[
 \mathfrak n_{0,\vartheta}(C_k^\nu-C_k^0)
 \lesssim\bigl(1+\mathfrak n_{1,\vartheta_1}(F_k^\nu)
                  +\mathfrak n_{0,\vartheta_1}(F_k^0)\bigr)
       \mathfrak n_{1,\vartheta_1}(F_k^\nu-F_k^0)
 \lesssim\nu.
\]
Since angular averaging is contractive in $\mathfrak n_{0,\vartheta}$
and $\Pi_0C_k^0=0$, the same bound holds for $\Pi_0C_k^\nu$.
\end{proof}
The radial defect is $\delta\Pi_0C_k^\nu=O(\nu\delta)$.
The source $S_k^\nu$ has zero radial component by oddness.

\subsection{Error estimate}
With the prescribed center velocity, we set
\begin{equation}\label{prof:prescribed-residual}
 \mathcal R_k(\xi,s)
 :=\mathcal R_k[w_{\app};v^{w_{\app}}](\xi,s)
       +\frac{e_k}{\delta^{\frac12}}\cdot\nabla w_{\app}^k.
\end{equation}

\begin{proof}[Proof of Theorem~\ref{prof:construction-theorem}]
Using the angular sectors and moment conditions of the profiles,
we obtain unit mass and zero first moments for \eqref{prof:approx}.
By Lemma~\ref{prof:J-profile}, we have the bounds for $J_k^\nu$
and its $s$ derivative. Using Proposition~\ref{prof:localization},
we obtain the corresponding bounds for the approximation.
In Proposition~\ref{prof:residualprop} below, we compute $\mathcal R_k$
and prove its bound.
\end{proof}

\paragraph{The transport associated with the approximation.}
By \eqref{prof:center-law-definition} and
\eqref{prof:prescribed-residual}, we have
$\mathcal R_k=(\partial_s-\eta\Lop)w_{\app}^k+B_k\cdot\nabla w_{\app}^k$,
where
\begin{equation}\label{prof:corrected-transport}
 B_k(\xi,s)={}\sum_j\frac{\Gamma_j}{\delta} v^{w_{\app}^j}(\xi+\frac{a_{kj}}{\delta^{\frac12}}) -\frac1{\delta^{\frac12}}\sum_{j\ne k}\Gamma_jK(a_{kj}) -\delta^{\frac32}V_{2,k}^\nu-\beta\xi^\perp.
\end{equation}
With $e$ from \eqref{prof:center-law-definition}, the transport field
obtained from $w_{\app}$ and $Z'$ is
\begin{equation}\label{eq:Bapp}
 B_k^\app=B_k-\frac{e_k}{\delta^{\frac12}}.
\end{equation}
Using \eqref{prof:prescribed-residual} and \eqref{eq:Bapp}, we obtain
\begin{equation}\label{prof:actual-prescribed-identity}
 \mathcal R_k[w_{\app};v^{w_{\app}}]
 =(\partial_s-\eta\Lop)w_{\app}^k
       +B_k^\app\cdot\nabla w_{\app}^k
 =\mathcal R_k
                  -\frac{e_k}{\delta^{\frac12}}\cdot\nabla w_{\app}^k.
\end{equation}
All transport fields here are divergence free.

\paragraph{Pointwise and weighted bounds.}
\begin{proposition}\label{prof:residualprop}
Under \eqref{prof:two-parameter-range}, \eqref{prof:gap},
and \eqref{prof:center-velocity-bounds}, the approximation
\eqref{prof:approx} satisfies, for every $0<\vartheta<1$,
\begin{equation}\label{prof:residualpointwise}
 \mathfrak n_{0,\vartheta}(\mathcal R_k)
             \lesssim \delta(\delta+\nu).
\end{equation}
For every positive weight $p_k\le C_p e^{\mu|\xi|^2/4}$ with
$\mu<2\vartheta$,
\begin{equation}\label{prof:residualweighted}
 \|\sqrt{p_k}\mathcal R_k\|_{L^2}
       \lesssim \delta(\delta+\nu).
\end{equation}
The constants are uniform in $\nu,\delta$. Those in the bounds
for $\mathcal R_k$ may depend on $C_Z$; the constant in the
weighted estimate may also depend on $C_p$ and $\mu$.
\end{proposition}
\begin{proof}
We define $\mathcal X_k$ and $\mathcal Q_k$ by
\begin{equation}\label{prof:externaldecomposition}
 \begin{split}
 &\sum_{j\ne k}\frac{\Gamma_j}{\delta}v^{w_{\app}^j}(\xi+\frac{a_{kj}}{\delta^{\frac12}})
 -\frac1{\delta^{\frac12}}\sum_{j\ne k}\Gamma_jK(a_{kj})
       -\delta^{\frac32}V_{2,k}^\nu\\
 &\quad=\Sigma_k\xi+\delta^{\frac12}T_{k,3}+\delta T_{k,4}
       +\delta^{\frac32}(T_{k,5}-b_k^\nu)+\mathcal X_k
 \end{split}
\end{equation}
and
\begin{equation}\label{prof:quadratic-remainder}
 \mathcal Q_k={}v^{w_{\app}^k-G}\cdot\nabla(w_{\app}^k-G)
 -\delta^2v^{F_k^\nu}\cdot\nabla F_k^\nu
 -\delta^{\frac52}\bigl\{v^{F_k^\nu}\cdot\nabla H_k^\nu
                       +v^{H_k^\nu}\cdot\nabla F_k^\nu\bigr\}.
\end{equation}
Using \eqref{prof:scale-expansion}, \eqref{prof:profile-system},
and $\mathcal U_k^\nu-V_{2,k}^\nu=-b_k^\nu$, we obtain by a direct calculation
\begin{align}\label{prof:residualidentity}
 \mathcal R_k={}&\delta\Pi_0C_k^\nu
 +\nu\delta(2-\Lop)P_k^\nu
 +\nu\delta^{\frac32}(\tfrac52-\Lop)J_k^\nu\notag\\
 &+\delta\partial_sF_k^\nu+\delta^{\frac32}\partial_sH_k^\nu
       +\delta^2\partial_sP_k^\nu+\delta^{\frac52}\partial_sJ_k^\nu\notag\\
 &+\delta^2(2P_k^\nu+(\Sigma_k\xi-\beta\xi^\perp)\cdot\nabla P_k^\nu)\notag\\
 &+\delta^{\frac52}(\tfrac52J_k^\nu+(\Sigma_k\xi-\beta\xi^\perp)\cdot\nabla J_k^\nu)\notag\\
 &+\delta^{\frac12}T_{k,3}\cdot\nabla(w_{\app}^k-G-\delta F_k^\nu)
       +\delta T_{k,4}\cdot\nabla(w_{\app}^k-G)\notag\\
 &+\delta^{\frac32}(T_{k,5}-b_k^\nu)\cdot\nabla(w_{\app}^k-G)
       +\frac{\Gamma_k}{\delta}\mathcal Q_k
       +\mathcal X_k\cdot\nabla w_{\app}^k.
\end{align}

We fix $0<\vartheta<\vartheta_1<1$ and $0<c_*<d_*/8$.
For $|\xi|\le c_*/\sqrt\delta$, we use homogeneity and expand
through degree four to obtain
\begin{equation}\label{prof:external-taylor-bound}
 \begin{split}
 &\left|\sum_{j\ne k}\frac{\Gamma_j}{\delta}K(\frac{a_{kj}}{\sqrt\delta}+\xi)
 -\frac1{\sqrt\delta}\sum_{j\ne k}\Gamma_jK(a_{kj})
 -\sum_{n=2}^5\delta^{\frac{n-2}{2}}T_{k,n}(\xi)\right|
 \lesssim\delta^2|\xi|^5.
 \end{split}
\end{equation}
For the Oseen correction, we have
\[
 \delta^{-1}|v^G(\frac{a_{kj}}{\sqrt\delta}+\xi)
                   -K(\frac{a_{kj}}{\sqrt\delta}+\xi)|
 \lesssim\delta^{-\frac12}e^{-c/\delta}.
\]
For the remote mode-two profile, we expand
$K(a_{kj}+\sqrt\delta\,\xi-\sqrt\delta\,y)$ to second order
for $|y|\le d_*/(8\sqrt\delta)$. For $0\le\tau\le1$, we have
$|a_{kj}+\sqrt\delta\,\xi-\tau\sqrt\delta\,y|>3d_*/4$.
Using the Gaussian bounds, we obtain
$\|\av y^mF_j^\nu\|_{L^1(|y|\ge d_*/(8\sqrt\delta))}
\lsim{m} e^{-c/\delta}$ for every fixed $m$.
For the omitted kernel term, we obtain
\[
 \begin{split}
 &\left|\int_{|y|\ge d_*/(8\sqrt\delta)}
 K(\frac{a_{kj}}{\sqrt\delta}+\xi-y)F_j^\nu(y)\dd y\right|\\
 &\qquad\lesssim
 \|F_j^\nu\|_{L^1(|y|\ge d_*/(8\sqrt\delta))}^{\frac12}
 \|F_j^\nu\|_{L^\infty(|y|\ge d_*/(8\sqrt\delta))}^{\frac12}
 \lesssim e^{-c/\delta}.
 \end{split}
\]
Using the zero mass and first moments of $F_j^\nu$, we therefore obtain
\[
 v^{F_j^\nu}(\frac{a_{kj}}{\sqrt\delta}+\xi)
 =\frac{\delta^{\frac32}}2\int_{\R^2}
 D^2K(a_{kj}+\sqrt\delta\,\xi)[y,y]F_j^\nu(y)\dd y+O(\delta^2).
\]
Using the third derivative bound and the uniformly bounded second
absolute moments, we obtain
\begin{equation}\label{prof:remote-quadrupole-expansion}
 \sum_{j\ne k}\Gamma_jv^{F_j^\nu}(\frac{a_{kj}}{\sqrt\delta}+\xi)
 =\delta^{\frac32}\mathcal U_k^\nu+O(\delta^2\av\xi).
\end{equation}
For the remaining profiles, we use \eqref{prof:profile-velocities} to obtain
\[
 \delta^{\frac12}|v^{H_j^\nu}(\frac{a_{kj}}{\sqrt\delta}+\xi)|
 +\delta|v^{P_j^\nu}(\frac{a_{kj}}{\sqrt\delta}+\xi)|
 +\delta^{\frac32}|v^{J_j^\nu}(\frac{a_{kj}}{\sqrt\delta}+\xi)|
 \lesssim\delta^{\frac52}.
\]
Combining these estimates with $\mathcal U_k^\nu-V_{2,k}^\nu=-b_k^\nu$,
we obtain $|\mathcal X_k|\lesssim\delta^2\av\xi^5$
on $|\xi|\le c_*/\sqrt\delta$.

On $|\xi|>c_*/\sqrt\delta$, we use the uniform $L^\infty$ velocity
bounds and \eqref{prof:externaldecomposition} to obtain
$|\mathcal X_k|\lesssim\delta^{-1}\av\xi^4
\lesssim\delta^2\av\xi^{10}$. Hence we have
$|\mathcal X_k(\xi)|\lesssim\delta^2\av\xi^{10}$ on $\R^2$,
and using \eqref{prof:derivativebounds}, we obtain
$\mathfrak n_{0,\vartheta}
(\mathcal X_k\cdot\nabla w_{\app}^k)\lesssim\delta^2$.
Using bilinearity and the same profile bounds, we obtain
$\mathfrak n_{0,\vartheta}(\mathcal Q_k)\lesssim\delta^3$.
In \eqref{prof:residualidentity}, we bound the radial and viscous terms
by $C\nu\delta$ using \eqref{prof:radial-defect-bound}
and \eqref{prof:derivativebounds}, the time-derivative terms by
$C\delta|Z'|$, and the remaining terms by $C\delta^2$.
Here we use $\delta\le1$, $\vartheta<\vartheta_1$,
$\mathfrak n_{1,\vartheta_1}(w_{\app}^k-G-\delta F_k^\nu)
\lesssim\delta^{\frac32}$, and
$\mathfrak n_{1,\vartheta_1}(w_{\app}^k-G)\lesssim\delta$.
Using \eqref{prof:center-velocity-bounds}, we therefore obtain
$\mathfrak n_{0,\vartheta}(\mathcal R_k)
\lsim{C_Z}\nu\delta+\delta^2+\delta|Z'|
\lesssim\delta(\delta+\nu)$.
Finally, for $p_k\le C_pe^{\mu|\xi|^2/4}$ with $\mu<2\vartheta$, we have
$\|\sqrt{p_k}\mathcal R_k\|_{L^2}
\lsim{C_p,\vartheta,\mu}\mathfrak n_{0,\vartheta}(\mathcal R_k)$,
and hence we obtain \eqref{prof:residualweighted}.
\end{proof}

\subsection{The remainder equation}
We define the remainders $f_k$ by
\begin{equation}\label{eq:error-definition}
 w_k=w_{\app}^k+\delta f_k,\qquad v^{w_k}=v^{w_{\app}^k}+\delta v^{f_k}.
\end{equation}
We then have
\begin{equation}\label{eq:error-zero-moments}
 \int_{\R^2} f_k\dd\xi=0,\qquad\int_{\R^2}\xi f_k\dd\xi=0,
\end{equation}
because the exact and approximate profiles have the same mass and
first moments. Substituting \eqref{eq:error-definition} into
\eqref{eq:inner-equation} and using
\eqref{prof:actual-prescribed-identity} and $\delta'/\delta=1+\eta$,
we obtain
\begin{equation}\label{eq:error-equation}
 \begin{split}
 \partial_sf_k={}&\eta\Lop f_k-(1+\eta)f_k-B_k^\app\cdot\nabla f_k -\frac1\delta\sum_j\Gamma_jv^{f_j}(\xi+\frac{a_{kj}}{\delta^{\frac12}})\cdot\nabla w_{\app}^k\\
 &-\sum_j\Gamma_jv^{f_j}(\xi+\frac{a_{kj}}{\delta^{\frac12}})\cdot\nabla f_k -\frac{\mathcal R_k}{\delta} +\frac{e_k}{\delta^{\frac32}}\cdot\nabla w_{\app}^k.
 \end{split}
\end{equation}

\section{Adapted weights and weighted norms}\label{wgt:section}
For the estimates in Section~\ref{sec:nonlinear}, we choose weights adapted
to the approximate Euler equilibria near the vortex centers. We extend
these weights to the whole plane to control moments, stream functions,
and velocities.

We work on the interval
\begin{equation}\label{wgt:parameter-range}
 0<s\le S_\nu=\gamma\log(1/\nu),
\end{equation}
where we have
\begin{equation}\label{wgt:scale-range}
 0<\delta(s)=\nu(e^s-1)\le\nu^{1-\gamma}-\nu\le\nu^{1-\gamma}.
\end{equation}
We consider a smooth center path $Z(s)$ satisfying the gap
\eqref{prof:gap} and the diameter bound
\begin{equation}\label{wgt:diameter}
 \max_{j\ne k}|Z_k-Z_j|
       \le\max_{j\ne k}|z_k^0-z_j^0|+1.
\end{equation}
The constants below may depend on the fixed configuration and gap,
on $\gamma$, and on $\kappa$, but are uniform in $s$ and $\nu$.
We use the stream function
\begin{equation}\label{wgt:stream-definition}
 \psi^f(\xi)=\frac1{2\pi}\int_{\R^2}\log|\xi-y|f(y)\dd y,
 \qquad\Delta\psi^f=f,\qquad v^f=\nabla^\perp\psi^f.
\end{equation}
For a positive weight $p$, we set
$\|f\|_{L^2_p}^2=\int_{\R^2}p|f|^2\dd\xi$.

\Needspace{6\baselineskip}

\begin{proposition}\label{wgt:main}
For sufficiently small $\nu$, there are positive weights $p_k$,
locally Lipschitz in space and time, such that
$\av\xi^{16}\lesssim p_k(\xi,s)\lesssim e^{|\xi|^2/4}$.
For $f_k\in L^2_{p_k}$ satisfying
$\int_{\R^2}f_k\dd\xi=0$ and
$\int_{\R^2}\xi f_k\dd\xi=0$, define
\begin{equation}\label{lin:energy-def}
 \mathcal E=\frac12\sum_k\left(
       \int_{\R^2}p_kf_k^2\dd\xi
       +\int_{\R^2}\psi^{f_k}f_k\dd\xi\right).
\end{equation}
Then
\begin{equation}\label{wgt:energy-equivalence}
 c\sum_k\|f_k\|_{L^2_{p_k}}^2\le\mathcal E
       \le\frac12\sum_k\|f_k\|_{L^2_{p_k}}^2
\end{equation}
for a fixed $c>0$. For every fixed $m\ge0$, every $2<q<\infty$,
and all $j,k$,
\begin{equation}\label{wgt:main-controls}
 \begin{aligned}
 \sum_k\bigl(\|\av\xi^mf_k\|_{L^1}
       +\|f_k\|_{L^{\frac{4}{3}}}
       +\|\psi^{f_k}\|_{H^2}
       +\|\psi^{f_k}\|_{L^\infty}\bigr)
       &\lsim{m}\sqrt{\mathcal E},\\
 \left\||\nabla\log p_k(\xi)|\,
       v^{f_j}\left(\xi+\frac{a_{kj}}{\sqrt\delta}\right)\right\|_{L^q_\xi}
       &\lsim{q}\sqrt{\mathcal E}.
 \end{aligned}
\end{equation}
The required smallness of $\nu$ is independent of $m$ and $q$.
\end{proposition}

We use the moment bounds to control the forcing of the vortex centers
and the stream-function and velocity bounds to control the terms
obtained by testing the remainder equation against $p_kf_k+\psi^{f_k}$.
We also define
\begin{equation}\label{lin:dissipation-def}
 \mathcal D=\sum_k\|\nabla(\sqrt{p_k}f_k)\|_{L^2}^2.
\end{equation}
The estimates above require only \eqref{prof:gap} and \eqref{wgt:diameter}.
For the time derivative estimates, we also assume
\begin{equation}\label{wgt:centerassumptions}
 |Z'|+|Z'-\mathcal F(Z)|\le C_Z\delta^2,
 \qquad \delta'=(1+\eta)\delta,\quad\eta=\nu/\delta;
\end{equation}
the corresponding constants may depend on $C_Z$.
In Section~\ref{sec:proof-main}, we verify this bound under the coupled bootstrap.

\subsection{From Gaussian cancellation to Arnold energy}\label{wgt:Arnold-structure}
The cancellation in the inverse-Gaussian energy relies on the oddness
and perpendicularity of the Biot--Savart kernel and on the Gaussian
identity $\nabla\log G=-\xi/2$. We consider the forced linearized Euler equation
\begin{equation}\label{wgt:Gaussian-toy}
 \partial_s f+\frac1\nu\Lambda f=g,\qquad
 \Lambda f=v^G\cdot\nabla f+v^f\cdot\nabla G.
\end{equation}
Here $g$ is prescribed and $\nu>0$ is fixed. The following identities
hold for smooth functions with sufficient decay to justify the integrals
and integrations by parts.

\paragraph{The inverse-Gaussian energy.}
Following Gallay and Wayne \cite{GW05}, we define
$E_G(f)=\frac12\int_{\R^2}\frac{f^2}{G}\dd\xi$.
Using $v^G\cdot\nabla G=0$ and $\nabla\log G=-\xi/2$,
we compute along \eqref{wgt:Gaussian-toy}
\[
 \frac{\dd}{\dd s}E_G(f)
 =\frac1{2\nu}\int f v^f\cdot\xi+\int\frac{fg}{G}.
\]
Using the oddness and perpendicularity of $K$, we obtain
\begin{equation}\label{wgt:Gaussian-antisymmetry}
 \int f v^f\cdot\xi
 =\frac12\iint f(\xi)f(y)K(\xi-y)\cdot(\xi-y)\dd y\dd\xi=0.
\end{equation}
Consequently, we have $\frac{\dd}{\dd s}E_G(f)=\int_{\R^2}\frac{fg}{G}\dd\xi$.
In particular, $E_G$ is conserved when $g=0$, and the forced estimate
\[
 \|f(s)\|_{L^2_{G^{-1}}}
 \le\|f(0)\|_{L^2_{G^{-1}}}
       +\int_0^s\|g(\sigma)\|_{L^2_{G^{-1}}}\dd\sigma
\]
has no loss as $\nu\to0$.

\paragraph{Beyond the Gaussian profile.}
For a positive profile $W$, the corresponding term in the energy
identity is
\[
 \begin{aligned}
 -\int\frac fW v^f\cdot\nabla W &=-\int f v^f\cdot\nabla\log W\\
 &=-\frac12\iint f(\xi)f(y)K(\xi-y) \cdot\{\nabla\log W(\xi)-\nabla\log W(y)\}\dd y\dd\xi.
 \end{aligned}
\]
For $W=G$, the difference in braces is $-(\xi-y)/2$, so we obtain zero
by the perpendicularity of the kernel. For a general $W$, this difference
need not be parallel to $\xi-y$, and the oddness of the kernel alone does not
cancel this term. 

\paragraph{The Arnold energy for a general equilibrium.}
We recall the cancellation underlying the adapted Arnold-energy
method for weakly viscous vortex asymptotics, developed by
Gallay and \v{S}ver\'ak \cite{GS24}, Dolce and Gallay \cite{DoGa26},
and Donati and Gallay \cite{DG26}.
We take a smooth stationary Euler velocity--vorticity pair $(U,W)$,
with $U=v^W=\nabla^\perp\psi^W$ and $U\cdot\nabla W=0$,
and consider the linearized equation
\begin{equation}\label{wgt:general-toy}
 \partial_s f+\frac1\nu
       (U\cdot\nabla f+v^f\cdot\nabla W)=g.
\end{equation}
For a smooth positive time-independent weight $p$, we define
$E_p(f)=\frac12\int_{\R^2}(pf^2+\psi^f f)\dd\xi$.
By a direct calculation, we obtain
\begin{equation}\label{wgt:general-Arnold-balance}
 \begin{aligned}
 \frac{\dd}{\dd s}E_p(f)
 &=\frac1{2\nu}\int f^2U\cdot\nabla p
   +\frac1\nu\int f(U+p\nabla^\perp W)\cdot\nabla\psi^f\\
 &\quad+\int(pf+\psi^f)g,
 \end{aligned}
\end{equation}
where we have used $v^f\cdot\nabla\psi^f=0$.
Thus the singular terms vanish if $p$ solves
\begin{equation}\label{wgt:cancellation-system}
 \begin{cases}
  U+p\nabla^\perp W=0,\\
  U\cdot\nabla p=0.
 \end{cases}
\end{equation}
For the approximate vortices below, we need only satisfy
\eqref{wgt:cancellation-system} up to errors controlled by the
energy estimate. We estimate contributions bounded by
$O(\nu^M)\|f\|_{L^2_p}^2$ directly, choosing $M$ sufficiently large
to account for all singular prefactors.
We can therefore modify the weight in the tails.
For exact cancellation, it suffices to impose the following
vorticity--stream-function relation.

\begin{proposition}[A weight from the equilibrium relation]
\label{wgt:Arnold-weight}
Let $W$ be a smooth time-independent profile and set $U=v^W$.
Assume that there exist a single-valued function $\mathcal H\in C^2(I)$
on an interval $I$ with $W(\R^2)\subset I$ and a constant $c$
such that
\begin{equation}\label{wgt:general-equilibrium}
 \psi^W+\mathcal H(W)=c,\qquad \mathcal H'(W)>0.
\end{equation}
Then $p=\mathcal H'(W)$ is a positive weight satisfying
\eqref{wgt:cancellation-system}. For every smooth solution of
\eqref{wgt:general-toy} for which the integrals and integrations by
parts above are justified,
\[
 \frac{\dd}{\dd s}E_p(f)=\int_{\R^2}(pf+\psi^f)g\dd\xi.
\]
In particular, $E_p$ is conserved when $g=0$.
\end{proposition}

\begin{proof}
Differentiating \eqref{wgt:general-equilibrium}, we obtain
$U+p\nabla^\perp W=0$ and
$U\cdot\nabla p=\mathcal H''(W)U\cdot\nabla W=0$.
Using \eqref{wgt:general-Arnold-balance}, we obtain the asserted energy identity.
\end{proof}
\paragraph{The Gaussian reference.}
For the Gaussian equilibrium, we have $W=G$ and $U=v^G$. Since $G$ is
radial and strictly decreasing, we define $\mathcal H:(0,G(0)]\to\R$
by $\psi^G+\mathcal H(G)=0$. With $r=|\xi|$, the corresponding
weight is
\begin{equation}\label{wgt:entropy-definition}
 \mathsf A(r):=\mathcal H'(G(r))
       =-\frac{(\psi^G)'(r)}{G'(r)}
       =\frac{4(e^{r^2/4}-1)}{r^2},\qquad \mathsf A(0)=1.
\end{equation}
It satisfies the cancellation system \eqref{wgt:cancellation-system}:
$v^G+\mathsf A\nabla^\perp G=0$ and
$v^G\cdot\nabla\mathsf A=0$.

For the deformed vortex $W_k=G+\delta F_k^0$, we include the strain
from the other vortices in the velocity:
\[
 U_k=v^{W_k}+\frac\delta{\Gamma_k}\Sigma_k\xi
     =v^G+\delta v^{F_k^0}+\frac\delta{\Gamma_k}\Sigma_k\xi.
\]
Accordingly, we replace the Gaussian weight
$\mathsf A=\mathcal H'(G)$ by the inner weight
$\mathsf A_k^\delta=\mathcal H'(W_k)=\mathcal H'(G+\delta F_k^0)$.
We extend the inner weight to a global weight $p_k$ that agrees with this
formula near the $k$th vortex center, is constant near the other vortex
centers, and grows at large distances to control moments.

\subsection{The inner weight}
We denote by $\mathsf S_k$ the homogeneous harmonic quadratic polynomial
with $\nabla^\perp\mathsf S_k=\Sigma_k\xi$, and we set
\begin{equation}\label{wgt:leading-stream}
 \Psi_k=\psi^{W_k}+\frac\delta{\Gamma_k}\mathsf S_k,
 \qquad U_k=\nabla^\perp\Psi_k.
\end{equation}
Using the Euler profile equation, we derive a scalar relation that cancels
the term of order $\delta$ in $\Psi_k+\mathcal H(W_k)$.
We also estimate the relative deformation $F_k^0/G$.

\begin{lemma}\label{wgt:Euler-profile}
The Euler profile satisfies
\begin{equation}\label{wgt:Euler-relation}
 \psi^{F_k^0}+\mathsf A F_k^0+\mathsf S_k/\Gamma_k=0.
\end{equation}
For $r=|\xi|$, it also satisfies
\begin{equation}\label{wgt:profile-ratio}
 \begin{aligned}
 \left|\frac{F_k^0}{G}\right|&\lesssim r^2\av r^2,\\
 \left|\nabla\left(\frac{F_k^0}{G}\right)\right|
       &\lesssim r\av r^2,\\
 \left|\partial_s\left(\frac{F_k^0}{G}\right)\right|
       &\lesssim |Z'(s)|r^2\av r^2.
 \end{aligned}
\end{equation}
\end{lemma}
\begin{proof}
Since $v^G=-\mathsf A\nabla^\perp G$ and $\mathsf A$ is radial,
we obtain from the profile equation
\[
 0=\Gamma_k\Lambda F_k^0+(\Sigma_k\xi)\cdot\nabla G
 =\frac{G(r)}2\partial_\theta
       (\Gamma_k\psi^{F_k^0}+\Gamma_k\mathsf A F_k^0+\mathsf S_k).
\]
The expression in parentheses is independent of the angle and has only
angular mode two, so it vanishes. Thus, using \eqref{wgt:Euler-relation}, we obtain
\[
 \frac{F_k^0}{G}
 =-\bigl(\psi^{F_k^0}+\mathsf S_k/\Gamma_k\bigr)(\mathsf A G)^{-1},
 \qquad
 (\mathsf A G)^{-1}=\frac{\pi r^2}{1-e^{-r^2/4}}.
\]
The last factor is smooth at the origin. By direct differentiation, we obtain
$(\mathsf A G)^{-1}\lesssim\av r^2$ and
$|\nabla(\mathsf A G)^{-1}|\lesssim r$.
Using the mode-two structure, \eqref{prof:mode-velocity-input}, and
\eqref{prof:F-bounds}, we obtain, for $r\ge1$,
\[
 \begin{aligned}
 |\psi^{F_k^0}|&\lesssim r^{-2},&
 |\nabla\psi^{F_k^0}|&\lesssim r^{-3},\\
 |\mathsf S_k|&\lesssim r^2,&
 |\nabla\mathsf S_k|&\lesssim r.
 \end{aligned}
\]
On $r\le1$, using smoothness and angular mode two, we obtain
$|F_k^0|\lesssim r^2$ and $|\nabla F_k^0|\lesssim r$.
Combining these bounds, we obtain the first two estimates in
\eqref{wgt:profile-ratio}. Differentiating the expression for the profile ratio, we obtain
\[
 \partial_s\left(\frac{F_k^0}{G}\right)
 =-\bigl(\partial_s\psi^{F_k^0}
             +\partial_s\mathsf S_k/\Gamma_k\bigr)(\mathsf A G)^{-1}.
\]
By \eqref{prof:strain-time-dependence}, we have
$|\partial_s\mathsf S_k|\lesssim|Z'|r^2$. Using the mode-two structure,
\eqref{prof:F-bounds}, and \eqref{prof:Gaussian-velocity-bound}, we obtain
$|\partial_s\psi^{F_k^0}|\lesssim|Z'|r^{-2}$ for $r\ge1$ and
$|\partial_s\psi^{F_k^0}|\lesssim|Z'|r^2$ for $r\le1$.
Substituting these bounds into the differentiated identity, we obtain the remaining estimate.
\end{proof}

We set
\begin{equation}\label{wgt:scales}
 r_\delta=\delta^{-\kappa}.
\end{equation}
We have $r_\delta\to\infty$,
$\delta r_\delta^4\to0$,
$\sqrt\delta\,r_\delta^3=\delta^{\frac{1}{2}-3\kappa}\to0$, and
$\delta r_\delta^6\to0$.

To define the inner weight, we verify that
$\{G(\xi)+\tau\delta(s)F_k^0(\xi,s):0\le\tau\le1\}\subset(0,G(0)]$
for $|\xi|\le2r_{\delta(s)}$.

\begin{lemma}\label{wgt:Gaussian-segment}
For sufficiently small $\nu$,
\begin{equation}\label{wgt:Gaussian-segment-range}
 \frac12G(\xi)
 \le G(\xi)+\tau\delta(s)F_k^0(\xi,s)\le G(0)
\end{equation}
for every $0<s\le S_\nu$, $|\xi|\le2r_{\delta(s)}$, and
$0\le\tau\le1$.
\end{lemma}
\begin{proof}
For $r\le2r_\delta$, we obtain from \eqref{wgt:profile-ratio}
\[
 \left|\delta\frac{F_k^0}{G}\right|
 \le C\delta r^2\av r^2\le C\delta^{1-4\kappa}\le\frac12
\]
for sufficiently small $\nu$. We thus obtain the lower bound. For the upper bound, we have
\[
 \log\frac{G+\tau\delta F_k^0}{G(0)} =-\frac{r^2}{4} +\log\left(1+\tau\delta\frac{F_k^0}{G}\right) \le-\frac{r^2}{4}+C\delta r^2\av r^2\le0,
\]
since $C\delta(1+4r_\delta^2)\le1/4$ for sufficiently small $\nu$.
\end{proof}

We therefore define
\begin{equation}\label{wgt:Hessian-definition}
 \mathsf A_k^\delta(\xi,s)
 =\mathcal H'\!\left(G(\xi)+\delta(s)F_k^0(\xi,s)\right),
 \qquad |\xi|\le2r_{\delta(s)}.
\end{equation}

The global weight $p_k$, defined below, will satisfy
$p_k=\mathsf A_k^\delta$ on $\mathcal I_k(s)$, defined in
\eqref{wgt:inner-region}.

We also define $\Theta_k:=\Psi_k+\mathcal H(W_k)$ for $|\xi|\le2r_\delta$.

\begin{lemma}\label{wgt:Theta-remainder}
For sufficiently small $\nu$, the identity
\begin{equation}\label{wgt:Theta-Taylor}
 \Theta_k=\delta^2(F_k^0)^2\int_0^1(1-\tau)
                 \mathcal H''(G+\tau\delta F_k^0)\dd\tau
\end{equation}
holds for $0<s\le S_\nu$ and $|\xi|\le2r_{\delta(s)}$.
\end{lemma}
\begin{proof}
By Lemma~\ref{wgt:Gaussian-segment}, we may apply Taylor's formula.
Using $\psi^G+\mathcal H(G)=0$ and \eqref{wgt:Euler-relation},
we cancel the constant and linear terms and obtain \eqref{wgt:Theta-Taylor}.
\end{proof}

On $\{|\xi|\le2r_{\delta(s)}\}$, the cancellation system
\eqref{wgt:cancellation-system} has the remainders
\begin{equation}\label{wgt:approximate-cancellation}
 \begin{cases}
 U_k+\mathsf A_k^\delta\nabla^\perp W_k=\nabla^\perp\Theta_k,\\
 U_k\cdot\nabla\mathsf A_k^\delta
       =\nabla^\perp\Theta_k\cdot\nabla\mathsf A_k^\delta.
 \end{cases}
\end{equation}

\begin{lemma}\label{wgt:inner-weight}
On $r\le2r_\delta$, the inner weight $\mathsf A_k^\delta$ satisfies
\begin{equation}\label{wgt:Hessian-bounds}
 \begin{split}
 |\mathsf A_k^\delta/\mathsf A-1|&\lesssim\delta r^2\av r^2,\\
 |\nabla\log\mathsf A_k^\delta-\nabla\log\mathsf A|
       &\lesssim\delta r\av r^2,\\
 |\partial_\theta\log\mathsf A_k^\delta|
       &\lesssim\delta r^2\av r^2,\\
 |\partial_s\log\mathsf A_k^\delta|
       &\lesssim(1+\eta)\delta r^2\av r^2.
 \end{split}
\end{equation}
Moreover,
\begin{equation}\label{wgt:stream-defect}
 |\nabla\Theta_k|\lesssim\delta^2\av r^5.
\end{equation}
\end{lemma}
\begin{proof}
We set $\rho_k=F_k^0/G$, $u=r^2/4$,
$b(t)=\log\frac{e^t-1}{t}$ with $b(0)=0$, and
$\widetilde u_k=u-\log(1+\delta\rho_k)$.
By Lemma~\ref{wgt:Gaussian-segment} and \eqref{wgt:profile-ratio}, we have
$\widetilde u_k\ge0$ and
$|\widetilde u_k-u|\lesssim\delta|\rho_k|\ll1$. Moreover,
by \eqref{wgt:entropy-definition}, we have
$\log\mathsf A=b(u)$ and $\log\mathsf A_k^\delta=b(\widetilde u_k)$.
By a direct calculation, we obtain
\[
 0\le b'(t)=\frac{e^t}{e^t-1}-\frac1t\le1,\qquad
 b''(t)=\frac1{t^2}-\frac1{4\sinh^2(t/2)},\qquad
 |b''(t)|\lesssim(1+t)^{-2}.
\]
Consequently, we obtain
\[
 \begin{aligned}
 |\log(\mathsf A_k^\delta/\mathsf A)|
 &\lesssim\delta|\rho_k|,\\
 |\nabla\log\mathsf A_k^\delta-\nabla\log\mathsf A|
 &\lesssim\frac{\delta|\rho_k|r}{(1+u)^2}
              +\delta|\nabla\rho_k|\lesssim\delta r\av r^2,\\
 |\partial_s\log\mathsf A_k^\delta|
 &\lesssim\delta\bigl((1+\eta)|\rho_k|+|\partial_s\rho_k|\bigr)
       \lesssim(1+\eta)\delta r^2\av r^2,
 \end{aligned}
\]
where we have used \eqref{wgt:profile-ratio} and
\eqref{wgt:centerassumptions}. Exponentiating the first estimate and
using $\partial_\theta=\xi^\perp\cdot\nabla$ in the second, we obtain
\eqref{wgt:Hessian-bounds}.

By \eqref{wgt:Theta-Taylor}, we have
\[
 \Theta_k=\delta^2\rho_k^2\int_0^1(1-\tau)
                 G^2\mathcal H''(G+\tau\delta F_k^0)\dd\tau.
\]
We set $\widetilde u_{k,\tau}=u-\log(1+\tau\delta\rho_k)$.
Using Lemma~\ref{wgt:Gaussian-segment} and \eqref{wgt:profile-ratio}, we obtain
$\widetilde u_{k,\tau}\ge0$,
$1+\widetilde u_{k,\tau}\sim1+u$, and
$|\nabla\widetilde u_{k,\tau}|\lesssim r$.
By a direct calculation, we obtain
\[
 G^2\mathcal H''(G+\tau\delta F_k^0)
 =-\frac1{4\pi(1+\tau\delta\rho_k)^2}
       \int_0^1(1-t)e^{-t\widetilde u_{k,\tau}}\dd t.
\]
We bound the integral and its derivative by
$C(1+\widetilde u_{k,\tau})^{-1}$ and
$C(1+\widetilde u_{k,\tau})^{-2}$, respectively. We therefore obtain
\[
 \begin{aligned}
 |G^2\mathcal H''(G+\tau\delta F_k^0)| &\lesssim\av r^{-2},\\
 |\nabla[G^2\mathcal H''(G+\tau\delta F_k^0)]|
 &\lesssim\frac{\delta|\nabla\rho_k|}{1+u}
       +\frac{|\nabla\widetilde u_{k,\tau}|}{(1+u)^2}
 \lesssim\av r^{-3},
 \end{aligned}
\]
using $\delta\av r^4\lesssim1$ for $r\le2r_\delta$. Combining these bounds
with the integral expression above, we obtain
\[
 |\nabla\Theta_k|
 \lesssim\delta^2\left(
     \frac{|\rho_k||\nabla\rho_k|}{\av r^2}
       +\frac{|\rho_k|^2}{\av r^3}\right)
 \lesssim\delta^2\av r^5.
\]
\end{proof}

\subsection{The global weight}
We set
\[
 h_\delta(s)=\log\mathsf A(r_{\delta(s)}),\qquad
 R(s)=\frac{L_0}{\sqrt{\delta(s)}},\qquad
 L_0>2\bigl(\max_{j\ne k}|z_k^0-z_j^0|+1\bigr)+2.
\]
We write $r_\delta$ for $r_{\delta(s)}=\delta(s)^{-\kappa}$ and
$R$ for $R(s)$.
For small $\nu$, these choices satisfy
\begin{equation}\label{wgt:separation-scales}
 R>6r_\delta,\qquad12r_\delta\sqrt\delta\le d_*.
\end{equation}
We write $\xi=r\mathbf{e}_r$, where $\mathbf{e}_r=(\cos\theta,\sin\theta)$.
We define the capping radius by the level $e^{h_\delta(s)}$.
This radius depends on time through both $\delta(s)$ and $F_k^0(\cdot,s)$.

\begin{lemma}\label{wgt:interface}
For sufficiently small $\nu$, the function
$r\mapsto\mathsf A_k^\delta(r\mathbf{e}_r,s)$ is strictly increasing on
$[0,2r_\delta]$. For every $k$, $\theta$, and $0<s\le S_\nu$, there
is a unique radius $r_k(\theta,s)$ satisfying
\begin{equation}\label{wgt:interface-definition}
 \mathsf A_k^{\delta(s)}\bigl(r_k(\theta,s)\mathbf{e}_r,s\bigr)
       =e^{h_\delta(s)},\qquad
 \frac{r_{\delta(s)}}2<r_k(\theta,s)<2r_{\delta(s)}.
\end{equation}
The function $r_k$ is locally $C^1$ in $(\theta,s)$.
\end{lemma}
\begin{proof}
Since $1/2\le b'\le1$, we have
\begin{equation}\label{wgt:Gaussian-log-derivative}
 \frac r4\le\partial_r\log\mathsf A(r)
 =\frac r2 b'\left(\frac{r^2}{4}\right)\le\frac r2.
\end{equation}
Combining this with \eqref{wgt:Hessian-bounds} and
$\delta\av{2r_\delta}^2\to0$, we obtain
\begin{equation}\label{wgt:interface-monotonicity}
 \frac r8\le\partial_r\log\mathsf A_k^\delta(r\mathbf{e}_r,s)
             \le\frac{5r}{8},\qquad 0<r\le2r_\delta.
\end{equation}
Moreover, using \eqref{wgt:Hessian-bounds} and
\eqref{wgt:Gaussian-log-derivative}, we obtain
\[
 \begin{aligned}
 \log\mathsf A_k^\delta\bigl((r_\delta/2)\mathbf{e}_r,s\bigr)-h_\delta(s)
 &\le C\delta r_\delta^4-\frac{3}{32}r_\delta^2<0,\\
 \log\mathsf A_k^\delta\bigl(2r_\delta \mathbf{e}_r,s\bigr)-h_\delta(s)
 &\ge\frac38 r_\delta^2-C\delta r_\delta^4>0,
 \end{aligned}
\]
for sufficiently small $\nu$. By continuity and
\eqref{wgt:interface-monotonicity}, we obtain existence and uniqueness.
Applying the implicit function theorem, we obtain local $C^1$ regularity.
\end{proof}

Using this radius, we define the global weight by
\begin{equation}\label{wgt:definition}
 p_k(\xi,s)=
 \begin{cases}
 \mathsf A_k^\delta(\xi,s),&0\le r\le r_k(\theta,s),\\
 e^{h_\delta(s)},&r_k(\theta,s)<r\le R,\\
 e^{h_\delta(s)+\sqrt\nu(r-R)},&r>R.
 \end{cases}
\end{equation}
By Lemma~\ref{wgt:interface} and \eqref{wgt:interface-definition}, we may
equivalently write
\[
 p_k(\xi,s)=
 \begin{cases}
 \min\{\mathsf A_k^\delta(\xi,s),e^{h_\delta(s)}\},&0\le r\le2r_\delta,\\
 e^{h_\delta(s)},&2r_\delta<r\le R,\\
 e^{h_\delta(s)+\sqrt\nu(r-R)},&r>R.
 \end{cases}
\]
We write $Q_k$ for the logarithmic weight and $\mathcal I_k(s)$ for its uncapped region:
\begin{equation}\label{wgt:inner-region}
 Q_k=\log p_k,\qquad
 \mathcal I_k(s)=\{\xi:r<r_k(\theta,s)\}.
\end{equation}
We denote the plateau and outer regions by
\begin{align}
 \mathcal P_k(s)&:=\{\xi=r\mathbf{e}_r:r_k(\theta,s)<r<R(s)\},
       \label{wgt:plateau-region}\\
 \mathcal O_k(s)&:=\{\xi:|\xi|>R(s)\}.
       \label{wgt:outer-region}
\end{align}
For comparison with $p_k$, we define the radial weight $p_0$
by replacing $\mathsf A_k^\delta$ with $\mathsf A$:
\begin{equation}\label{wgt:radial-weight}
 p_0(r,s)=
 \begin{cases}
 \mathsf A(r),&0\le r\le r_\delta,\\
 \mathsf A(r_\delta),&r_\delta<r\le R,\\
 \mathsf A(r_\delta)e^{\sqrt\nu(r-R)},&r>R.
 \end{cases}
\end{equation}
We use $p_0$ for radial estimates of growth and inverse moments, then
transfer these bounds to $p_k$ by \eqref{wgt:comparison}.
In the energy, we retain $p_k$ because its dependence on
$F_k^0$ is needed for the adapted Arnold cancellation.

\begin{lemma}\label{wgt:regularity}
The weights $p_0$ and $p_k$ are locally Lipschitz in space and time
for $s>0$, and satisfy
\begin{equation}\label{wgt:comparison}
 \tfrac12p_0\le p_k\le\tfrac32p_0.
\end{equation}
\end{lemma}
\begin{proof}
The smooth pieces in \eqref{wgt:definition} and
\eqref{wgt:radial-weight} agree at their locally $C^1$ interfaces,
so both weights are locally Lipschitz. On $r\le2r_\delta$,
we obtain from \eqref{wgt:Hessian-bounds}
$\frac12\mathsf A\le\mathsf A_k^\delta\le\frac32\mathsf A$.
We preserve this comparison by capping both weights at the common level
$e^{h_\delta}$. For $r>2r_\delta$, we have $p_k=p_0$.
\end{proof}

\begin{lemma}\label{wgt:geometry}
At each fixed $s$, the boundary radius $r_k(\theta,s)$ and the inner
region $\mathcal I_k(s)$ satisfy
\begin{equation}\label{wgt:interface-location}
 |r_k-r_\delta|+|\partial_\theta r_k|\lesssim\delta r_\delta^3,
 \qquad \{r<r_\delta/2\}\subset\mathcal I_k(s)\subset\{r<2r_\delta\}.
\end{equation}
For $j\ne k$,
\begin{equation}\label{wgt:plateau-inclusion}
 \left\{\xi:\left|\xi+\frac{a_{kj}}{\sqrt\delta}\right|
                      \le r_\delta\right\}
 \subset\mathcal P_k(s)
 \subset\{\xi:\nabla p_k(\xi,s)=0\}.
\end{equation}
\end{lemma}
\begin{proof}
By \eqref{wgt:Hessian-bounds} and
\eqref{wgt:interface-monotonicity}, we obtain
\[
 \frac{r_\delta}{16}|r_k-r_\delta|
 \le\left|\log\mathsf A_k^\delta(r_\delta\mathbf e_r,s)
                         -h_\delta(s)\right|
 \lesssim\delta r_\delta^4.
\]
Differentiating \eqref{wgt:interface-definition}, we obtain
\[
 \partial_\theta r_k
 =-\left.
 \frac{\partial_\theta\log\mathsf A_k^\delta(r\mathbf e_r,s)}
      {\partial_r\log\mathsf A_k^\delta(r\mathbf e_r,s)}
 \right|_{r=r_k},
\]
and, using the same bounds, we obtain $|\partial_\theta r_k|\lesssim\delta r_\delta^3$.
By \eqref{wgt:interface-definition}, we also obtain
$\{r<r_\delta/2\}\subset\mathcal I_k(s)\subset\{r<2r_\delta\}$.

Finally, we assume $|\xi+a_{kj}/\sqrt\delta|\le r_\delta$ with $j\ne k$.
Using the gap \eqref{prof:gap}, the diameter bound \eqref{wgt:diameter},
and \eqref{wgt:separation-scales}, we obtain
\[
 11r_\delta
 \le\frac{d_*}{\sqrt\delta}-r_\delta
 \le|\xi|
 \le\frac{|a_{kj}|}{\sqrt\delta}+r_\delta
 <\frac R2+\frac R6<R.
\]
Thus we have
$\{\xi:|\xi+a_{kj}/\sqrt\delta|\le r_\delta\}\subset\mathcal P_k(s)$,
and $\nabla p_k=0$ on $\mathcal P_k(s)$.
\end{proof}

The transport term contains $B_k^\app\cdot\nabla Q_k$.

\begin{lemma}\label{wgt:derivatives}
The function $Q_k=\log p_k$ satisfies, almost everywhere,
\begin{equation}\label{wgt:radial-gradient}
 \nabla Q_k=
 \begin{cases}
 \nabla\log\mathsf A_k^\delta,&\xi\in\mathcal I_k(s),\\
 0,&\xi\in\mathcal P_k(s),\\
 \sqrt\nu\,\xi/r,&\xi\in\mathcal O_k(s).
 \end{cases}
\end{equation}
Consequently, for any vector field $B$,
$B\cdot\nabla Q_k=0$ on $\mathcal P_k(s)$ and
$|B\cdot\nabla Q_k|\le\sqrt\nu\,|B|$ on $\mathcal O_k(s)$.
\end{lemma}
\begin{proof}
Differentiating \eqref{wgt:definition} directly, we obtain the assertions.
\end{proof}

On $\mathcal I_k(s)$, using \eqref{wgt:Hessian-bounds}, we have
$\partial_s Q_k=\partial_s\log\mathsf A_k^\delta$ and
$|\partial_s Q_k|\lesssim(1+\eta)\delta r^2\av r^2$.
On the plateau, we use the sign of the time derivative.

\begin{lemma}\label{wgt:plateau-time}
On $\mathcal P_k(s)$, $\partial_s Q_k=h_\delta'(s)$ and
$\partial_s p_k=h_\delta'(s)p_k$.
The cap height satisfies
\begin{equation}\label{wgt:cap-decay}
 \frac\kappa4(1+\eta)r_\delta^2
       \le-h_\delta'(s)\le\frac\kappa2(1+\eta)r_\delta^2.
\end{equation}
\end{lemma}
\begin{proof}
By direct differentiation, we obtain
\[
 h_\delta'=-\kappa(1+\eta)r_\delta\partial_r\log\mathsf A(r_\delta),
 \qquad \partial_s Q_k=h_\delta',\qquad
 \partial_s p_k=h_\delta'p_k\quad\text{on }\mathcal P_k(s).
\]
We obtain the bound \eqref{wgt:cap-decay} from
\eqref{wgt:Gaussian-log-derivative}.
\end{proof}

On $\mathcal O_k(s)$, the time derivative alone need not be negative.
We combine it with the terms arising from diffusion.

\begin{lemma}\label{wgt:exterior-time}
For sufficiently small $\nu$, uniformly for $0<s\le S_\nu$ and
$\xi\in\mathcal O_k(s)$,
\begin{equation}\label{wgt:exterior-damping}
 \frac\eta4\bigl(\xi\cdot\nabla Q_k-|\nabla Q_k|^2\bigr)
       -\frac12\partial_s Q_k
 \ge\frac\kappa{16}(1+\eta)r_\delta^2.
\end{equation}
\end{lemma}
\begin{proof}
Differentiating $R=L_0\delta^{-\frac{1}{2}}$ and
\eqref{wgt:definition} directly, we obtain
$\partial_s Q_k=h_\delta'+\frac12(1+\eta)\sqrt\nu R$.
Using \eqref{wgt:radial-gradient}, we obtain
\begin{equation}\label{wgt:outer-coefficient}
 \begin{aligned}
 &\frac\eta4\bigl(\xi\cdot\nabla Q_k-|\nabla Q_k|^2\bigr)
       -\frac12\partial_s Q_k\\
 &\qquad=-\frac{h_\delta'}2-\frac{\sqrt\nu R}{4}
       +\frac{\eta\sqrt\nu}{4}(r-R)-\frac{\eta\nu}{4}.
 \end{aligned}
\end{equation}
Since $0<\frac12-2\kappa<1$ and $\delta\le1$, we have
\begin{equation}\label{wgt:outer-scale-ratio}
 \begin{aligned}
 \frac{\sqrt\nu R}{(1+\eta)r_\delta^2}
 &=L_0\nu^{2\kappa}\frac{\eta^{\frac{1}{2}-2\kappa}}{1+\eta}
 \le L_0\nu^{2\kappa},\\
 \frac{\eta\nu}{(1+\eta)r_\delta^2}&\le\nu.
 \end{aligned}
\end{equation}
Combining these bounds with \eqref{wgt:cap-decay} and $r>R$, we obtain
\eqref{wgt:exterior-damping} whenever
$L_0\nu^{2\kappa}+\nu\le\kappa/4$.
\end{proof}

\subsection{Moment and derivative estimates}
Using the bounds on $p_k$ and $\nabla Q_k$, we obtain moment and tail estimates
in terms of $\|f_k\|_{L^2_{p_k}}$.

\begin{lemma}\label{wgt:bounds}
The weight satisfies
\begin{equation}\label{wgt:basicbounds}
 \av\xi^{16}\lesssim p_k\lesssim\mathsf A(r).
\end{equation}
\end{lemma}
\begin{proof}
By \eqref{wgt:comparison}, it suffices to estimate the radial weight
$p_0$. By \eqref{wgt:entropy-definition}, we have
\begin{equation}\label{wgt:Gaussian-weight-comparison}
 \frac{e^{r^2/4}}{\av r^2}\le\mathsf A(r)
       \le8\frac{e^{r^2/4}}{\av r^2}.
\end{equation}

By monotonicity, we have $p_0\le\mathsf A$ on $r\le R$. On $r>R$,
we obtain the same inequality using
$\partial_r\log\mathsf A\ge r/4\ge\sqrt\nu$ and
$p_0(R)\le\mathsf A(R)$. Since
$\mathsf A(r)=\int_0^1e^{tr^2/4}\dd t$, we obtain
\begin{equation}\label{wgt:sharpgaussianupper}
 p_k\lesssim\mathsf A(r)\le e^{r^2/4}.
\end{equation}

For the lower bound, we use \eqref{wgt:Gaussian-weight-comparison} to obtain
\[
 \sup_{r\le r_\delta}\frac{\av r^{16}}{p_0(r)}
 \le\sup_{r\ge0}\av r^{18}e^{-r^2/4}<\infty.
\]
For $r\ge r_\delta$, we use \eqref{wgt:scale-range} to obtain
\begin{equation}\label{wgt:radius-decay-scales}
 \nu^{-\frac{1}{2}}\le\delta^{-\frac{1}{2(1-\gamma)}},\qquad
 R=L_0\delta^{-\frac{1}{2}}\le L_0\delta^{-\frac{1}{2(1-\gamma)}}.
\end{equation}
Consequently, we obtain
\[
 \sup_{r\ge r_\delta}\av r^{16}e^{-\sqrt\nu(r-R)_+}
 \lesssim R^{16}+\nu^{-8}\lesssim\delta^{-\frac{8}{1-\gamma}}.
\]
Using \eqref{wgt:radial-weight} and
\eqref{wgt:Gaussian-weight-comparison}, we conclude that
\[
 \sup_{r\ge r_\delta}\frac{\av r^{16}}{p_0(r)}
 \lesssim\frac{\delta^{-\frac{8}{1-\gamma}}}{\mathsf A(r_\delta)}
 \lesssim\delta^{-\frac{8}{1-\gamma}-2\kappa}e^{-\delta^{-2\kappa}/4}
 \lesssim1.
\]
We have thus proved \eqref{wgt:basicbounds}.
\end{proof}

\begin{lemma}[Bounds for $\nabla Q_k$]\label{wgt:gradient-bounds}
The logarithmic weight in \eqref{wgt:inner-region} satisfies
\begin{equation}\label{wgt:log-gradient-bounds}
 |\nabla Q_k|\lesssim\av\xi,\qquad
 \|\nabla Q_k\|_{L^\infty(\R^2)}\lesssim r_\delta.
\end{equation}
\end{lemma}
\begin{proof}
On $\mathcal I_k(s)$, using \eqref{wgt:radial-gradient},
\eqref{wgt:Hessian-bounds}, and
\eqref{wgt:Gaussian-log-derivative}, we obtain
\begin{equation}\label{wgt:inner-gradient-bound}
 |\nabla Q_k|
 \le|\nabla\log\mathsf A|+C\delta r\av r^2
 \le\frac r2+C\delta r\av r^2\lesssim r.
\end{equation}
Here we have $r<2r_\delta$ and $\delta\av{2r_\delta}^2\lesssim1$.
On $\mathcal P_k(s)$ and $\mathcal O_k(s)$, we have $|\nabla Q_k|=0$
and $|\nabla Q_k|=\sqrt\nu\le1$, respectively, so we obtain both assertions.
\end{proof}

\begin{lemma}[Moment estimates]\label{wgt:moments}
For every fixed $j\ge0$ and $b>0$,
\begin{equation}\label{wgt:inversemoments}
 \int_{\R^2}\frac{\av\xi^j}{p_k^b}\dd\xi\lsim{j,b,\gamma}1.
\end{equation}
Consequently, for every $f_k\in L^2_{p_k}$ and fixed $j\ge0$,
\begin{equation}\label{wgt:functionmoments}
 \|\av\xi^j f_k\|_{L^1}\lsim{j}\|f_k\|_{L^2_{p_k}},\qquad
 \|f_k\|_{L^{\frac{4}{3}}(\R^2)}\lesssim\|f_k\|_{L^2_{p_k}}.
\end{equation}
\end{lemma}
\begin{proof}
By \eqref{wgt:comparison}, it suffices to estimate $p_0^{-b}$.
Integrating \eqref{wgt:radial-weight} radially and using
\eqref{wgt:Gaussian-weight-comparison}, we obtain
\[
 \begin{aligned}
 \int_{r\le r_\delta}\frac{\av r^j}{p_0^b}\dd\xi
 &\lesssim_b\int_0^\infty r\av r^{j+2b}e^{-br^2/4}\dd r
       \lsim{j,b}1,\\
 \int_{r_\delta<r\le R}\frac{\av r^j}{p_0^b}\dd\xi
 &\lsim{j,b}r_\delta^{2b}e^{-br_\delta^2/4}R^{j+2},\\
 \int_{r>R}\frac{\av r^j}{p_0^b}\dd\xi
 &\lsim{j,b}r_\delta^{2b}e^{-br_\delta^2/4}
       \bigl(R^{j+1}\nu^{-\frac{1}{2}}+\nu^{-\frac{j+2}{2}}\bigr).
 \end{aligned}
\]
Thus, using \eqref{wgt:radius-decay-scales}, we obtain
\begin{equation}\label{wgt:outer-inverse-moments}
 \int_{r\ge r_\delta}\frac{\av r^j}{p_0^b}\dd\xi
 \lsim{j,b,\gamma}
 r_\delta^{2b}\delta^{-\frac{j+2}{2(1-\gamma)}}e^{-br_\delta^2/4}.
\end{equation}
The right-hand side is uniformly bounded, so we have
\eqref{wgt:inversemoments}. Applying this bound to the weighted norms, we obtain
\begin{equation}\label{wgt:weighted-holder}
 \begin{aligned}
 \|\av\xi^j f_k\|_{L^1(\R^2)}
 &\le\|f_k\|_{L^2_{p_k}}
       \left(\int_{\R^2}\frac{\av\xi^{2j}}{p_k}\dd\xi\right)^{\frac{1}{2}}
 \lsim{j}\|f_k\|_{L^2_{p_k}},\\
 \|f_k\|_{L^{\frac{4}{3}}(\R^2)}
 &\le\|f_k\|_{L^2_{p_k}}
       \left(\int_{\R^2}p_k^{-2}\dd\xi\right)^{\frac{1}{4}}
 \lesssim\|f_k\|_{L^2_{p_k}}.
 \end{aligned}
\end{equation}
\end{proof}

\begin{lemma}[Tail estimates]\label{wgt:tails}
For every fixed $j\ge0$, $c_*,N>0$, and every $f_k\in L^2_{p_k}$,
\begin{equation}\label{wgt:functiontails}
 \|f_k\|_{L^{\frac{4}{3}}(|\xi|\ge c_*r_\delta)}
 +\|\av\xi^j f_k\|_{L^1(|\xi|\ge c_*r_\delta)}
 \lsim{j,c_*,N}\delta^N\|f_k\|_{L^2_{p_k}}.
\end{equation}
The same estimate holds on $\{|\xi|\ge c_*/\sqrt\delta\}$.
\end{lemma}
\begin{proof}
For $b>0$, we obtain from \eqref{wgt:Gaussian-weight-comparison}
\[
 \int_{c_*r_\delta\le r\le r_\delta}\frac{\av r^j}{p_0^b}\dd\xi
 \lesssim_b\int_{c_*r_\delta}^\infty r\av r^{j+2b}e^{-br^2/4}\dd r
 \lsim{j,b}e^{-bc_*^2r_\delta^2/8}.
\]
Combining this estimate with \eqref{wgt:outer-inverse-moments} and
\eqref{wgt:comparison}, and taking $c_1=\frac b8\min\{1,c_*^2\}$, we obtain
\begin{equation}\label{wgt:tailmoments}
 \int_{|\xi|\ge c_*r_\delta}\frac{\av\xi^j}{p_k^b}\dd\xi
 \lsim{j,b,c_*,\gamma}
 r_\delta^{2b}\delta^{-\frac{j+2}{2(1-\gamma)}}e^{-c_1r_\delta^2}
 \lsim{j,b,c_*,\gamma,N}\delta^N.
\end{equation}
On $\{|\xi|\ge c_*r_\delta\}$, using the estimates in
\eqref{wgt:weighted-holder} together with \eqref{wgt:tailmoments}
at powers $2N$ and $4N$, we obtain \eqref{wgt:functiontails}.
Finally, we have $c_*/\sqrt\delta\ge c_*r_\delta$ because $\delta\le1$ and
$\kappa<1/2$. Using
$\{|\xi|\ge c_*/\sqrt\delta\}\subset\{|\xi|\ge c_*r_\delta\}$,
we obtain the last assertion.
\end{proof}

\subsection{Stream-function and velocity bounds}
Using the mass and first-moment conditions, we obtain $H^2$ control of
$\psi^{f_k}$.

\begin{lemma}\label{wgt:stream-H2}
If $f_k\in L^2_{p_k}$ satisfies
$\int_{\R^2}f_k\dd\xi=0$ and $\int_{\R^2}\xi f_k\dd\xi=0$, then
\begin{equation}\label{wgt:stream-H2-bound}
 \|\psi^{f_k}\|_{H^2(\R^2)}+\|\psi^{f_k}\|_{L^\infty(\R^2)}
       \lesssim\|f_k\|_{L^2_{p_k}}.
\end{equation}
In particular, $\|\nabla\psi^{f_k}\|_{L^q(\R^2)}
\lsim{q}\|f_k\|_{L^2_{p_k}}$ for each fixed $2\le q<\infty$.
\end{lemma}
\begin{proof}
Using the two moment conditions, we obtain
\[
 |\widehat f_k(\zeta)|
 =\left|\int_{\R^2}(e^{-i\zeta\cdot\xi}-1+i\zeta\cdot\xi)f_k(\xi)\dd\xi\right|
 \le\tfrac12|\zeta|^2\|\av\xi^2 f_k\|_{L^1}.
\]
Using $\widehat{\psi^{f_k}}(\zeta)=-|\zeta|^{-2}\widehat f_k(\zeta)$
and splitting the Fourier integral at $|\zeta|=1$, we obtain
\[
 \|\psi^{f_k}\|_{H^2}^2
 \lesssim\|\av\xi^2f_k\|_{L^1}^2+\|f_k\|_{L^2}^2
 \lesssim\|f_k\|_{L^2_{p_k}}^2
\]
by \eqref{wgt:basicbounds} and \eqref{wgt:functionmoments}.
The solution defined in Fourier variables agrees with \eqref{wgt:stream-definition}:
both solve the Poisson equation with source $f_k$ and tend to zero at infinity.
Their difference is therefore harmonic and vanishes. Using the two-dimensional
Sobolev embeddings, we obtain
\[
 \|\psi^{f_k}\|_{L^\infty}+\|\nabla\psi^{f_k}\|_{L^q}
 \lsim{q}\|\psi^{f_k}\|_{H^2}\lesssim\|f_k\|_{L^2_{p_k}},
 \qquad 2\le q<\infty.
\]
\end{proof}

\begin{lemma}[\cite{GW02,DG26}]\label{lem:weighted-BS}
Fix $2<q<\infty$. For every $f\in L^{\frac{2q}{q+2}}(\R^2)$,
\begin{equation}\label{eq:HLS-q}
 \|v^f\|_{L^q}\le\bigl\||K|*|f|\bigr\|_{L^q}
 \lsim{q}\|f\|_{L^{\frac{2q}{q+2}}}.
\end{equation}
If $p\ge c\av\xi^{16}$ for some $c>0$ and $f\in L^2_p$ satisfies
$\int_{\R^2}f\dd\xi=0$, then
\begin{equation}\label{eq:weighted-BS-q}
 \|\av\xi v^f\|_{L^q}\lsim{q,c}\|f\|_{L^2_p}.
\end{equation}
\end{lemma}
\begin{proof}[Proof sketch]
Since $|K(\xi)|=(2\pi|\xi|)^{-1}$, we obtain
\eqref{eq:HLS-q} by the Hardy--Littlewood--Sobolev inequality.
To prove \eqref{eq:weighted-BS-q}, we use $\int f=0$ to subtract $K(\xi)$
from the kernel in the Biot--Savart integral. For $r=|\xi|\ge1$, we split at
$|y|=r/2$. Using $|K(\xi-y)-K(\xi)|\lesssim|y|/r^2$
for $|y|\le r/2$ and $r\le2|y|$ for $|y|>r/2$, we obtain
\[
 r|v^f(\xi)|\lesssim\frac1r\||\xi|f\|_{L^1}
       +\bigl(|K|*(|\cdot|\,|f|)\bigr)(\xi).
\]
Since $r^{-1}\in L^q(\{r\ge1\})$, we use \eqref{eq:HLS-q} to estimate
the convolution term and the velocity on $r<1$.
Using $p\ge c\av\xi^{16}$, we obtain
\[
 \|\av\xi v^f\|_{L^q}
 \lsim{q}\|\av\xi f\|_{L^1}
       +\|\av\xi f\|_{L^{\frac{2q}{q+2}}}
 \lsim{q,c}\|f\|_{L^2_p}.
\]
\end{proof}

\begin{lemma}
\label{lem:nonlinear-velocity}
Fix $2<q<\infty$ and assume that $f_j\in L^2_{p_j}$ satisfies
$\int_{\R^2}f_j\dd\xi=0$ for all $j$. Then
\begin{equation}\label{eq:nonlinear-velocity}
 \big\||\nabla Q_k(\xi)|\,|v^{f_j}(\xi+\frac{a_{kj}}{\delta^{\frac12}})|\big\|_{L^q_\xi}
       \lsim{q}\|f_j\|_{L^2_{p_j}}.
\end{equation}
\end{lemma}
\begin{proof}
For $j=k$, we use $|\nabla Q_k|\lesssim \av\xi$ and the preceding lemma.
For $j\ne k$, using \eqref{wgt:inversemoments} and
\eqref{wgt:tailmoments}, we obtain
\[
 \begin{aligned}
 \|f_j\|_{L^{\frac{2q}{q+2}}}&\lsim{q}\|f_j\|_{L^2_{p_j}},\\
 \|f_j\|_{L^{\frac{2q}{q+2}}(|y|\ge d_*/(4\delta^{\frac12}))}
 &\le\|f_j\|_{L^2_{p_j}}
 \left(\int_{|y|\ge d_*/(4\delta^{\frac12})}p_j^{-\frac{q}{2}}\dd y\right)^{\frac{1}{q}}.
 \end{aligned}
\]
In particular, for every fixed $N>0$, we have
\begin{equation}\label{eq:nonlinear-source-tail-q}
 \|f_j\|_{L^{\frac{2q}{q+2}}(|y|\ge d_*/(4\delta^{\frac12}))}
       \lsim{q,N,d_*}\delta^N\|f_j\|_{L^2_{p_j}}.
\end{equation}
The required smallness of $\nu$ is independent of $q,N$.

On $\mathcal I_k(s)\subset\{r<2r_\delta\}$,
using \eqref{wgt:inner-gradient-bound}, we obtain
$\|\nabla Q_k\|_{L^q(\mathcal I_k(s))}
\lsim{q} r_\delta^{1+\frac{2}{q}}$.
By \eqref{wgt:separation-scales}, for $|\xi|\le3r_\delta$ and
$|y|\le d_*/(4\delta^{\frac12})$, we have
\[
 |\xi+\frac{a_{kj}}{\delta^{\frac12}}-y|
 \ge\frac{|Z_k-Z_j|-3r_\delta\delta^{\frac12}-d_*/4}{\delta^{\frac12}}
 \ge\frac{d_*}{2\delta^{\frac12}}.
\]
We therefore bound the contribution from $\{|y|\le d_*/(4\delta^{\frac12})\}$ by
$C_q\delta^{\frac{1}{2}}r_\delta^{1+\frac{2}{q}}\|f_j\|_{L^2_{p_j}}$.
For the prefactor, we have
\begin{equation}\label{eq:nonlinear-near-source-factor}
 \delta^{\frac{1}{2}}r_\delta^{1+\frac{2}{q}}
 =\delta^{\frac{1}{2}-\kappa(1+\frac{2}{q})}
 \le\delta^{\frac12-2\kappa}\le1,
 \qquad q>2.
\end{equation}
For the integral over $\{|y|\ge d_*/(4\delta^{\frac12})\}$, we use
\eqref{eq:HLS-q}, translation invariance,
\eqref{eq:nonlinear-source-tail-q} with $N=1$, and
$\|\nabla Q_k\|_{L^\infty}\lesssim r_\delta$ to obtain
\begin{equation}\label{eq:nonlinear-far-source-q}
 \left\||\nabla Q_k|\, \left|\int_{|y|\ge d_*/(4\delta^{\frac12})} K(\,\cdot+\frac{a_{kj}}{\delta^{\frac12}}-y)f_j(y)\dd y\right|\right\|_{L^q} \lsim{q} r_\delta\delta\|f_j\|_{L^2_{p_j}} =\delta^{1-\kappa}\|f_j\|_{L^2_{p_j}}.
\end{equation}
Finally, on $\mathcal O_k(s)$, the logarithmic gradient has magnitude
$\sqrt\nu\le1$. By translation invariance and \eqref{eq:HLS-q}, we obtain
\[
 \big\||\nabla Q_k|\,|v^{f_j}(\,\cdot+\frac{a_{kj}}{\delta^{\frac12}})|\big\|_{L^q(\mathcal O_k(s))}
 \le\sqrt\nu\|v^{f_j}\|_{L^q}
 \lsim{q}\sqrt\nu\|f_j\|_{L^2_{p_j}}\lsim{q}\|f_j\|_{L^2_{p_j}}.
\]
The gradient vanishes on $\mathcal P_k(s)$, so we obtain
\eqref{eq:nonlinear-velocity}.
\end{proof}

\subsection{Coercivity of the energy}\label{wgt:coercivity-section}
\begin{lemma}\label{wgt:coercivity}
There is a constant $c>0$, independent of $\nu$, such that every $f_k\in L^2_{p_k}$ satisfying
$\int_{\R^2}f_k\dd\xi=0$ and $\int_{\R^2}\xi f_k\dd\xi=0$ obeys
\[
 c\|f_k\|_{L^2_{p_k}}^2
 \le E_{p_k}(f_k):=\frac12\left(\|f_k\|_{L^2_{p_k}}^2
       +\int_{\R^2}\psi^{f_k} f_k\dd\xi\right)
 \le\frac12\|f_k\|_{L^2_{p_k}}^2.
\]
Consequently, \eqref{wgt:energy-equivalence} holds.
\end{lemma}
\begin{proof}
Using \eqref{wgt:stream-H2-bound} and integrating by parts, we obtain
$\int_{\R^2}\psi^{f_k} f_k\dd\xi=-\|\nabla\psi^{f_k}\|_{L^2}^2\le0$,
and hence we obtain the upper bound. For the lower bound, we set
$f_{k,\mathrm{near}}=\boldsymbol1_{\{r\le r_\delta/2\}}f_k$ and
$f_{k,\mathrm{far}}=\boldsymbol1_{\{r>r_\delta/2\}}f_k$.
Their supports are disjoint, so we have
\begin{equation}\label{wgt:coercivity-norm-splitting}
 \|f_k\|_{L^2_{p_k}}^2
 =\|f_{k,\mathrm{near}}\|_{L^2_{p_k}}^2+\|f_{k,\mathrm{far}}\|_{L^2_{p_k}}^2.
\end{equation}
Using the symmetry of the logarithmic kernel, we then obtain
\begin{equation}\label{wgt:coercivity-splitting}
 \begin{aligned}
 E_{p_k}(f_k)
 &=E_{p_k}(f_{k,\mathrm{near}})+\frac12\|f_{k,\mathrm{far}}\|_{L^2_{p_k}}^2
       +\int_{\R^2}\psi^{f_{k,\mathrm{near}}}f_{k,\mathrm{far}}\dd\xi+\frac12\int_{\R^2}\psi^{f_{k,\mathrm{far}}}f_{k,\mathrm{far}}\dd\xi.
 \end{aligned}
\end{equation}

By monotonicity of $p_0$ and \eqref{wgt:comparison}, we have
$p_k\ge\frac12\mathsf A(r_\delta/2)$ on $r\ge r_\delta/2$.
Thus, using \eqref{wgt:Gaussian-weight-comparison}, we obtain
\[
 \|f_{k,\mathrm{far}}\|_{L^2}
 \lesssim r_\delta e^{-r_\delta^2/32}\|f_k\|_{L^2_{p_k}}
 \lesssim\delta\|f_k\|_{L^2_{p_k}}.
\]
Combining this with \eqref{wgt:functiontails}, we obtain
\begin{equation}\label{wgt:coercivity-tail}
 \|f_{k,\mathrm{far}}\|_{L^2}+\|\av\xi f_{k,\mathrm{far}}\|_{L^1}
 \lesssim\delta\|f_k\|_{L^2_{p_k}}.
\end{equation}
Using \eqref{wgt:basicbounds} and \eqref{wgt:functionmoments}, we obtain
\begin{equation}\label{wgt:coercivity-near-norms}
 \|f_{k,\mathrm{near}}\|_{L^2}+\|\av\xi f_{k,\mathrm{near}}\|_{L^1}
 \lesssim\|f_k\|_{L^2_{p_k}}.
\end{equation}

Using $\log|\cdot|\in L^1(\{|\xi|<1\})$ and
$\log|\xi-y|\le\av\xi+\av y$ for $|\xi-y|\ge1$, together with
\eqref{wgt:coercivity-tail} and
\eqref{wgt:coercivity-near-norms}, we obtain
\begin{equation}\label{wgt:coercivity-pairings}
 \left|\int_{\R^2}\psi^{f_{k,\mathrm{near}}}f_{k,\mathrm{far}}\dd\xi\right|
 +\left|\int_{\R^2}\psi^{f_{k,\mathrm{far}}}f_{k,\mathrm{far}}\dd\xi\right|
 \lesssim\delta\|f_k\|_{L^2_{p_k}}^2.
\end{equation}
Substituting \eqref{wgt:coercivity-pairings} into
\eqref{wgt:coercivity-splitting}, we obtain
\begin{equation}\label{wgt:coercivity-reduction}
 E_{p_k}(f_k)\ge E_{p_k}(f_{k,\mathrm{near}})
       +\frac12\|f_{k,\mathrm{far}}\|_{L^2_{p_k}}^2
       -C\delta\|f_k\|_{L^2_{p_k}}^2.
\end{equation}

On $r\le r_\delta/2$, we have $p_k=\mathsf A_k^\delta$. Hence, by
\eqref{wgt:Hessian-bounds}, we have $f_{k,\mathrm{near}}\in L^2_{\mathsf A}$,
\begin{equation}\label{wgt:coercivity-weight-comparison}
 \frac23\|f_{k,\mathrm{near}}\|_{L^2_{p_k}}^2
 \le\|f_{k,\mathrm{near}}\|_{L^2_{\mathsf A}}^2
 \le2\|f_{k,\mathrm{near}}\|_{L^2_{p_k}}^2,
\end{equation}
and
\begin{equation}\label{wgt:coercivity-weight-error}
 \left|\int_{\R^2}(p_k-\mathsf A)f_{k,\mathrm{near}}^2\dd\xi\right|
 \lesssim\delta r_\delta^4\|f_{k,\mathrm{near}}\|_{L^2_{\mathsf A}}^2
 \lesssim\delta r_\delta^4\|f_k\|_{L^2_{p_k}}^2.
\end{equation}
Using the moment conditions on $f_k$ and \eqref{wgt:coercivity-tail}, we obtain
\begin{equation}\label{wgt:coercivity-moments}
 \left|\int_{\R^2}f_{k,\mathrm{near}}\dd\xi\right|
 +\left|\int_{\R^2}\xi f_{k,\mathrm{near}}\dd\xi\right|
 \lesssim\|\av\xi f_{k,\mathrm{far}}\|_{L^1}
 \lesssim\delta\|f_k\|_{L^2_{p_k}}.
\end{equation}
Using the Gaussian Arnold estimate \cite{GS24,DoGa26} and
\eqref{wgt:coercivity-moments}, we obtain, for some $0<c<1/8$,
\begin{equation}\label{wgt:Gaussian-coercivity}
 \begin{aligned}
 E_{\mathsf A}(f_{k,\mathrm{near}})
 &\ge4c\|f_{k,\mathrm{near}}\|_{L^2_{\mathsf A}}^2
       -C\left(\left|\int_{\R^2}f_{k,\mathrm{near}}\dd\xi\right|^2
       +\left|\int_{\R^2}\xi f_{k,\mathrm{near}}\dd\xi\right|^2\right)\\
 &\ge4c\|f_{k,\mathrm{near}}\|_{L^2_{\mathsf A}}^2
       -C\delta^2\|f_k\|_{L^2_{p_k}}^2.
 \end{aligned}
\end{equation}
Combining \eqref{wgt:Gaussian-coercivity},
\eqref{wgt:coercivity-weight-comparison}, and
\eqref{wgt:coercivity-weight-error}, we obtain
\begin{equation}\label{wgt:coercivity-perturbation}
 \begin{aligned}
 E_{p_k}(f_{k,\mathrm{near}})
 &=E_{\mathsf A}(f_{k,\mathrm{near}})+\frac12\int_{\R^2}(p_k-\mathsf A)f_{k,\mathrm{near}}^2\dd\xi\\
 &\ge2c\|f_{k,\mathrm{near}}\|_{L^2_{p_k}}^2
       -C(\delta r_\delta^4+\delta^2)\|f_k\|_{L^2_{p_k}}^2.
 \end{aligned}
\end{equation}
By \eqref{wgt:coercivity-reduction},
\eqref{wgt:coercivity-perturbation},
\eqref{wgt:coercivity-norm-splitting}, and \eqref{wgt:scales}, we obtain
\[
 \begin{aligned}
 E_{p_k}(f_k)&\ge2c\|f_{k,\mathrm{near}}\|_{L^2_{p_k}}^2
       +\frac12\|f_{k,\mathrm{far}}\|_{L^2_{p_k}}^2
       -C(\delta r_\delta^4+\delta)\|f_k\|_{L^2_{p_k}}^2\\
 &\ge\left(2c-C(\delta^{1-4\kappa}+\delta)\right)
       \|f_k\|_{L^2_{p_k}}^2.
 \end{aligned}
\]
By \eqref{wgt:scale-range}, we obtain the lower bound for sufficiently
small $\nu$. Summing over $k$ in \eqref{lin:energy-def}, we obtain
\eqref{wgt:energy-equivalence}.
\end{proof}

\begin{proof}[Proof of Proposition~\ref{wgt:main}]
Using Lemmas~\ref{wgt:regularity} and \ref{wgt:bounds}, we obtain the stated regularity
and pointwise bounds for the weights defined in \eqref{wgt:definition}.
By Lemma~\ref{wgt:coercivity}, we have \eqref{wgt:energy-equivalence}.
Combining \eqref{wgt:functionmoments}, Lemma~\ref{wgt:stream-H2},
and Lemma~\ref{lem:nonlinear-velocity} with this equivalence, we obtain
\eqref{wgt:main-controls}.
\end{proof}

\section{Center dynamics and moment cancellations}\label{ctr:section}
We assume \eqref{wgt:parameter-range}, \eqref{prof:gap}, and
\eqref{wgt:diameter}, and use the profiles \eqref{prof:approx},
the weights \eqref{wgt:definition}, and the energy \eqref{lin:energy-def}.
Constants are uniform in $s$, $\nu$, and $\delta$.

\begin{proposition}\label{ctr:energy-force}
Assume that $w_k=w_{\app}^k+\delta f_k$, with
$f_k\in L^2(\R^2,p_k\dd\xi)$,
$\int_{\R^2}f_k\dd\xi=0$, and
$\int_{\R^2}\xi f_k\dd\xi=0$.
The exact viscous centers satisfy
\begin{equation}\label{ctr:refined-center-equation}
 Z'=\mathcal F(Z)+\delta^2V_2^\nu(s)+e,\qquad
 |e|\lesssim\delta^3+\delta^2\sqrt{\mathcal E}+\delta^4\mathcal E,
\end{equation}
where $V_2^\nu$ is given by \eqref{prof:center-coefficient}.
\end{proposition}

\subsection{The exact center law and pair interactions}
For $a\in\R^2$ and $\delta>0$, we write
\begin{equation}\label{ctr:pair-definition}
 \mathcal I_{a,\delta}(f,g)
 :=\iint_{\R^2\times\R^2}
 K(a+\delta^{\frac12}(\xi-y))f(\xi)g(y)\dd\xi\dd y.
\end{equation}
\begin{lemma}\label{ctr:law}
For $s>0$, the exact centers satisfy
\begin{equation}\label{ctr:exact-law}
 Z_k'(s)=\frac12Z_k(s)-\beta Z_k(s)^\perp+
 \sum_{j\ne k}\Gamma_j
 \mathcal I_{a_{kj}(s),\delta(s)}(w_k(s),w_j(s)).
\end{equation}
\end{lemma}
\begin{proof}
In collapse coordinates, $W_k$ satisfies
\[
 \partial_sW_k+\nabla_X\cdot
 \bigl[(V+X/2-\beta X^\perp)W_k\bigr]=\nu\Delta_XW_k,
 \qquad V=\sum_j K*W_j.
\]
Taking the first moment, we obtain
$Z_k'=\frac12Z_k-\beta Z_k^\perp
       +\frac1{\Gamma_k}\int_{\R^2}VW_k\dd X$.
The self-interaction vanishes by oddness of $K$:
\[
 \iint_{\R^2\times\R^2}K(X-Y)W_k(s,X)W_k(s,Y)\dd X\dd Y=0.
\]
In the remaining terms, we change variables to
$X=Z_k+\sqrt\delta\,\xi$ and $Y=Z_j+\sqrt\delta\,y$.
Using $W_k(s,X)\dd X=\Gamma_kw_k(s,\xi)\dd\xi$, we obtain \eqref{ctr:exact-law}.
\end{proof}

\begin{lemma}\label{ctr:pair-expansion}
Take an integer $N\ge0$ and assume that $|a|\ge d_*>0$ and $0<\delta\le1$.
If $f,g\in L^{\frac{4}{3}}(\R^2)$ and
$\av\xi^{N+1}f,\av\xi^{N+1}g\in L^1(\R^2)$, then
\begin{equation}\label{ctr:pair-remainder}
 \begin{aligned}
 &\left|\begin{aligned}
 &\mathcal I_{a,\delta}(f,g)\\
 &\quad-\sum_{|\alpha|+|\beta|\le N}
 \frac{(-1)^{|\beta|}\delta^{\frac{|\alpha|+|\beta|}{2}}}{\alpha!\beta!}
 \,\partial^{\alpha+\beta}K(a)
 \int_{\R^2}\xi^\alpha f(\xi)\dd\xi
 \int_{\R^2}y^\beta g(y)\dd y
 \end{aligned}\right|\\
 &\quad\lsim{N,d_*}
  \delta^{\frac{N+1}{2}}\bigl(
      \|\av\xi^{N+1}f\|_{L^1}\|g\|_{L^1}
     +\|f\|_{L^1}\|\av\xi^{N+1}g\|_{L^1}\bigr)\\
 &\qquad+\|\av\xi^N f\|_{L^1(|\xi|\ge d_*/(8\delta^{\frac12}))}
                         \|\av\xi^N g\|_{L^1}\\
 &\qquad+\|\av\xi^N f\|_{L^1}
             \|\av\xi^N g\|_{L^1(|\xi|\ge d_*/(8\delta^{\frac12}))}\\
 &\qquad+\delta^{-\frac{1}{2}}
       \|f\|_{L^{\frac{4}{3}}(|\xi|\ge d_*/(8\delta^{\frac12}))}\|g\|_{L^{\frac{4}{3}}}\\
 &\qquad+\delta^{-\frac{1}{2}}\|f\|_{L^{\frac{4}{3}}}
                    \|g\|_{L^{\frac{4}{3}}(|\xi|\ge d_*/(8\delta^{\frac12}))}.
 \end{aligned}
\end{equation}
\end{lemma}

\begin{proof}
By a direct calculation, we obtain
\[
 \frac1{n!}D^nK(a)[(\xi-y)^{\otimes n}]
 =\sum_{|\alpha|+|\beta|=n}
 \frac{(-1)^{|\beta|}}{\alpha!\beta!}
 \partial^{\alpha+\beta}K(a)\,\xi^\alpha y^\beta.
\]
Thus the finite sum in \eqref{ctr:pair-remainder} is obtained by integrating
the Taylor polynomial of $K(a+\delta^{\frac12}(\xi-y))$.
On $\max\{|\xi|,|y|\}\le d_*/(8\delta^{\frac12})$, the Taylor
segment stays at distance at least $3d_*/4$ from the origin, so we obtain
\[
 \left|K(a+\delta^{\frac12}(\xi-y))
 -\sum_{n=0}^N\frac{\delta^{\frac{n}{2}}}{n!}
 D^nK(a)[(\xi-y)^{\otimes n}]\right|
 \lsim{N,d_*}\delta^{\frac{N+1}{2}}(|\xi|+|y|)^{N+1}.
\]
Integrating against $|f(\xi)g(y)|$, we obtain the first term on the
right of \eqref{ctr:pair-remainder}. For
$\max\{|\xi|,|y|\}>d_*/(8\delta^{\frac12})$, we have
\[
 \left|\sum_{n=0}^N\frac{\delta^{\frac{n}{2}}}{n!}
 D^nK(a)[(\xi-y)^{\otimes n}]\right|
 \lsim{N,d_*}\av\xi^N\av y^N.
\]
Splitting $\{(\xi,y):\max\{|\xi|,|y|\}>d_*/(8\delta^{\frac12})\}$ according to
$|\xi|\ge d_*/(8\delta^{\frac12})$ or
$|y|\ge d_*/(8\delta^{\frac12})$, we obtain the two polynomial tail
terms. Finally, using homogeneity and the Hardy--Littlewood--Sobolev
inequality, we obtain
\begin{equation}\label{ctr:pair-HLS}
 \iint_{\R^2\times\R^2}
 |K(a+\delta^{\frac12}(\xi-y))|\,|h(\xi)\ell(y)|\dd\xi\dd y
 \lesssim\delta^{-\frac{1}{2}}
 \|h\|_{L^{\frac{4}{3}}}\|\ell\|_{L^{\frac{4}{3}}}.
\end{equation}
Applying \eqref{ctr:pair-HLS} first with $f$ restricted to
$\{|\xi|\ge d_*/(8\delta^{\frac12})\}$ and then with $g$ restricted to
$\{|y|\ge d_*/(8\delta^{\frac12})\}$, we obtain the remaining terms.
\end{proof}

\begin{lemma}\label{ctr:pair-identities}
The interaction form satisfies
\begin{equation}\label{ctr:pair-exchange}
 \mathcal I_{a,\delta}(g,f)=-\mathcal I_{-a,\delta}(f,g).
\end{equation}
For $G(\xi)=\frac1{4\pi}e^{-|\xi|^2/4}$ and $a\ne0$,
\begin{equation}\label{ctr:Gaussian-pair}
 \mathcal I_{a,\delta}(G,G)
 =K(a)\left(1-e^{-|a|^2/(8\delta)}\right).
\end{equation}
\end{lemma}
\begin{proof}
Exchanging $\xi$ and $y$ and using the oddness of $K$, we obtain
\eqref{ctr:pair-exchange}. Using
$\int_{\R^2}G(z+y)G(y)\dd y=\frac1{8\pi}e^{-|z|^2/8}$,
we change variables to $z=\xi-y$ and apply the radial Biot--Savart formula
to obtain \eqref{ctr:Gaussian-pair}.
\end{proof}

\subsection{The quadrupole correction}
\begin{lemma}\label{ctr:corrected-approximate}
The coefficient \eqref{prof:center-coefficient}, with
$a_{ij}=Z_i-Z_j$, is given by
\begin{equation}\label{ctr:quadrupole-center-field}
 V_{2,i}^\nu(s)
 =\frac12\sum_{j\ne i}\Gamma_j\int_{\R^2}
 D^2K(a_{ij})[\xi,\xi]
       \bigl(F_i^\nu(\xi,s)+F_j^\nu(\xi,s)\bigr)\dd\xi.
\end{equation}
The interaction of the approximate profiles \eqref{prof:approx}
satisfies
\begin{equation}\label{ctr:corrected-approximate-force}
 \left|\sum_{j\ne i}\Gamma_j
 \bigl\{\mathcal I_{a_{ij},\delta}(w_{\app}^i,w_{\app}^j)
                                      -K(a_{ij})\bigr\}
              -\delta^2V_{2,i}^\nu(s)\right|
 \lesssim\delta^3.
\end{equation}
Thus substitution of these profiles into
\eqref{ctr:exact-law} gives the corrected velocity
$\mathcal F_i(Z)+\delta^2V_{2,i}^\nu$
up to an error $O(\delta^3)$.
\end{lemma}
\begin{proof}
We fix $s$ and $i\ne j$, and set $a=a_{ij}(s)$; all profiles are evaluated at $s$.
We expand the kernel and use angular orthogonality \cite{DoGa26,DG26}.
For $k=i,j$, we have
\[
 w_{\app}^k-G
 =\delta\underbrace{F_k^\nu}_{\text{mode }2}
  +\delta^{\frac{3}{2}}\underbrace{H_k^\nu}_{\text{mode }3}
  +\delta^2\underbrace{P_k^\nu}_{\text{modes }2,4}
  +\delta^{\frac{5}{2}}\underbrace{J_k^\nu}_{\text{modes }1,3,5}.
\]
For $q\in\{G,F_k^\nu,H_k^\nu,P_k^\nu,J_k^\nu\}$,
we apply Proposition~\ref{prof:localization} to obtain
\[
 \begin{aligned}
 \|\av\xi^4q\|_{L^1}+\|q\|_{L^{\frac{4}{3}}}&\lesssim1,\\
 \|\av\xi^4q\|_{L^1(|\xi|\ge d_*/(8\delta^{\frac12}))}
 +\|q\|_{L^{\frac{4}{3}}(|\xi|\ge d_*/(8\delta^{\frac12}))}
 &\lesssim e^{-c/\delta}.
 \end{aligned}
\]
Applying Lemma~\ref{ctr:pair-expansion}, we obtain, for any two such
profiles $q,g$ and $N=1,2,3$,
\begin{equation}\label{ctr:localized-pair-Taylor}
 \mathcal I_{a,\delta}(q,g)
 =\sum_{n=0}^N\frac{\delta^{\frac{n}{2}}}{n!}
 \iint_{\R^2\times\R^2}D^nK(a)[(\xi-y)^{\otimes n}]
 q(\xi)g(y)\dd\xi\dd y+O(\delta^{\frac{N+1}{2}}).
\end{equation}
The mode-two profiles $F_k^\nu$ have zero mass and first moments.
Since $G$ has unit mass and zero first moment, we obtain the quadratic term by direct calculation:
\begin{equation}\label{ctr:quadratic-pair-moment}
 \iint_{\R^2\times\R^2}D^2K(a)[\xi-y,\xi-y]
 F_i^\nu(\xi)G(y)\dd\xi\dd y
 =\int_{\R^2}D^2K(a)[\xi,\xi]F_i^\nu(\xi)\dd\xi.
\end{equation}
The cubic term vanishes because $F_i^\nu$ and $G$ are even.
Exchanging the variables leaves the quadratic term unchanged.
Thus, using \eqref{ctr:localized-pair-Taylor} with $N=3$, we obtain
\begin{equation}\label{ctr:linear-quadrupole-expansion}
 \begin{aligned}
 \mathcal I_{a,\delta}(F_i^\nu,G)
 &=\frac\delta2\int_{\R^2}D^2K(a)[\xi,\xi]F_i^\nu(\xi)\dd\xi
       +O(\delta^2),\\
 \mathcal I_{a,\delta}(G,F_j^\nu)
 &=\frac\delta2\int_{\R^2}D^2K(a)[\xi,\xi]F_j^\nu(\xi)\dd\xi
       +O(\delta^2).
 \end{aligned}
\end{equation}

By angular orthogonality, we find that the moments through degree two of
$H_k^\nu$ vanish, as do the mass and first moments of $P_k^\nu$.
The mass and second moments of $J_k^\nu$ vanish by parity, while
its first moment vanishes because $\Pi_1J_k^\nu\in Y_1^{\mathrm c}$.
Applying \eqref{ctr:localized-pair-Taylor}
with $N=2,1,2$, respectively, we obtain
\[
 \begin{aligned}
 |\mathcal I_{a,\delta}(H_i^\nu,G)|
  +|\mathcal I_{a,\delta}(G,H_j^\nu)|&\lesssim\delta^{\frac{3}{2}},\\
 |\mathcal I_{a,\delta}(P_i^\nu,G)|
  +|\mathcal I_{a,\delta}(G,P_j^\nu)|&\lesssim\delta,\\
 |\mathcal I_{a,\delta}(J_i^\nu,G)|
  +|\mathcal I_{a,\delta}(G,J_j^\nu)|&\lesssim\delta^{\frac{3}{2}}.
 \end{aligned}
\]
Every corrector has zero mass and first moments. For any pair
$q_i,q_j$ of correctors, each Taylor monomial of total degree at
most three contains a mass or first moment of one of the correctors,
and therefore vanishes. Using \eqref{ctr:localized-pair-Taylor} with $N=3$,
we obtain $|\mathcal I_{a,\delta}(q_i,q_j)|\lesssim\delta^2$.
Since each corrector enters $w_{\app}^k-G$ with a coefficient
$O(\delta)$, by bilinearity we obtain
$|\mathcal I_{a,\delta}(w_{\app}^i-G,w_{\app}^j-G)|\lesssim\delta^4$.
Combining these estimates with \eqref{ctr:Gaussian-pair}, we obtain
\begin{equation}\label{ctr:approximate-pair-expansion}
 \mathcal I_{a,\delta}(w_{\app}^i,w_{\app}^j)-K(a)
 =\frac{\delta^2}{2}\int_{\R^2}D^2K(a)[\xi,\xi]
 \bigl(F_i^\nu(\xi,s)+F_j^\nu(\xi,s)\bigr)\dd\xi+O(\delta^3).
\end{equation}
Multiplying by the fixed circulations, summing, and using
\eqref{ctr:quadrupole-center-field}, we obtain
\eqref{ctr:corrected-approximate-force}.
\end{proof}

\subsection{Moment cancellation in the center equation}
\begin{lemma}
\label{ctr:pair-cancellation}
Assume that $|a|\ge d_*>0$ and let $p_i,p_j$ be the weights
\eqref{wgt:definition} under \eqref{wgt:parameter-range}, for
sufficiently small $\nu$. Assume that $f\in L^2_{p_i}$ satisfies
$\int_{\R^2}f(\xi)\dd\xi
 =\int_{\R^2}\xi f(\xi)\dd\xi=0$.
If $g$ has uniform Gaussian localization, then
\begin{equation}\label{ctr:single-cancellation}
 |\mathcal I_{a,\delta}(f,g)|+
 |\mathcal I_{a,\delta}(g,f)|
 \lesssim\delta\|f\|_{L^2_{p_i}}.
\end{equation}
If $g\in L^2_{p_j}$ satisfies
$\int_{\R^2}g(y)\dd y
 =\int_{\R^2}y g(y)\dd y=0$,
then
\begin{equation}\label{ctr:double-cancellation}
 |\mathcal I_{a,\delta}(f,g)|
 \lesssim\delta^2\|f\|_{L^2_{p_i}}\|g\|_{L^2_{p_j}}.
\end{equation}
The constants depend on $d_*$, $\gamma$, and the fixed weight parameters;
for \eqref{ctr:single-cancellation}, they also depend on the Gaussian bounds for $g$.
They are independent of $s$, $\delta$, and $\nu$.
\end{lemma}
\begin{proof}
Using the inverse-moment bounds \eqref{wgt:inversemoments}, we have
$\int_{\R^2}\frac{\av{\xi}^{8}}{p_k}\dd\xi+
 \int_{\R^2}\frac1{p_k^2}\dd\xi\lesssim1$ for $k=i,j$,
and, for every $N>0$, we obtain
\begin{equation}\label{ctr:dual-tails}
 \begin{aligned}
 &\int_{|\xi|\ge d_*/(8\delta^{\frac12})}
 \left(\frac{\av\xi^8}{p_k}+\frac1{p_k^2}\right)\dd\xi\\
 &\quad\lsim{N,d_*}\delta^N\int_{\R^2}
 \left(\frac{\av\xi^{8+2N}}{p_k}+\frac{\av\xi^{2N}}{p_k^2}\right)\dd\xi
 \lsim{N}\delta^N,\qquad k=i,j.
 \end{aligned}
\end{equation}
Consequently, we have
$\|\av\xi^4f\|_{L^1}+\|f\|_{L^{\frac{4}{3}}}
 \lesssim\|f\|_{L^2_{p_i}}$.
Applying \eqref{ctr:dual-tails} with $N=12$, we also obtain
\begin{equation}\label{ctr:explicit-tail-norms}
 \|\av\xi^4f\|_{L^1(|\xi|\ge d_*/(8\delta^{\frac12}))}
 \lesssim\delta^6\|f\|_{L^2_{p_i}},\qquad
 \|f\|_{L^{\frac{4}{3}}(|\xi|\ge d_*/(8\delta^{\frac12}))}
 \lesssim\delta^3\|f\|_{L^2_{p_i}}.
\end{equation}
If $g$ has uniform Gaussian localization, the norms
$\|\av\xi^2g\|_{L^1}$ and $\|g\|_{L^{\frac{4}{3}}}$ are bounded,
and their tails on
$|\xi|\ge d_*/(8\delta^{\frac12})$ are $O(e^{-c/\delta})$.
The zero mass and first moments of $f$ annihilate the constant
and linear Taylor terms. Thus, using \eqref{ctr:pair-remainder} with
$N=1$, we obtain
\[
 |\mathcal I_{a,\delta}(f,g)|
 \lesssim\bigl(\delta+\delta^6+\delta^{\frac{5}{2}}\bigr)
 \|f\|_{L^2_{p_i}}\lesssim\delta\|f\|_{L^2_{p_i}}.
\]
Using \eqref{ctr:pair-exchange}, we obtain the same estimate with the arguments reversed,
and hence \eqref{ctr:single-cancellation}.

If $g\in L^2_{p_j}$ has zero mass and first moments, the same
norm and tail bounds hold for $g$. Every Taylor monomial of total
degree at most three contains a moment of degree at most one of
$f$ or $g$ and therefore vanishes. Applying
\eqref{ctr:pair-remainder} with $N=3$, we obtain
\[
 |\mathcal I_{a,\delta}(f,g)|
 \lesssim\bigl(\delta^2+\delta^6+\delta^{\frac{5}{2}}\bigr)
 \|f\|_{L^2_{p_i}}\|g\|_{L^2_{p_j}}
 \lesssim\delta^2\|f\|_{L^2_{p_i}}\|g\|_{L^2_{p_j}}.\qedhere
\]
\end{proof}

\begin{proof}[Proof of Proposition~\ref{ctr:energy-force}]
For $a=a_{ij}$, we use bilinearity to write
\begin{equation}\label{ctr:exact-pair-decomposition}
 \begin{aligned}
 \mathcal I_{a,\delta}(w_i,w_j)-K(a)
 ={}&\mathcal I_{a,\delta}(w_{\app}^i,w_{\app}^j)-K(a)
 +\delta\mathcal I_{a,\delta}(f_i,w_{\app}^j)\\
 &+\delta\mathcal I_{a,\delta}(w_{\app}^i,f_j)
 +\delta^2\mathcal I_{a,\delta}(f_i,f_j).
 \end{aligned}
\end{equation}
Substituting the estimates of Lemmas~\ref{ctr:corrected-approximate}
and \ref{ctr:pair-cancellation} into \eqref{ctr:exact-law}, and
using coercivity \eqref{wgt:energy-equivalence}, we obtain
\[
 \begin{aligned}
 |Z'-\mathcal F(Z)-\delta^2V_2^\nu(s)|
 &\lesssim\delta^3+\delta^2\sum_i\|f_i\|_{L^2_{p_i}}
 +\delta^4\sum_{i,j}\|f_i\|_{L^2_{p_i}}\|f_j\|_{L^2_{p_j}}\\
 &\lesssim\delta^3+\delta^2\sqrt{\mathcal E}+\delta^4\mathcal E.
 \end{aligned}
\]
We have proved \eqref{ctr:refined-center-equation}. Since $V_2^\nu$
is uniformly bounded by Proposition~\ref{prof:localization},
$\mathfrak r=\delta^2V_2^\nu+e$ satisfies
\begin{equation}\label{ctr:closed-force}
 Z'=\mathcal F(Z)+\mathfrak r,\qquad
 |\mathfrak r|\lesssim\delta^2
       \bigl(1+\sqrt{\mathcal E}+\delta^2\mathcal E\bigr).
\end{equation}
\end{proof}

\begin{proposition}
\label{ctr:continuation}
Fix the collapsing configuration $z^0$ from
\eqref{eq:configuration}. There is $\rho>0$, depending only on this
configuration, with the following property. Assume that $Z$ solves
$Z'=\mathcal F(Z)+\mathfrak r$, with
$|\mathfrak r(s)|\lesssim\delta(s)^2$ and
$\delta(s)=\nu(e^s-1)$,
on $[s_0,S)$, where $s_0=\log2$, $s_0<S\le S_\nu$, and
$\mathfrak r$ is continuous. Assume also that, for some $M\ge1$,
$|Z(s_0)-z^0|\lesssim\nu^2$ and
$|Z(s)-z^0|\le\min\{\rho,M\delta(s)^2\}$,
with the implicit constants independent of $M$, $s$, and $\nu$.
Then
\begin{equation}\label{ctr:center-bootstrap-estimate}
 |Z(s)-z^0|\le C\bigl(\delta(s)^2+M^2\delta(s)^4\bigr),
 \qquad s_0\le s<S,
\end{equation}
where $C$ is independent of $M$, $s$, and $\nu$.
\end{proposition}
\begin{proof}
We set $h=Z-z^0$. By \eqref{geo:nonlinear-linearization}, we choose
$\rho>0$ depending only on the fixed configuration such that
$4\rho\le\min_{i\ne j}|z_i^0-z_j^0|$ and
$|\mathcal N(h)|\lesssim|h|^2$ for $|h|\le\rho$.
Then we have $\min_{i\ne j}|Z_i-Z_j|\ge2\rho$ and
$h'=\mathscr Ah+\mathcal N(h)+\mathfrak r$, from which we obtain
\begin{equation}\label{ctr:center-duhamel}
 h(s)=e^{(s-s_0)\mathscr A}h(s_0)
 +\int_{s_0}^s e^{(s-s')\mathscr A}
       \bigl(\mathfrak r(s')+\mathcal N(h(s'))\bigr)\dd s'.
\end{equation}
For $s_0\le s'\le s$, we have
$\frac{\delta(s')}{\delta(s)}
 \le e^{-(s-s')}$.
Combining this ratio bound with \eqref{geo:semigroup}, we obtain, for $p>1$,
\begin{equation}\label{ctr:forced-convolution}
 \int_{s_0}^s\|e^{(s-s')\mathscr A}\|\delta(s')^p\dd s'
 \lesssim\delta(s)^p\int_0^\infty(1+u)e^{-(p-1)u}\dd u
 =\left(\frac1{p-1}+\frac1{(p-1)^2}\right)\delta(s)^p.
\end{equation}
Since $\delta(s_0)=\nu$, we obtain
\[
 |e^{(s-s_0)\mathscr A}h(s_0)|
 \lesssim(1+s-s_0)e^{-(s-s_0)}\delta(s)^2
 \lesssim\delta(s)^2.
\]
Finally, we have $|\mathcal N(h(s'))|\lesssim M^2\delta(s')^4$.
Applying \eqref{ctr:forced-convolution} with $p=2,4$ in
\eqref{ctr:center-duhamel}, we obtain \eqref{ctr:center-bootstrap-estimate},
with $C$ depending only on the fixed configuration and the
implicit constants in the initial and forcing bounds.
\end{proof}

\section{The weighted energy estimate}\label{sec:nonlinear}

The remainders $f_k=(w_k-w_{\app}^k)/\delta$ satisfy
\eqref{eq:error-equation}. With the weights from
Section~\ref{wgt:section}, we recall the energy \eqref{lin:energy-def} and the dissipation
\eqref{lin:dissipation-def}:
\[
 \mathcal E(s)=\frac12\sum_{k=1}^3\int_{\R^2}
       \bigl(p_kf_k^2+\psi^{f_k}f_k\bigr)\dd\xi,
 \qquad
 \mathcal D(s)=\sum_{k=1}^3\|\nabla(\sqrt{p_k}f_k)\|_{L^2}^2.
\]

\begin{proposition}\label{prop:energy}
Assume that
\eqref{prof:configuration-neighborhood} and
\eqref{prof:center-velocity-bounds} hold on an interval
$I\subset(0,S_\nu]$, where $S_\nu=\gamma\log(1/\nu)$.
For sufficiently small $\nu$, there are positive constants
$c,C$ such that, for every fixed $2<q<\infty$,
\begin{equation}\label{eq:complete-energy}
 \begin{aligned}
 &\mathcal E'+c\eta(\mathcal D+\mathcal E)+\mathcal E
       +c(1+\eta)r_\delta^2
          \sum_k\int_{\mathcal O_k(s)}p_kf_k^2\dd\xi\\
 &\qquad\le C\delta(\delta+\nu)+C_q\eta^{-\frac1{q-1}}
       \mathcal E^{\frac{3q-2}{2(q-1)}}+C\mathcal E^{\frac32}.
 \end{aligned}
\end{equation}
\end{proposition}

\subsection{The energy identity}
\begin{proposition}\label{prop:energy-identity}
For the remainders satisfying \eqref{eq:error-equation}, the energy
\eqref{lin:energy-def} satisfies
\begin{equation}\label{eq:energy-five-terms}
 \mathcal{E}'=R_1+R_2+R_3+R_4+R_5,
\end{equation}
where
\begin{align}
 R_1={}&\sum_k\int_{\R^2}(p_kf_k+\psi^{f_k})
                 \{\eta\Lop f_k-(1+\eta)f_k\}\dd\xi
       +\frac12\sum_k\int_{\R^2}(\partial_s p_k)f_k^2\dd\xi,
       \label{energy:R1}\\
 R_2={}&-\sum_k\int_{\R^2}(p_kf_k+\psi^{f_k})
       \left(B_k^\app\cdot\nabla f_k
       +\frac{\Gamma_k}{\delta}v^{f_k}\cdot\nabla w_{\app}^k\right)\dd\xi,
       \label{energy:R2}\\
 R_3={}&-\sum_k\sum_{j\ne k}\frac{\Gamma_j}{\delta}
       \int_{\R^2}(p_kf_k+\psi^{f_k})
           v^{f_j}(\xi+\frac{a_{kj}}{\delta^{\frac12}})\cdot\nabla w_{\app}^k\dd\xi,
       \label{energy:R3}\\
 R_4={}&-\sum_{j,k}\Gamma_j
       \int_{\R^2}(p_kf_k+\psi^{f_k})
           v^{f_j}(\xi+\frac{a_{kj}}{\delta^{\frac12}})\cdot\nabla f_k\dd\xi,
       \label{energy:R4}\\
 R_5={}&-\delta^{-1}\sum_k\int_{\R^2}
                    (p_kf_k+\psi^{f_k})\left\{\mathcal R_k
       -\frac{e_k}{\delta^{\frac12}}\cdot\nabla w_{\app}^k\right\}\dd\xi.
       \label{energy:R5}
\end{align}
\end{proposition}
\begin{proof}
By Lemma~\ref{wgt:regularity} and the symmetry of the logarithmic
kernel in \eqref{wgt:stream-definition}, we differentiate directly to obtain
\begin{equation}\label{energy:time-derivative}
 \mathcal E'=\sum_k\int_{\R^2}(p_kf_k+\psi^{f_k})\partial_sf_k\dd\xi
                +\frac12\sum_k\int_{\R^2}(\partial_s p_k)f_k^2\dd\xi.
\end{equation}
Substituting \eqref{eq:error-equation}, we obtain
\eqref{eq:energy-five-terms}.
\end{proof}

\subsection{Diffusion and damping}
For the estimates below, we assume
\eqref{prof:configuration-neighborhood} and \eqref{prof:center-velocity-bounds}.

\begin{lemma}\label{lin:diffusion}
There is a fixed constant $c>0$, which may depend on $\kappa$,
such that $R_1$ defined in \eqref{energy:R1} satisfies
\begin{equation}\label{energy:R1-bound}
 R_1\le{}-c\eta(\mathcal D+\mathcal E)-2(1+\eta)\mathcal E+C\delta^{1-4\kappa}\mathcal E-c(1+\eta)r_\delta^2\sum_k\int_{\mathcal O_k(s)}p_kf_k^2\dd\xi.
\end{equation}
\end{lemma}
\begin{proof}
Using $\Delta\psi^{f_k}=f_k$ and \eqref{eq:error-zero-moments},
we integrate by parts to obtain
$\int_{\R^2}\psi^{f_k}\Lop f_k\dd\xi
 =\|f_k\|_{L^2}^2$.
Using $Q_k=\log p_k$, we obtain
\begin{equation}\label{lin:conjugation}
 \begin{aligned}
 \mathcal Q_{p_k}(f_k)
 &:=-\int_{\R^2}(p_kf_k+\psi^{f_k})\Lop f_k\dd\xi\\
 &=\|\nabla(\sqrt{p_k}f_k)\|_{L^2}^2+\int_{\R^2}
 \left\{\frac14(\xi\cdot\nabla Q_k-|\nabla Q_k|^2)
                        -\frac12-\frac1{p_k}\right\}p_kf_k^2\dd\xi.
 \end{aligned}
\end{equation}
Thus, using \eqref{energy:R1} and \eqref{lin:energy-def}, we obtain
\begin{equation}\label{lin:diffusionbound}
 R_1=-\sum_k\left\{\eta\mathcal Q_{p_k}(f_k)
                    -\frac12\int_{\R^2}(\partial_s p_k)f_k^2\dd\xi\right\}
                         -2(1+\eta)\mathcal E.
\end{equation}

We compare the coefficient in \eqref{lin:conjugation} with
\[
 \mathsf V_{\mathsf A}
 :=\frac14\bigl(\xi\cdot\nabla\log\mathsf A
                    -|\nabla\log\mathsf A|^2\bigr)
       -\frac12-\frac1{\mathsf A}.
\]
Using $\mathsf A\ge1$ and the bounds
$r/4\le\partial_r\log\mathsf A\le r/2$ from
\eqref{wgt:Gaussian-log-derivative}, we obtain
\begin{equation}\label{lin:Gaussian-potential}
 \mathsf V_{\mathsf A}\ge\frac{r^2}{32}-\frac32.
\end{equation}
On $\mathcal I_k(s)$, using \eqref{wgt:Hessian-bounds} and
\eqref{wgt:interface-location}, we obtain
\begin{equation}\label{lin:diffusion-uncapped}
 \frac\eta4(\xi\cdot\nabla Q_k-|\nabla Q_k|^2)
       -\frac12\partial_s Q_k-\frac\eta2-\frac\eta{p_k}
 \ge\eta(\mathsf V_{\mathsf A}-C\delta\av r^4)
       -C\delta r_\delta^4.
\end{equation}
On $\mathcal P_k(s)$, we have $Q_k=h_\delta(s)$.
Using \eqref{wgt:cap-decay}, we therefore obtain
$-\partial_s Q_k/2\ge\kappa(1+\eta)r_\delta^2/8$.
Combining this bound with \eqref{wgt:exterior-damping} on
$\mathcal O_k(s)$ and $p_k\gtrsim1$ from \eqref{wgt:basicbounds},
we obtain, for sufficiently small $\nu$ and throughout
$\mathcal P_k(s)\cup\mathcal O_k(s)$,
\[
 \frac\eta4(\xi\cdot\nabla Q_k-|\nabla Q_k|^2)
       -\frac12\partial_s Q_k-\frac\eta2-\frac\eta{p_k}
 \ge\frac\kappa{32}(1+\eta)r_\delta^2.
\]
We choose a sufficiently large constant $C$ in \eqref{lin:diffusion-uncapped},
keep it fixed, and define
\begin{equation}\label{lin:diffusion-potential-definition}
 \mathsf V_k=
 \begin{cases}
 \mathsf V_{\mathsf A}-C\delta\av r^4,
       &\xi\in\mathcal I_k(s),\\
 \kappa r_\delta^2/64,
       &\xi\in\mathcal P_k(s)\cup\mathcal O_k(s).
 \end{cases}
\end{equation}
Using the preceding bounds, we obtain
\begin{equation}\label{lin:diffusion-global-lower-bound}
 \begin{aligned}
 &\eta\mathcal Q_{p_k}(f_k)
       -\frac12\int_{\R^2}(\partial_s p_k)f_k^2\dd\xi\\
 &\quad\ge\eta\int_{\R^2}
       \bigl\{|\nabla(\sqrt{p_k}f_k)|^2
                         +\mathsf V_kp_kf_k^2\bigr\}\dd\xi
       -C\delta r_\delta^4\|f_k\|_{L^2_{p_k}}^2\\
 &\qquad+\frac\kappa{64}(1+\eta)r_\delta^2
       \int_{\mathcal O_k(s)}p_kf_k^2\dd\xi.
 \end{aligned}
\end{equation}
Using \eqref{lin:Gaussian-potential}, the limit
$\delta r_\delta^2\to0$, and \eqref{wgt:interface-location}, we obtain
\begin{equation}\label{lin:diffusion-potential-bounds}
 \mathsf V_k\ge c\min\{r^2,r_\delta^2\}-C.
\end{equation}
For every fixed $L>0$, since $\delta\le\nu^{1-\gamma}\to0$ and
$r_\delta\to\infty$, we have $\{\xi:|\xi|\le L\}\subset\mathcal I_k(s)$
for sufficiently small $\nu$. On $\{|\xi|\le L\}$, we have
$p_k=\mathsf A_k^\delta$ and
$|p_k/\mathsf A-1|\lesssim\delta r^2\av r^2$. Hence we obtain
\begin{equation}\label{lin:diffusion-potential-limits}
 p_k\longrightarrow\mathsf A,\qquad
 \mathsf V_k\longrightarrow\mathsf V_{\mathsf A}
\end{equation}
uniformly for $|\xi|\le L$, $1\le k\le3$, and $0<s\le S_\nu$ as $\nu\to0$.

For $f_k\in L^2_{p_k}$ with zero mass and first moments and
$\sqrt{p_k}f_k\in H^1(\R^2)$, we claim that, for sufficiently small $\nu$,
\begin{equation}\label{lin:diffusion-global-coercivity}
 \int_{\R^2}\bigl\{|\nabla(\sqrt{p_k}f_k)|^2
                              +\mathsf V_kp_kf_k^2\bigr\}\dd\xi
 \ge c\left\{\|\nabla(\sqrt{p_k}f_k)\|_{L^2}^2
                         +\|f_k\|_{L^2_{p_k}}^2\right\}.
\end{equation}
We first establish a uniform lower bound of $c\|f_k\|_{L^2_{p_k}}^2$.
If this bound failed, we could fix $k$ along a subsequence and find
$\nu_n\to0$, $0<s_n\le S_{\nu_n}$, and functions
$u_{k,n}\in H^1(\R^2)$ such that
\[
 \|u_{k,n}\|_{L^2}=1,\qquad
 \int_{\R^2}\bigl(|\nabla u_{k,n}|^2
                    +\mathsf V_{k,n}u_{k,n}^2\bigr)\dd\xi\le\frac1n,
\]
with $\int_{\R^2}\frac{u_{k,n}}{\sqrt{p_{k,n}}}\dd\xi=0$ and
$\int_{\R^2}\frac{\xi u_{k,n}}{\sqrt{p_{k,n}}}\dd\xi=0$.
Here $p_{k,n}$ and $\mathsf V_{k,n}$ are the weight and potential
at $(\nu_n,s_n)$.
Writing $\delta_n=\nu_n(e^{s_n}-1)$ and using
\eqref{lin:diffusion-potential-bounds}, we obtain
\[
 \|\nabla u_{k,n}\|_{L^2}^2
       +\int_{\R^2}\min\{r^2,r_{\delta_n}^2\}u_{k,n}^2\dd\xi
 \le C.
\]
For every fixed $L\ge1$, we have $r_{\delta_n}\ge L$ for sufficiently large
$n$, and hence we obtain $\int_{r\ge L}u_{k,n}^2\dd\xi\le\frac C{L^2}$.
Using compactness in $L^2(\{|\xi|<L\})$ for each fixed $L\ge1$
and this tail bound, we pass to a subsequence for which
$u_{k,n}\rightharpoonup u_k$ in $H^1(\R^2)$,
$u_{k,n}\longrightarrow u_k$ in $L^2(\R^2)$, and
$\|u_k\|_{L^2}=1$.
Using the same bound, we have $ru_k\in L^2$.
Using $p_{k,n}\gtrsim\av r^{16}$ from \eqref{wgt:basicbounds}
and the local convergence above, we also obtain
$\left\|\av r\bigl(p_{k,n}^{-\frac{1}{2}}-\mathsf A^{-\frac{1}{2}}\bigr)\right\|_{L^2}
 \longrightarrow0$.
Passing to the limit in the two moment identities, we obtain
$\int_{\R^2}\frac{u_k}{\sqrt{\mathsf A}}\dd\xi=0$ and
$\int_{\R^2}\frac{\xi u_k}{\sqrt{\mathsf A}}\dd\xi=0$.
Moreover, we have
\[
 \int_{\R^2}\mathsf A
       \left|\nabla\left(\frac{u_k}{\sqrt{\mathsf A}}\right)\right|^2\dd\xi
 \le2\|\nabla u_k\|_{L^2}^2+C\|ru_k\|_{L^2}^2<\infty.
\]
By \eqref{lin:diffusion-potential-bounds}, we can fix $L>0$, independent of $n$,
such that $\mathsf V_{k,n}(\xi)\ge0$ for $|\xi|\ge L$ and sufficiently large $n$.
Using the local convergence in \eqref{lin:diffusion-potential-limits}
and the weak $H^1$ convergence, we obtain the lower
semicontinuity of the quadratic forms.
Applying the Gaussian diffusion estimate \cite{GS24,DoGa26}
to $u_k/\sqrt{\mathsf A}$, we then obtain
\[
 \begin{aligned}
 0<c\|u_k\|_{L^2}^2
 &\le\int_{\R^2}\bigl(|\nabla u_k|^2
                           +\mathsf V_{\mathsf A}u_k^2\bigr)\dd\xi\\
 &\le\liminf_{n\to\infty}\int_{\R^2}
       \bigl(|\nabla u_{k,n}|^2+\mathsf V_{k,n}u_{k,n}^2\bigr)\dd\xi
 \le0,
 \end{aligned}
\]
a contradiction. Using $\mathsf V_k\ge-C$ and this uniform
$L^2_{p_k}$ lower bound, we also control
$\|\nabla(\sqrt{p_k}f_k)\|_{L^2}^2$ and obtain
\eqref{lin:diffusion-global-coercivity}.
Substituting this coercivity estimate into \eqref{lin:diffusion-global-lower-bound} and
\eqref{lin:diffusionbound}, and then using
\eqref{wgt:energy-equivalence}, we obtain \eqref{energy:R1-bound},
since $\delta r_\delta^4=\delta^{1-4\kappa}$.
\end{proof}

\subsection{Transport and self-reaction}
\begin{lemma}
\label{lin:transport}\label{lin:selfreaction}
For every fixed $\varepsilon>0$ and sufficiently small $\nu$,
the term $R_2$ defined in \eqref{energy:R2} satisfies
\begin{equation}\label{energy:R2-bound}
 |R_2|\le{}C(\delta^{\frac12-3\kappa}+\nu)\mathcal{E} +\varepsilon(1+\eta)r_\delta^2 \sum_k\int_{\mathcal O_k(s)}p_kf_k^2\dd\xi.
\end{equation}
\end{lemma}
\begin{proof}
Since $B_k^\app$ is divergence free and
$v^{f_k}\cdot\nabla\psi^{f_k}=0$, we obtain the following identity by direct calculation:
\begin{equation}\label{energy:transport-reaction-identity}
 \begin{aligned}
 R_2={}&\frac12\sum_k\int_{\R^2}
             p_kf_k^2 B_k^\app\cdot\nabla Q_k\dd\xi\\
 &+\sum_k\int_{\R^2}f_k
       \left(B_k^\app+\frac{\Gamma_k}{\delta}
                    p_k\nabla^\perp w_{\app}^k\right)
                    \cdot\nabla\psi^{f_k}\dd\xi.
 \end{aligned}
\end{equation}
By \eqref{wgt:Gaussian-antisymmetry}, we have
$\int_{\R^2}f_k\,\xi^\perp\cdot\nabla\psi^{f_k}\dd\xi=0$,
so we may replace $B_k^\app$ by $B_k^\app+\beta\xi^\perp$
in the second integral.
On $r\le2r_\delta$, using \eqref{prof:shifted-Gaussian-expansion},
\eqref{lin:profileinputs}, and
$|Z'|+|Z'-\mathcal F(Z)|\le C_Z\delta^2$ in
\eqref{prof:center-velocity-bounds}, we obtain
\begin{equation}\label{lin:transportexpansion}
 \begin{aligned}
 B_k^\app
 &=\frac{\Gamma_k}{\delta}v^G+\Sigma_k\xi
       +\Gamma_kv^{F_k^0}-\beta\xi^\perp+\mathcal K_k,\\
 |\mathcal K_k|
 &\lesssim\frac{\nu+\delta^{\frac12}}{\av r}
                         +\delta^{\frac12}\av r^2.
 \end{aligned}
\end{equation}
On $\mathcal I_k(s)$, we have $p_k=\mathcal H'(G+\delta F_k^0)$.
Combining \eqref{lin:transportexpansion} with
\eqref{wgt:leading-stream} and \eqref{wgt:approximate-cancellation},
we write the two coefficients in
\eqref{energy:transport-reaction-identity} as
\begin{equation}\label{lin:transport-inner-identities}
 \begin{aligned}
 B_k^\app\cdot\nabla Q_k
 &=\frac{\Gamma_k}{\delta}\nabla^\perp\Theta_k\cdot\nabla Q_k
       +\mathcal K_k\cdot\nabla Q_k-\beta\partial_\theta Q_k,\\
 B_k^\app+\beta\xi^\perp
       +\frac{\Gamma_k}{\delta}p_k\nabla^\perp w_{\app}^k
 &=\frac{\Gamma_k}{\delta}\nabla^\perp\Theta_k+\mathcal K_k
       +\frac{\Gamma_k}{\delta}p_k
            \nabla^\perp(w_{\app}^k-G-\delta F_k^0).
 \end{aligned}
\end{equation}
Using \eqref{wgt:stream-defect} and $|\nabla Q_k|\lesssim r$, we obtain
\begin{equation}\label{lin:transport-defect}
 \frac1\delta\left|\int_{\mathcal I_k(s)}
       f_k^2\nabla^\perp\Theta_k\cdot\nabla p_k\dd\xi\right|
 \lesssim\delta r_\delta^6\|f_k\|_{L^2_{p_k}}^2.
\end{equation}
On $\mathcal I_k(s)$, we also have $\frac12\mathsf A\le p_k\le\frac32\mathsf A$.
Using also $\av\xi^5/\sqrt{\mathsf A}\in L^4(\R^2)$, we obtain
\[
 \begin{aligned}
 \frac1\delta\left|\int_{\mathcal I_k(s)}
       f_k\nabla^\perp\Theta_k\cdot\nabla\psi^{f_k}\dd\xi\right|
 &\le\frac1\delta\|f_k\|_{L^2_{p_k}}
       \|\nabla\psi^{f_k}\|_{L^4}
       \|\nabla\Theta_k/\sqrt{p_k}\|_{L^4(\mathcal I_k(s))}\\
 &\lesssim\delta\|f_k\|_{L^2_{p_k}}^2.
 \end{aligned}
\]
Here and below we use the velocity bound
$\|\nabla\psi^{f_k}\|_{L^4}\lesssim\|f_k\|_{L^2_{p_k}}$
from \eqref{wgt:stream-H2-bound}.

The two contributions of $\mathcal K_k$ satisfy
$\|\mathcal K_k\cdot\nabla Q_k\|_{L^\infty(\mathcal I_k(s))}
       \lesssim\nu+\delta^{\frac12}r_\delta^3$
and
\[
 \begin{aligned}
 \left|\int_{\mathcal I_k(s)}
           f_k\mathcal K_k\cdot\nabla\psi^{f_k}\dd\xi\right|
 &\le\|f_k\|_{L^2_{p_k}}\|\nabla\psi^{f_k}\|_{L^4}
       \|\mathcal K_k/\sqrt{p_k}\|_{L^4(\mathcal I_k(s))}\\
 &\lesssim(\nu+\delta^{\frac12})\|f_k\|_{L^2_{p_k}}^2.
 \end{aligned}
\]
By \eqref{lin:profileinputs} and \eqref{wgt:sharpgaussianupper},
the remaining reaction profile satisfies
$\|\sqrt{p_k}\nabla(w_{\app}^k-G-\delta F_k^0)\|_{L^4}
                   \lesssim\delta(\nu+\delta^{\frac12})$.
We therefore estimate its contribution as follows:
\[
 \begin{aligned}
 &\frac{|\Gamma_k|}{\delta}
 \left|\int_{\mathcal I_k(s)}p_kf_k
       \nabla^\perp(w_{\app}^k-G-\delta F_k^0)
                         \cdot\nabla\psi^{f_k}\dd\xi\right|\\
 &\qquad\lesssim\frac1\delta\|f_k\|_{L^2_{p_k}}
       \|\nabla\psi^{f_k}\|_{L^4}
       \|\sqrt{p_k}\nabla(w_{\app}^k-G-\delta F_k^0)\|_{L^4}
 \lesssim(\nu+\delta^{\frac12})\|f_k\|_{L^2_{p_k}}^2.
 \end{aligned}
\]
Finally, the weighted rotation satisfies
\[
 \left|\frac\beta2\int_{\mathcal I_k(s)}
                p_kf_k^2\partial_\theta Q_k\dd\xi\right|
 \lesssim\delta r_\delta^4\|f_k\|_{L^2_{p_k}}^2,
\]
by the angular bound in \eqref{wgt:Hessian-bounds}.
Using $r_\delta=\delta^{-\kappa}$, we bound the inner
contribution by
$C(\nu+\delta^{\frac12-3\kappa})\|f_k\|_{L^2_{p_k}}^2$.

On $\mathcal P_k(s)$, we have $\nabla Q_k=0$. Since
$\mathcal P_k(s)\cup\mathcal O_k(s)\subset\{r\ge r_\delta/2\}$,
using \eqref{lin:profileinputs} and \eqref{wgt:sharpgaussianupper}, we obtain
\[
 \frac1\delta
 \|\sqrt{p_k}\nabla w_\app^k\|_{L^4(\mathcal P_k(s)\cup\mathcal O_k(s))}
 \lesssim\delta^{-1}e^{-cr_\delta^2}\lesssim\delta^2.
\]
By \eqref{prof:corrected-transport}, \eqref{eq:Bapp}, and
\eqref{prof:center-law-definition}, we have
\begin{equation}\label{lin:transport-field}
 \begin{aligned}
 B_k^\app+\beta\xi^\perp
 ={}&\sum_j\frac{\Gamma_j}{\delta}
       v^{w_\app^j}\left(\xi+\frac{a_{kj}}{\sqrt\delta}\right)
       -\frac1{\sqrt\delta}\sum_{j\ne k}\Gamma_jK(a_{kj})\\
 &+\frac{\mathcal F_k(Z)-Z_k'}{\sqrt\delta}.
 \end{aligned}
\end{equation}
Using the uniform velocity bounds and
$|Z'-\mathcal F(Z)|\lesssim\delta^2$, we have
$|B_k^\app+\beta\xi^\perp|\lesssim\delta^{-1}$ on $\R^2$.
Thus, using \eqref{wgt:tailmoments}, we obtain
\begin{equation}\label{lin:transport-tails}
 \|(B_k^\app+\beta\xi^\perp)/\sqrt{p_k}\|_
      {L^4(\mathcal P_k(s)\cup\mathcal O_k(s))}
 \lesssim\delta^{-1}
       \|p_k^{-\frac12}\|_{L^4(r\ge r_\delta/2)}
 \lesssim\delta^2.
\end{equation}
Using these estimates together with \eqref{wgt:stream-H2-bound}, we bound
the second integral in \eqref{energy:transport-reaction-identity}
over $\mathcal P_k(s)\cup\mathcal O_k(s)$ by
$C\delta^2\|f_k\|_{L^2_{p_k}}^2$.
It remains to estimate the first integral in
\eqref{energy:transport-reaction-identity}.
For $r>R$, we have $|\xi+a_{kj}/\sqrt\delta|\ge r/2$.
Using \eqref{lin:transport-field} and the profile velocity bounds, we obtain
$|B_k^\app+\beta\xi^\perp|
 \lesssim(\delta r)^{-1}+\delta^{-\frac12}
 \lesssim\delta^{-\frac12}$.
Since $\nabla Q_k=\sqrt\nu\,\xi/r$ on $\mathcal O_k(s)$, we obtain
$|B_k^\app\cdot\nabla Q_k|
 =|(B_k^\app+\beta\xi^\perp)\cdot\nabla Q_k|
 \lesssim\sqrt\eta$.
Using $\delta=\nu/\eta$ and $r_\delta=\delta^{-\kappa}$, we have
$\frac{\sqrt\eta}{(1+\eta)r_\delta^2}
 =\nu^{2\kappa}\frac{\eta^{\frac12-2\kappa}}{1+\eta}
 \le\nu^{2\kappa}$.
For sufficiently small $\nu$, we therefore obtain
\[
 \frac12\left|\int_{\mathcal O_k(s)}
           p_kf_k^2 B_k^\app\cdot\nabla Q_k\dd\xi\right|
 \le\varepsilon(1+\eta)r_\delta^2
              \int_{\mathcal O_k(s)}p_kf_k^2\dd\xi.
\]

Summing the estimates on $\mathcal I_k(s)$, $\mathcal P_k(s)$, and
$\mathcal O_k(s)$ over $k$ and using
\eqref{wgt:energy-equivalence}, we obtain \eqref{energy:R2-bound}.
\end{proof}

\subsection{Reaction from the other vortices}
\begin{lemma}\label{lin:remotereaction}
For remainders satisfying \eqref{eq:error-zero-moments}, the term
$R_3$ defined in \eqref{energy:R3} satisfies
\begin{equation}\label{energy:R3-bound}
 |R_3|\lesssim\delta^{\frac12} \mathcal{E}.
\end{equation}
\end{lemma}
\begin{proof}
We fix $j\ne k$ and $\vartheta\in(\frac12,1)$. By
\eqref{lin:profileinputs} and \eqref{wgt:sharpgaussianupper}, we have
\[
 \sqrt{p_k}|\nabla w_\app^k|\lesssim e^{-(2\vartheta-1)r^2/8},
 \qquad
 \|\sqrt{p_k}\nabla w_\app^k\|_{L^2}
 +\|\sqrt{p_k}\nabla w_\app^k\|_{L^4}\lesssim1.
\]
On $|\xi|>|a_{kj}|/(4\sqrt\delta)$, using this Gaussian bound,
\eqref{eq:HLS-q}, and \eqref{wgt:functionmoments}, we obtain
\[
 \left\|\sqrt{p_k}|\nabla w_\app^k|
       \left|v^{f_j}\left(\cdot+\frac{a_{kj}}{\sqrt\delta}\right)\right|
       \right\|_{L^2(|\xi|>|a_{kj}|/(4\sqrt\delta))}
 \lesssim e^{-c/\delta}\|f_j\|_{L^2_{p_j}}
 \lesssim\delta^{\frac32}\|f_j\|_{L^2_{p_j}}.
\]

For $|\xi|\le|a_{kj}|/(4\sqrt\delta)$, we set
$x=\xi+a_{kj}/\sqrt\delta$. Using the zero moments
\eqref{eq:error-zero-moments} and Taylor's formula, we obtain
\[
 \begin{aligned}
 v^{f_j}(x)
 ={}&\underbrace{\int_{|y|\le|a_{kj}|/(4\sqrt\delta)}
       \int_0^1(1-t)D^2K(x-ty)[y,y]f_j(y)\dd t\dd y}_{\text{(I)}}\\
 &+\underbrace{\int_{|y|>|a_{kj}|/(4\sqrt\delta)}
       [-K(x)+DK(x)y]f_j(y)\dd y}_{\text{(II)}}\\
 &+\underbrace{\int_{|y|>|a_{kj}|/(4\sqrt\delta)}
       K(x-y)f_j(y)\dd y}_{\text{(III)}}.
 \end{aligned}
\]
For $|y|\le|a_{kj}|/(4\sqrt\delta)$ and $0\le t\le1$, we have
$|x-ty|\ge|a_{kj}|/(2\sqrt\delta)$. Hence, using
\eqref{prof:kernel-derivatives} and \eqref{wgt:functionmoments}, we obtain
\[
 |\text{(I)}|
 \lesssim\delta^{\frac32}
       \int_{|y|\le|a_{kj}|/(4\sqrt\delta)}|y|^2|f_j(y)|\dd y
 \lesssim\delta^{\frac32}\|f_j\|_{L^2_{p_j}}.
\]
Using $|x|\ge3|a_{kj}|/(4\sqrt\delta)$ and
\eqref{wgt:functiontails} with $N=3$, we obtain
\[
 \begin{aligned}
 |\text{(II)}|
 &\lesssim\sqrt\delta\,
       \|f_j\|_{L^1(|y|>|a_{kj}|/(4\sqrt\delta))}
       +\delta\|\av y f_j\|_{L^1(|y|>|a_{kj}|/(4\sqrt\delta))}\\
 &\lesssim\delta^{\frac72}\|f_j\|_{L^2_{p_j}}.
 \end{aligned}
\]
For the singular term $\text{(III)}$, using \eqref{eq:HLS-q} and
\eqref{wgt:functiontails}, we obtain
\[
 \|\text{(III)}\|_{L^4_\xi}
 \lesssim\|f_j\|_{L^{\frac43}(|y|>|a_{kj}|/(4\sqrt\delta))}
 \lesssim\delta^3\|f_j\|_{L^2_{p_j}}.
\]
Multiplying the first two bounds by
$\|\sqrt{p_k}\nabla w_\app^k\|_{L^2}$ and the last by
$\|\sqrt{p_k}\nabla w_\app^k\|_{L^4}$, we combine these bounds with the
estimate on $\{|\xi|>|a_{kj}|/(4\sqrt\delta)\}$ to obtain
\begin{equation}\label{lin:remotereactionbound}
 \|\sqrt{p_k}|\nabla w_\app^k|
                  |v^{f_j}(\cdot+\frac{a_{kj}}{\delta^{\frac12}})|\|_{L^2}
             \lesssim\delta^{\frac32}\|f_j\|_{L^2_{p_j}}.
\end{equation}
Finally, using \eqref{wgt:basicbounds} and \eqref{wgt:stream-H2-bound},
we have $\|\psi^{f_k}/\sqrt{p_k}\|_{L^2}\lesssim\|f_k\|_{L^2_{p_k}}$.
Using \eqref{lin:remotereactionbound} in \eqref{energy:R3}, we obtain
\[
 |R_3|\lesssim\delta^{\frac12}\sum_{j\ne k}
       \|f_k\|_{L^2_{p_k}}\|f_j\|_{L^2_{p_j}}
 \lesssim\delta^{\frac12}\mathcal E.
\]
\end{proof}
\subsection{Quadratic transport}

\begin{lemma}\label{lem:nonlinear-energy}
For every $2<q<\infty$ and every $\varepsilon>0$,
the term $R_4$ defined in \eqref{energy:R4} satisfies
\begin{equation}\label{eq:nonlinear-energy}
 |R_4|\le\varepsilon\eta\mathcal D
       +C_{q,\varepsilon}\eta^{-\frac{1}{q-1}}
              \mathcal E^{\frac{3q-2}{2(q-1)}}+C\mathcal E^{\frac32}.
\end{equation}
\end{lemma}
\begin{proof}
The translated velocities are divergence free. By a direct calculation, we obtain
\[
 \begin{aligned}
 R_4={}&\underbrace{\frac12\sum_{j,k}\Gamma_j\int_{\R^2}
       p_kf_k^2v^{f_j}\!\left(\xi+\frac{a_{kj}}{\sqrt\delta}\right)
       \cdot\nabla Q_k\dd\xi}_{\text{(I)}}\\
 &+\underbrace{\sum_{j,k}\Gamma_j\int_{\R^2}
       f_kv^{f_j}\!\left(\xi+\frac{a_{kj}}{\sqrt\delta}\right)
       \cdot\nabla\psi^{f_k}\dd\xi}_{\text{(II)}}.
 \end{aligned}
\]
By \eqref{eq:nonlinear-velocity}, \eqref{wgt:energy-equivalence},
and interpolation in $\R^2$, we obtain
\[
 \begin{aligned}
 |\text{(I)}|
 &\le C_q\mathcal E^{\frac12}
       \sum_k\|\sqrt{p_k}f_k\|_{L^{\frac{2q}{q-1}}}^2\\
 &\le C_q\mathcal E^{\frac12}
       \sum_k\|f_k\|_{L^2_{p_k}}^{2-\frac2q}
       \|\nabla(\sqrt{p_k}f_k)\|_{L^2}^{\frac2q}
 \le C_q\mathcal E^{\frac32-\frac1q}\mathcal D^{\frac1q}.
 \end{aligned}
\]
The terms $j=k$ in $\text{(II)}$ vanish because
$v^{f_k}\cdot\nabla\psi^{f_k}=0$. Using \eqref{eq:HLS-q},
\eqref{wgt:functionmoments}, and $p_k\gtrsim1$, we obtain
\[
 |\text{(II)}|
 \lesssim\sum_{j\ne k}\|f_k\|_{L^2}
       \|v^{f_j}\|_{L^4}\|\nabla\psi^{f_k}\|_{L^4}
 \lesssim\sum_{j\ne k}\|f_k\|_{L^2_{p_k}}^2\|f_j\|_{L^2_{p_j}}
 \lesssim\mathcal E^{\frac32}.
\]
Hence, for every $\varepsilon>0$, we obtain
\[
 |R_4|\le C_q\mathcal E^{\frac32-\frac1q}\mathcal D^{\frac1q}
       +C\mathcal E^{\frac32}
 \le\varepsilon\eta\mathcal D
       +C_{q,\varepsilon}\eta^{-\frac1{q-1}}
              \mathcal E^{\frac{3q-2}{2(q-1)}}
       +C\mathcal E^{\frac32}.
\]
\end{proof}

\subsection{The error term}
\begin{lemma}\label{lem:residual-energy}
There is a fixed constant $C_R$ such that the term $R_5$ defined in
\eqref{energy:R5} satisfies
\begin{equation}\label{energy:R5-bound}
 |R_5|\le\tfrac14(1+\eta)\mathcal{E}+C\delta^{\frac12} \mathcal{E}
       +C\delta^{\frac52}\mathcal E^{\frac32}
                              +C_R\delta(\delta+\nu).
\end{equation}
\end{lemma}
\begin{proof}
Using \eqref{wgt:basicbounds}, \eqref{wgt:stream-H2-bound},
\eqref{wgt:energy-equivalence}, \eqref{lin:profileinputs}, and
\eqref{wgt:sharpgaussianupper}, we obtain
\[
 \sum_k\left\|\sqrt{p_k}f_k+
          \frac{\psi^{f_k}}{\sqrt{p_k}}\right\|_{L^2}^2
       \lesssim\mathcal E,\qquad
 \|\sqrt{p_k}\nabla w_\app^k\|_{L^2}\lesssim1.
\]
Thus, using \eqref{prof:residualweighted} and
\eqref{ctr:refined-center-equation}, we obtain
\[
 \begin{aligned}
 |R_5|
 &\lesssim(\delta+\nu)\sqrt{\mathcal E}
              +\delta^{-\frac32}|e|\sqrt{\mathcal E}\\
 &\lesssim(\delta+\nu)\sqrt{\mathcal E}
       +\delta^{\frac12}\mathcal E
       +\delta^{\frac52}\mathcal E^{\frac32},
 \end{aligned}
\]
where we used $\delta\le1$. Since $\eta=\nu/\delta$, we obtain
\[
 C(\delta+\nu)\sqrt{\mathcal E}
 \le\frac14(1+\eta)\mathcal E
       +C_R\frac{(\delta+\nu)^2}{1+\eta}
 =\frac14(1+\eta)\mathcal E+C_R\delta(\delta+\nu),
\]
and we conclude \eqref{energy:R5-bound}.
\end{proof}

\subsection{Proof of the energy estimate}
\begin{proof}[Proof of Proposition~\ref{prop:energy}]
We choose $\varepsilon=c/4$ in \eqref{energy:R2-bound} and
\eqref{eq:nonlinear-energy}, with $c$ from \eqref{energy:R1-bound}.
Combining these bounds with \eqref{energy:R3-bound},
\eqref{energy:R5-bound}, and \eqref{eq:energy-five-terms}, we obtain
\[
 \begin{aligned}
 &\mathcal E'+\frac{3c}{4}\eta\mathcal D+c\eta\mathcal E
       +\frac74(1+\eta)\mathcal E
       +\frac{3c}{4}(1+\eta)r_\delta^2
          \sum_k\int_{\mathcal O_k(s)}p_kf_k^2\dd\xi\\
 &\quad\le C\bigl(\delta^{1-4\kappa}
       +\delta^{\frac12-3\kappa}+\delta^{\frac12}+\nu\bigr)\mathcal E
       +C\delta(\delta+\nu)\\
 &\qquad+C_q\eta^{-\frac1{q-1}}
       \mathcal E^{\frac{3q-2}{2(q-1)}}
       +C(1+\delta^{\frac52})\mathcal E^{\frac32}.
 \end{aligned}
\]
Since $\delta\le\nu^{1-\gamma}\le1$, we have
\[
 C\bigl(\delta^{1-4\kappa}+\delta^{\frac12-3\kappa}
       +\delta^{\frac12}+\nu\bigr)
 \le C\bigl(\nu^{(1-\gamma)(\frac12-3\kappa)}+\nu\bigr)
 \le\frac12
\]
for sufficiently small $\nu$. Absorbing this term and using
$1+\delta^{\frac52}\le2$, we obtain \eqref{eq:complete-energy}
after decreasing $c$. The choice of $\varepsilon$ and the
required smallness of $\nu$ are independent of $q$.
\end{proof}

\section{Proof of the main theorems}
\label{sec:proof-main}

\begin{proposition}\label{prop:coupled-continuation}
For sufficiently small $\nu>0$, the solution satisfies
\begin{equation}\label{eq:energy-growth}
 \mathcal E(s)\lsim{\gamma}\delta(s)^2,\qquad
 |Z(s)-z^0|\lsim{\gamma}\delta(s)^2,
 \qquad 0<s\le S_\nu,
\end{equation}
with implicit constants independent of $s$ and $\nu$.
\end{proposition}

\begin{proof}
We set $s_0=\log2$, corresponding to $t_0=\frac12$, and take $\nu$
small enough that $S_\nu>s_0$. Since the reference centers remain
separated on $[0,t_0]$, we apply the finite-time weighted expansion
in \cite{Ga11} and get
\begin{equation}\label{match:coarse-conclusion}
 |Z(s)-z^0|\lesssim\delta(s)^2,\qquad
 0\le\mathcal E(s)\lesssim\delta(s),\qquad 0<s\le s_0.
\end{equation}
Using \eqref{ctr:closed-force} and $\mathcal F(z^0)=0$, we also obtain
$|Z'|+|Z'-\mathcal F(Z)|\lesssim\delta^2$.
We apply Proposition~\ref{prop:energy} with $q=4$.
Since $\mathcal E\lesssim\delta\le\nu$ and $\eta\ge1$, we bound
its nonlinear terms by $C(\nu^{\frac23}+\nu^{\frac12})\mathcal E$.
Absorbing these terms for small $\nu$, we obtain
\[
 \mathcal E'+\frac12\mathcal E
 \le C\delta(\delta+\nu)\le C\nu^2s.
\]
Integrating from $\varepsilon$ to $s\le s_0$ and letting
$\varepsilon\to0$, we use \eqref{match:coarse-conclusion}
to obtain $\mathcal E(s)\lesssim\nu^2s^2\le\delta(s)^2$.
In particular, since $\delta(s_0)=\nu$, we have
\begin{equation}\label{eq:matching-use}
 \mathcal E(s_0)\lesssim\nu^2,\qquad
 |Z(s_0)-z^0|\lesssim\nu^2.
\end{equation}

For fixed $\nu>0$, using positive-time regularity and the Gaussian
bounds for the exact and approximate profiles, we obtain continuity
of $Z$ and $\mathcal E$ while the centers remain separated.
For large $M,K>1$, let $[s_0,S)$, with $S\le S_\nu$, be the maximal
interval on which
\begin{equation}\label{eq:nonlinear-bootstrap}
 \mathcal E(s)\le K\delta(s)^2,\qquad
 |Z(s)-z^0|\le M\delta(s)^2.
\end{equation}
By \eqref{eq:matching-use}, we have $S>s_0$.
Using $\nu\le\delta(s)\le\nu^{1-\gamma}$ and taking $\nu$ small
for fixed $M,K$, we have $|Z-z^0|\le\rho/2$, where $\rho$ is chosen
in Proposition~\ref{ctr:continuation}. In particular, the gap
\eqref{prof:gap} and diameter bound \eqref{wgt:diameter} hold.
Using \eqref{ctr:closed-force}, we obtain
\[
 |\mathfrak r(s)|
 \le C\delta(s)^2\bigl(1+\sqrt K\nu^{1-\gamma}
                         +K\nu^{4(1-\gamma)}\bigr)
 \le C\delta(s)^2,
\]
with $C$ independent of $M,K$.
By Proposition~\ref{ctr:continuation}, we then have
$|Z-z^0|\le C\delta^2+CM^2\delta^4$.
Fixing $M\ge4C$ and decreasing $\nu$, we obtain
$|Z-z^0|\le\frac M2\delta^2$ and, using $\mathcal F(z^0)=0$,
$|Z'|+|Z'-\mathcal F(Z)|\lsim{M}\delta^2$.
Thus we can apply Proposition~\ref{prop:energy} with constants
depending on $M$, but not on $K$.

Fix $q>\max\{4,\gamma/(1-\gamma)\}$.
Using \eqref{eq:nonlinear-bootstrap} and $\eta=\nu/\delta$, we obtain
\begin{align}
 \eta^{-\frac{1}{q-1}}\mathcal E^{\frac q{2(q-1)}}
 &\le K^{\frac q{2(q-1)}}
       \left(\frac{\delta^{q+1}}\nu\right)^{\frac{1}{q-1}}
 \le K^{\frac q{2(q-1)}}
       \nu^{\frac{q(1-\gamma)-\gamma}{q-1}},
       \label{eq:ratio-bound}\\*
 \sqrt{\mathcal E}&\le\sqrt K\nu^{1-\gamma}.
       \label{eq:energy-small}
\end{align}
Taking $\nu$ sufficiently small, we absorb the nonlinear terms in
\eqref{eq:complete-energy} and use $\nu\le\delta$ to obtain
\begin{equation}\label{eq:absorbed-energy}
 \mathcal E'+c\eta(\mathcal D+\mathcal E)
       +\frac12\mathcal E\le C\delta^2.
\end{equation}
Using $\delta(s')\le e^{-(s-s')}\delta(s)$ for $s_0\le s'\le s$
and \eqref{eq:matching-use}, we integrate to obtain
\[
 \mathcal E(s)
 \le e^{-(s-s_0)/2}\mathcal E(s_0)
       +C\int_{s_0}^s e^{-(s-s')/2}\delta(s')^2\dd s'
 \le C\delta(s)^2,
\]
with $C$ independent of $K$. Fixing $K\ge4C$ and decreasing $\nu$
if necessary, we strictly improve both bounds in
\eqref{eq:nonlinear-bootstrap}. By continuity, we conclude that
$S=S_\nu$ and obtain \eqref{eq:energy-growth}.
\end{proof}

\begin{proof}[Proof of Theorem~\ref{thm:main}]
Using \eqref{prof:approx} and \eqref{prof:derivativebounds}, we obtain
$\sum_k\|w_{\app}^k-G\|_{L^1}\lesssim\delta$.
By \eqref{wgt:inversemoments}, \eqref{wgt:energy-equivalence}, and
\eqref{eq:energy-growth}, we have
\begin{equation}\label{eq:vortex-conclusion}
 \begin{aligned}
 \sum_k\|w_k-G\|_{L^1}
 &\le\delta\sum_k\left(\int_{\R^2}p_k^{-1}\dd\xi\right)^{\frac12}
       \|f_k\|_{L^2_{p_k}}+\sum_k\|w_{\app}^k-G\|_{L^1}\\
 &\lesssim\delta\sqrt{\mathcal E}+\delta\lsim{\gamma}\delta,
 \qquad 0<s\le S_\nu.
 \end{aligned}
\end{equation}
Applying Proposition~\ref{lem:reduction}, we conclude the proof.
\end{proof}

\begin{proof}[Proof of Theorem~\ref{thm:centers}]
Using \eqref{eq:physical-centroids}, \eqref{eq:reference-centers}, and
\eqref{eq:energy-growth}, we obtain, for $s=-\log(1-t)$,
\[
 \begin{aligned}
 \max_{1\le k\le3}
 \left|\frac1{\Gamma_k}\int_{\R^2}x\omega_k^\nu(t,x)\dd x-z_k(t)\right|
 &=\sqrt{1-t}\,\max_{1\le k\le3}|Z_k(s)-z_k^0|\\
 &\lsim{\gamma}\sqrt{1-t}\,\delta(s)^2
 =\frac{\nu^2t^2}{(1-t)^{\frac32}},
 \qquad 0<t\le T_\nu.
 \end{aligned}
\]
Since $T_\nu=1-\nu^\gamma$, we also have
$\sup_{0<t\le T_\nu}\nu^2t^2/(1-t)^{\frac32}
\le\nu^2/(\nu^\gamma)^{\frac32}=\nu^{2-\frac32\gamma}$,
which proves \eqref{eq:main-centers-uniform}.
\end{proof}

\subsection*{Acknowledgments}
The author is supported by the NSFC Excellent Young Scientists
Fund (Overseas) and the Research Start-up Fund of Shanghai Jiao Tong
University.

\subsection*{Data availability}
No datasets were generated or analysed during the current study.

\subsection*{Competing interests}
The author declares no competing interests.

\end{document}